\documentclass{m2an}

\newcommand\Hline{%
\noalign{\xdef\origarrayrulewidth{\the\arrayrulewidth}%
\global\arrayrulewidth 2\arrayrulewidth}%
\hline%
\noalign{\global\arrayrulewidth\origarrayrulewidth}%
}

\usepackage{graphicx}        
\usepackage{url}             
\usepackage{psfrag}          
\usepackage{wrapfig}         
\usepackage[hang,centerlast]{subfigure}

\newcommand{\RR}{\mathbb{R}}

\newcommand{\NN}{\mathbb{N}}

\renewcommand{\t}{\tilde}

\newcommand{\e}{\textbf{e}}
\newcommand{\M}{\mathcal{M}}
\newcommand{\A}{\mathcal{A}}
\newcommand{\B}{\mathcal{B}}

\renewcommand{\O}{\mathcal{O}}
\newcommand{\N}{\mathcal{N}}
\newcommand{\id}[1]{_{\vert_{#1}}}

\theoremstyle{definition}
\newtheorem{assmptm}[thrm]{Assumption}

\begin{document}
\title{Boussinesq/Boussinesq systems for internal waves with a free surface, and the KdV approximation}

\author{Vincent Duch\^ene}
\address{D\'epartement de Math\'ematiques et Applications, UMR 8553, \'Ecole normale sup\'erieure, 45 rue d'Ulm, F 75230
Paris cedex 05, France ; e-mail: duchene@dma.ens.fr}
\date{\today}
\begin{abstract}
We study here some asymptotic models for the propagation of internal and surface waves in a two-fluid system. We focus on the so-called long wave regime for one-dimensional waves, and consider the case of a flat bottom. Following the method presented in~\cite{BonaColinLannes05} for the one-layer case, we introduce a new family of symmetric hyperbolic models, that are equivalent to the classical Boussinesq/Boussinesq system displayed in~\cite{ChoiCamassa96}. We study the well-posedness of such systems, and the asymptotic convergence of their solutions towards solutions of the full Euler system. Then, we provide a rigorous justification of the so-called KdV approximation, stating that any bounded solution of the full Euler system can be decomposed into four propagating waves, each of them being well approximated by the solutions of uncoupled Korteweg-de Vries equations. Our method also applies for models with the rigid lid assumption, using the Boussinesq/Boussinesq models introduced in~\cite{BonaLannesSaut08}. Our explicit and simultaneous decomposition allows to study in details the behavior of the flow depending on the depth and density ratios, for both the rigid lid and free surface configurations. In particular, we consider the influence of the rigid lid assumption on the evolution of the interface, and specify its domain of validity. Finally, solutions of the Boussinesq/Boussinesq systems and the KdV approximation are numerically computed, using a Crank-Nicholson scheme with a predictive step inspired from~\cite{Besse98,BesseBruneau98}.
\end{abstract}

\maketitle

\section{Introduction}
\label{Sec:Introduction}
\subsection{Motivation of the problem}
\begin{figure} \centering
\psfrag{z1}[Bc][Br]{\begin{footnotesize}$\zeta_1(t,x)$                          \end{footnotesize}}
\psfrag{z2}[Bc][Br]{\begin{footnotesize}$\zeta_2(t,x)$                          \end{footnotesize}}
\psfrag{a1}{\begin{footnotesize}$a_2$                          \end{footnotesize}}
\psfrag{a2}{\begin{footnotesize}$a_1$                          \end{footnotesize}}
\psfrag{h1}[][]{\begin{footnotesize}$d_1$                          \end{footnotesize}}
\psfrag{0}[][]{\begin{footnotesize}$0$                          \end{footnotesize}}
\psfrag{h2}[][]{\begin{footnotesize}$-d_2\ \ $                          \end{footnotesize}}
\psfrag{n1}{\begin{footnotesize}$n_1$                          \end{footnotesize}}
\psfrag{n2}{\begin{footnotesize}$n_2$                          \end{footnotesize}}
\psfrag{z}{\begin{footnotesize}$z$                   \end{footnotesize}}
\psfrag{X}{\begin{footnotesize}$x$                   \end{footnotesize}}
\psfrag{O1}{\begin{footnotesize}$\Omega_1^t$                   \end{footnotesize}}
\psfrag{O2}{\begin{footnotesize}$\Omega_2^t$                   \end{footnotesize}}
 \includegraphics[width=0.8\textwidth,keepaspectratio=true]{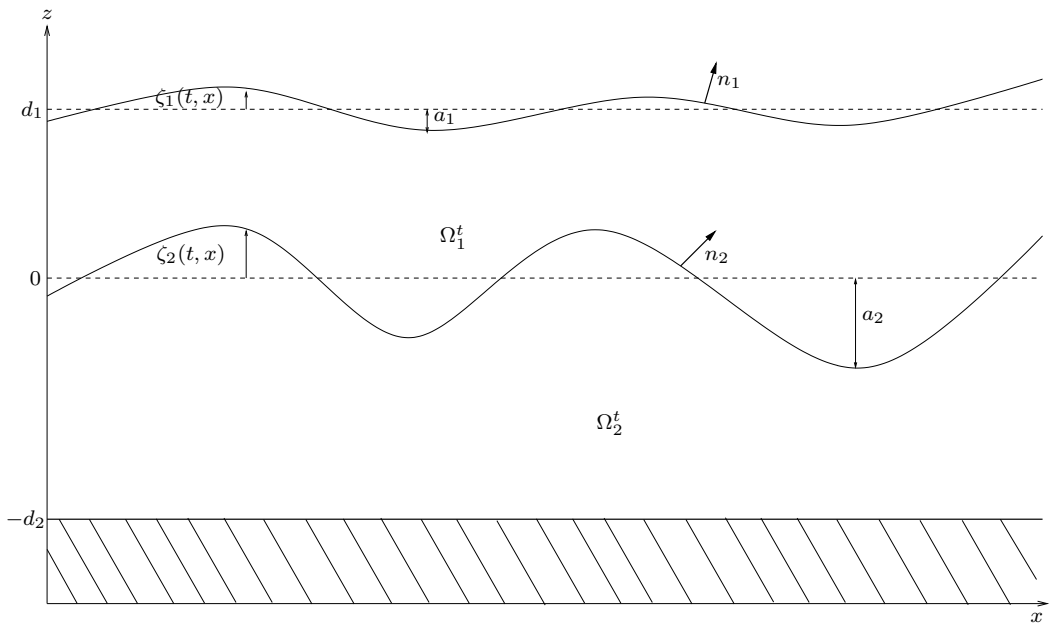}
 \caption{Sketch of the domain.}\label{FigDomain}
\end{figure}
This paper deals with different asymptotic models for the propagation of weakly nonlinear internal waves in a two-fluid system. The system we study consists in two layers of immiscible, homogeneous, ideal, incompressible and irrotationnal fluids under the only influence of gravity. Since we are interested in KdV equations, which are unidirectional, we focus on the one-dimensional case, and the bottom is assumed to be flat (see Figure~\ref{FigDomain}). 

Let us denote by $\rho_1$ the density of the upper fluid and $d_1$ its depth, $\rho_2$ the density of the lower fluid, and $d_2$ its depth, $a_1$ the typical amplitude of the deformation of the surface and $a_2$ the one of the interface, and finally $\lambda$ is a characteristic horizontal length, say the typical wavelength of the interface. Then the regime of the system is described by the following dimensionless parameters:
\[ \epsilon_1=\dfrac{a_1}{d_1}\in(0,1), \quad \epsilon_2 =\dfrac{a_2}{d_1}\in(0,1), \quad \mu =\dfrac{d_1^2}{\lambda^2}\in(0,+\infty), \quad \delta=\dfrac{d_1}{d_2}\in(0,+\infty), \quad  \gamma=\dfrac{\rho_1}{\rho_2} \in(0,1).
\]
The governing equations of such a system, that we call ``full Euler'', have been obtained in~\cite{Duchene10}; we briefly recall the system in Section~\ref{sec:fullEuler} below. 
This system is strongly nonlinear, and remains complex, for its direct study as well as for numerical computations. That is why the derivation of approximate asymptotic models has attracted lots of interests in the past decades. We focus here on the so-called {\em long wave regime}, where the dimensionless parameters 
$\epsilon_1$, $\epsilon_2$, and $\mu$ are small and of the same order:
\[\epsilon_1\ \sim\ \epsilon_2\ \sim\ \mu \ \ll\ 1.\]

The long wave regime for the one-fluid system (water wave) has been considerably studied, and has led to many approximate equations. Among them, of particular interest are the Boussinesq systems (from the work of Boussinesq~\cite{Boussinesq71,Boussinesq72}), and the KdV approximation (from Korteweg and de Vries~\cite{KortewegDe95}). The latter model states that any solution of the one-layer water wave problem in the long-wave limit splits up into two counter-propagating waves, each of them evolving independently as a solution of a KdV equation. A justification of such models has been investigated among others by Craig~\cite{Craig85}, Schneider-Wayne~\cite{SchneiderWayne00} and Bona-Colin-Lannes~\cite{BonaColinLannes05}. The study of internal waves followed quickly. When the surface is assumed to be fixed as a rigid lid, the KdV equations related to such a system have been formally obtained by Keulegan~\cite{Keulegan53} and Long~\cite{Long56}. The related Boussinesq-type systems have been justified (among many other asymptotic models) by Bona-Lannes-Saut~\cite{BonaLannesSaut08}. When the surface is not rigid and allowed to move as a free surface, it is known that there exist two different modes of wave motion, corresponding to different linear phase speeds (see Kakutani-Yamasaki~\cite{KakutaniYamasaki78}, Leone-Segur-Hammack~\cite{LeoneSegurHammack82}, Michallet-Barth\'elemy~\cite{MichalletBarth'elemy98} and Craig-Guyenne-Kalisch~\cite{CraigGuyenneKalisch05} for example). Accordingly, the KdV approximation states that any deformation of the surface and/or the interface will split up into four waves, each of them being lead by KdV equations. 

However, to our knowledge, the four KdV equations related to the problem have never been revealed {\it simultaneously}, neither rigorously justified with a convergence theorem. Such a precise and complete decomposition allows to directly compare the different models. In particular, as an application of our results, we present an in-depth study of the influence of the rigid lid assumption on the evolution of the interface, and therefore assert the domain of validity of such a hypothesis. 

An intermediate result for the construction of the KdV approximation is the full justification (in the sense of well-posedness, and convergence of the solutions towards solutions of the full Euler system) of symmetric coupled models, that are equivalent to the Boussinesq-type models derived by Camassa-Choi~\cite{ChoiCamassa96}, Craig-Guyenne-Kalisch~\cite{CraigGuyenneKalisch05} and Duch\^ene~\cite{Duchene10}. The latter systems are justified only by a consistency result, and the stronger properties of the symmetric models make them interesting by their own.

The construction and rigorous justification of the symmetric Boussinesq-type models and the KdV approximation are the main motivation of this article.

\subsection{Main results and outline of the paper}\label{sec:Outline}
Our study applies the methods of~\cite{BonaColinLannes05} to the case of internal waves, and accordingly uses as direct supports the full Euler system and the Boussinesq/Boussinesq model obtained in~\cite{Duchene10}. The derivation of the full Euler system as governing equations of our problem is quickly explained in Section~\ref{sec:fullEuler}. The full Euler system is consistent with the Boussinesq/Boussinesq model at order $\O(\mu^2)$, provided that $\epsilon_1,\ \epsilon_2  =  \O(\mu)$; this result is recalled and precisely stated in Proposition~\ref{Prop:ConsistencyBouss}, below.

\medskip

The first step of our study lies in the construction of symmetric systems, obtained from the original Boussinesq/Boussinesq model by using a first order symmetrizer, and withdrawing the $\O(\epsilon^2)$ terms. Section~\ref{Sec:DerivationAndAnalysisOfANewSymmetricModel} is dedicated to the construction, analysis and justification of such models. The systems we obtain are of the form 
\[\Big(S_0+\epsilon \big(S_1(U)-S_2\partial_x^2\big)\Big)\partial_t U+\Big(\Sigma_0+\epsilon\big(\Sigma_1(U)- \Sigma_2\partial_x^2\big)\Big)\partial_x U=0,\]
with $U(t,x)\in\RR^4$, symmetric matrices $S_0$, $\Sigma_0$, $S_2$, $\Sigma_2 \in \M_4(\RR)$ ($S_0$ and $S_2$ being definite positive), and linear mappings $S_1(\cdot)$ and $\Sigma_1(\cdot)$ with values in real symmetric matrices. The well-posedness of such systems over times of order $\O(1/\epsilon)$ is stated in Proposition~\ref{Prop:WPSBOUSS}. The convergence of their solutions towards solutions of the full Euler system is then proved to be of order $\O(\epsilon)$, for times of order $\O(1/\epsilon)$, in Proposition~\ref{Prop:ConvBoussEuler}.  

\medskip

From these models, using the classical WKB method, we prove in Section~\ref{Sec:TheKdVApproximation} that a rougher approximation consists in \emph{four uncoupled KdV equations}, that is to say that any bounded solution of the full Euler system in the long-wave limit splits up into four propagating waves, two of them moving to the right at different speed, and the other two moving to the left, each of them approximated by independent solutions of KdV equations. Our main result is Theorem~\ref{Th:CVKdVEuler}, which states explicitly the decomposition, and precise the convergence rate between bounded solutions of the full Euler system and the solutions of the KdV approximation. The case of the rigid-lid configuration can be treated in the same way, and is quickly tackled in Section~\ref{sSec:RigidLid}. 

The complete, simultaneous decomposition of the flow, and its rigorous justification, is a key point for the comparison with different models. In particular, we present in Section~\ref{sSec:AnalysisCoefs} an in-depth study of the behavior of the flow predicted by the KdV approximation, depending of the density ratio $\gamma$ and the depth ratio $\delta$, for both the rigid lid and free surface configurations. As a result, we show that the rigid lid hypothesis is valid only for small density differences between the two fluids. To our knowledge, this fact has never been established, thought it has been addressed for example in~\cite{Keulegan53,MichalletBarth'elemy98}.

\medskip

Finally, the Boussinesq/Boussinesq models and the KdV approximation are numerically computed and compared for different initial data and parameters in Section~\ref{Sec:NumericalComparison}. The numerical schemes we use are based on Crank-Nicholson methods, with a predictive step in order to deal with the nonlinearities. Such a scheme has been introduced by Besse and Bruneau in~\cite{Besse98,BesseBruneau98}; it is formally of order 2 in space and time, and appears to be unconditionally stable. The precise schemes in our framework are presented in Section~\ref{sSec:NumericalSchemes}.

\subsection{The full Euler system}\label{sec:fullEuler}
Let us recall here briefly the governing equations of our system (see~\cite{Duchene10} for more details). The velocity fields of an irrotational flow can be expressed as gradients of potentials (that we call $\phi_1$ for the upper fluid and $\phi_2$ for the lower fluid). The conservation of mass, together with the incompressibility of the fluids, yields Laplace equation for the potentials. The momentum equations of the Euler equations can then be integrated, which yields the Bernoulli equation in terms of potentials. The system is finally closed by kinematic boundary conditions (stating that no particles of fluid cross the bottom, the surface or the interface) and pressure laws (the pressure is assumed to be constant at the surface, and continuous at the interface, see Remark~\ref{rem:SurfTension} below). Thanks to an appropriate scaling, the two-layer full Euler system can be written in dimensionless form. 

\medskip

The key point is then to remark that the system can be reduced into four evolution equations coupling Zakharov's canonical variables~\cite{Zakharov68}, specifically the deviation of the free surface and interface from their rest position (respectively $\zeta_1$, $\zeta_2$), and the trace of the velocity potentials of the upper fluid evaluated at the surface (namely $\psi_1$) and of the lower fluid evaluated at the interface (namely $\psi_2$). The system is then given by
\begin{equation}
\label{FullEuler}
\left\{ \begin{array}{l}
\displaystyle\alpha\partial_{ t}{\zeta_1}-\frac{1}{\mu}G_1(\psi_1,\psi_2)=0,  \\
\displaystyle\partial_{ t}{\zeta_2}-\frac{1}{\mu}G_2\psi_2=0,  \\
\displaystyle\partial_{ t} \partial_x{\psi_1}+\alpha\partial_x{\zeta_1}+\frac{\epsilon_2}{2} \partial_x(|\partial_x {\psi_1}|^2)-\mu\epsilon_2\partial_x\N_1=0, \\
\displaystyle\partial_{ t} (\partial_x{\psi_2}-\gamma H(\psi_1,\psi_2))+(1-\gamma)\partial_x{\zeta_2} +\frac{\epsilon_2}{2} \partial_x(|\partial_x{\psi_2}|^2-\gamma |H(\psi_1,\psi_2)|^2) -\mu\epsilon_2\partial_x\N_2=0,
\end{array} 
\right. 
\end{equation}
with $\alpha=\frac{\epsilon_1}{\epsilon_2}$ and where $\N_1$ and $\N_2$ are given by the formulae
\[\N_1  \equiv \dfrac{(\frac{1}{\mu}G_1(\psi_1,\psi_2)+\epsilon_1\partial_x{\zeta_1} \partial_x{\psi_1})^2}{2(1+\mu|\epsilon_1\partial_x{\zeta_1}|^2)},\quad
       \N_2 \equiv \dfrac{(\frac{1}{\mu}G_2\psi_2+\epsilon_2\partial_x{\zeta_2} \partial_x{\psi_2})^2-\gamma(\frac{1}{\mu}G_2\psi_2+\epsilon_2\partial_x{\zeta_2}  H(\psi_1,\psi_2))^2}{2(1+\mu|\epsilon_2\partial_x{\zeta_2}|^2)},\]
and with $G_1$ and $G_2$ the Dirichlet-Neumann operators and $H$ the interface operator, defined as follows.

\begin{dfntn}
For $\zeta_1$, $\zeta_2 \in W^{2,\infty}(\RR)$ and $\partial_x \psi_1$, $\partial_x \psi_2\in H^{1/2}(\RR)$, the operators $G_1$, $G_2$ and $H$ are defined by
\begin{align*}
G_2[\epsilon_2\zeta_2]\psi_2&\equiv(\partial_z\phi_2-\mu\epsilon_2\partial_x\zeta_2  \partial_x\phi_2)\id{z=\epsilon_2\zeta_2}.\\
G_1[\epsilon_1\zeta_1,\epsilon_2\zeta_2](\psi_1,\psi_2)&\equiv(\partial_z\phi_1-\mu\epsilon_1\partial_x\zeta_1  \partial_x\phi_1)\id{z=1+\epsilon_1\zeta_1},\\
H[\epsilon_1\zeta_1,\epsilon_2\zeta_2](\psi_1,\psi_2)& \equiv \partial_x \big({\phi_1}\id{z=\epsilon_2\zeta_2}\big),\end{align*}
with $\phi_1$ and $\phi_2$ the unique solutions of the boundary problems 
\begin{equation}\label{phi}\left\{\begin{array}{ll}
 \Delta^\mu_{x,z}\phi_2=0 & \text{in } \Omega_2^t, \\
\phi_2 =\psi_2 & \text{on } \{z=\epsilon_2 \zeta_2\}, \\
 \partial_{z}\phi_2 =0 & \text{on } \{z=-\frac{1}{\delta}\}, \\
\end{array}
 \right. \quad \text{and} \quad \left\{\begin{array}{ll}
 \Delta^\mu_{x,z}\phi_1=0 & \text{in } \Omega_1^t, \\
\phi_1 =\psi_1 & \text{on } \{z=1+\epsilon_1 \zeta_1\}, \\
 \partial_{n}\phi_1 =\frac{G_2[\epsilon_2\zeta_2]\psi_2}{\sqrt{1+\mu \epsilon_2^2|\partial_x \zeta_2|^2}} & \text{on } \{z=\epsilon_2 \zeta_2\}.
\end{array}
\right.\end{equation}
Here, we denote by \[\Omega_1^t=\{(x,z)\in\RR^2, \epsilon_2\zeta_2(t,x)<z<1+\epsilon_1\zeta_1(t,x)\} \quad \text{ and }\quad \ {\Omega_2^t=\{(x,z)\in\RR^2, -\frac1\delta<z<\epsilon_2\zeta_2(t,x)\}}\] the domains of the fluids, by $\Delta^\mu_{x,z}\equiv\mu\partial_x^2+\partial_z^2$ the scaled Laplace operator, and by $\partial_n$  the upward conormal derivative:
\[(\partial_n\phi)\id{z=\epsilon_2\zeta_2}=\frac{1}{\sqrt{1+\mu\epsilon_2^2|\partial_x \zeta_2|^2}}(\partial_z\phi-\mu\epsilon_2\partial_x\zeta_2 \partial_x \phi)\id{z=\epsilon_2\zeta_2}.\]
\end{dfntn}
The domains of the two fluids are assumed to remain strictly connected, {\em i.e.} there exists $h_{\text{min}}>0$ such that 
\begin{equation}
\label{eq:h}
 \forall x\in\RR, \qquad h_1(x)\ \equiv \ 1+\epsilon_1\eta_1(x) \ \geq\ h_{\text{min}}\ >\ 0 \qquad \text{and} \qquad h_2(x)\ \equiv\ \frac{1}{\delta}+\epsilon_2\eta_2(x)\ \geq\ h_{\min}\ >\ 0.
\end{equation}
This assumption is necessary in order to obtain the consistency of the full Euler system~\eqref{FullEuler} with the Boussinesq/Boussinesq model~\eqref{Bouss}, as seen in~\cite{Duchene10}, and recalled in Proposition~\ref{Prop:ConsistencyBouss} below. We do not always precise this assumption thereafter.

\begin{rmrk}\label{rem:SurfTension}
 Even if the Cauchy problem associated to the Euler system at the interface of two fluids of different positive densities is known to be ill-posed in Sobolev spaces in the absence of surface tension (as Kelvin-Helmholtz instabilities appear), Lannes~\cite{Lannes10} proved thanks to a stability criterion that adding a small amount of surface tension guarantees the well-posedness of such a problem, with a time of existence that does not vanish as the surface tension goes to zero, and thus is consistent with the observations. The stability criterion states that the Kelvin-Helmholtz instabilities appear for high frequencies, where the regularization effect of the surface tension is relevant, while the main profile of the wave that we want to capture is located in lower frequencies, and is unaffected by surface tension. Consequently, we decide to neglect the surface tension term, as its effect does not appear in our asymptotic models.
 
 Furthermore, we know from Theorem 5.8 of~\cite{Lannes10} that in the long wave regime, there exists uniformly bounded solutions of the full Euler system for times of order $\O(1/\mu)$ (again with a small amount of surface tension) in the rigid lid case. With this result in mind, we assume in the following that smooth, uniformly bounded families of solutions to~\eqref{FullEuler}, whose existence is assumed in Proposition~\ref{Prop:ConvBoussEuler} and in Theorem~\ref{Th:CVKdVEuler}, indeed exist.
 \end{rmrk}

\section{Derivation and analysis of Boussinesq/Boussinesq models}
\label{Sec:DerivationAndAnalysisOfANewSymmetricModel}
As said in the introduction, the starting point of our study is the Boussinesq/Boussinesq model, obtained from~\eqref{FullEuler} thanks to an asymptotic expansion of the operators $G_1$, $G_2$ and $H$ (see~\cite{Duchene10}). In order to simplify the notations, we assume that the small parameters of the long wave regime are equal (the general case can easily be obtained by modifying some constants) and set
\[\epsilon_1\ =\ \epsilon_2\ =\ \mu\ \ \equiv \ \epsilon \ \ll\ 1.\]
The Boussinesq/Boussinesq system, in the one-dimensional case, can then be written using the convenient set of variables $(\eta_1,\eta_2,u_1,u_2)\ \equiv\ (\zeta_1-\zeta_2,\zeta_2,\partial_x\psi_1,\partial_x\psi_2)$ as 
\begin{equation}
\vspace{1mm}\label{Bouss} \left\{ \begin{array}{l}
\displaystyle\partial_t \eta_1 + \partial_x(h_1u_1) = 0,  \\ 
\displaystyle\partial_t \eta_2 +\partial_x(h_2 u_2)= 0,\\
\displaystyle\partial_t u_1 + \partial_x (\eta_1+\eta_2) +\epsilon\big(u_1\partial_x u_1 -\frac13  \partial_x^2\partial_t u_1- \frac{1}{2\delta} \partial_x^2\partial_t u_2  \big)= 0, \\ 
\displaystyle\partial_t u_2 + \partial_x (\gamma\eta_1+\eta_2) +\epsilon\big(u_2\partial_x u_2 -\frac{1+3\gamma\delta}{3\delta^2} \partial_x^2\partial_t u_2 - \frac{\gamma}{2}  \partial_x^2\partial_t u_1 \big)= 0,
\end{array} \right.  \end{equation}
where $h_1\ \equiv\ 1+\epsilon \eta_1$ and $h_2\ \equiv \frac{1}{\delta}+\epsilon\eta_2$ are the respective nondimensionalized depths of the upper and lower layer. Let us recall the following consistency result, that have been obtained in~\cite[Proposition 2.15]{Duchene10}
\begin{prpstn}\label{Prop:ConsistencyBouss}
 Let $U \ \equiv \ (\eta_1,\eta_2,u_1,u_2)$ be a strong solution of system~\eqref{FullEuler}, bounded in $L^{1,\infty}_t([0,T];H^{s+t_0})$ with $s>1$ and $t_0\geq9/2$, and such that~\eqref{eq:h} is satisfied. Then $U$ satisfies~\eqref{Bouss} up to a residual $R$ bounded by 
 \[\big|R \big|_{L^{\infty}H^{s}} \ \leq \ \epsilon^2 C_0 \Big(\frac1{h_{\min}},\big|U\big|_{L^{1,\infty}_tH^{s+t_0}}\Big).\] 
\end{prpstn}
\begin{rmrk}
 Here, and in the following, we denote by $C_0(\lambda_1,\lambda_2,\ldots)$ any positive constant, depending on the parameters ${\lambda_1,\lambda_2,\ldots }$, and whose dependence on $\lambda_j$ is assumed to be nondecreasing. 
 Moreover, for $0 < T \leq \infty$ and $f(t,x)$ a function defined on $[0,T]\times\RR$, we write $f\in L^\infty([0,T];H^s)$ if $f$ is uniformly (with respect to $t\in [0,T]$) bounded in $H^s=H^s(\RR)$ the $L^2$-based Sobolev space. Finally, one has $f\in L^{1,\infty}_t([0,T];H^s)$ if $f\in L^\infty([0,T];H^s)$ and $\partial_t f \in L^\infty([0,T];H^{s-1})$. Their respective norm is denoted $\big\vert\cdot\big\vert_{L^\infty H^s}$ and $\big\vert\cdot\big\vert_{L^{1,\infty}_tH^s}$.
\end{rmrk}

\subsection{A new family of symmetric models}
\label{sSec:ANewSymmetricModel}
In order to prove the well-posedness of hyperbolic systems as~\eqref{Bouss}, one can use energy methods, that require symmetries of the system. Although our system is not entirely symmetrizable, we show in this section that it is equivalent (in the sense of consistency) at order $\O(\epsilon^2)$ to a system of the form
\begin{equation} \label{SBOUSS}\Big(S_0+\epsilon \big(S_1(U)-S_2\partial_x^2\big)\Big)\partial_t U+\Big(\Sigma_0+\epsilon\big(\Sigma_1(U)- \Sigma_2\partial_x^2\big)\Big)\partial_x U=0,\end{equation}
which satisfies the following crucial properties:
\begin{assmptm}
 \label{Hyp}
 \begin{enumerate}
 \item The matrices $S_0$, $\Sigma_0$, $S_2$, $\Sigma_2 \in \M_4(\RR)$ are symmetric.
 \item $S_1(\cdot)$ and $\Sigma_1(\cdot)$ are linear mappings with values in $\M_4(\RR)$, and for all $U\in \RR^4$, the matrices $S_1(U)$ and $\Sigma_1(U)$ are symmetric.
 \item $S_0$ and $S_2$ are definite positive.
\end{enumerate}
\end{assmptm}

\begin{rmrk}
We sometimes write the system~\eqref{SBOUSS} under the form
\begin{equation} \label{OpBOUSS} P_\epsilon(U,\partial_x)\partial_t U + Q_\epsilon(U,\partial_x) \partial_x U=0,\end{equation}
with $P_\epsilon(U,\partial)=S_0+\epsilon \big(S_1(U)-S_2\partial^2\big)$ and $Q_\epsilon(U,\partial)=\Sigma_0+\epsilon\big(\Sigma_1(U)- \Sigma_2\partial^2\big)$.
\end{rmrk}

\bigskip

First of all, let us point out that system~\eqref{Bouss} is only one among many other Boussinesq-type systems. In~\cite{BonaChenSaut02}, Bona, Chen and Saut studied in the one-layer case a three-parameter family of Boussinesq systems, which are all approximations to the full Euler equations at the same order (in the sense of consistency). The same structure applies also in the two-layer case, and we describe it quickly below. We then exhibit first order symmetrizers {\em adapted to each of the Boussinesq-type systems}, and leading to systems of the form~\eqref{SBOUSS}. All of these models are equivalent in the sense that the original Boussinesq/Boussinesq system~\eqref{Bouss} is consistent at order $\O(\epsilon^2)$ with any of the systems presented in this section, as stated in Proposition~\ref{Prop:ConsistencyBoussSBOUSS}.

\medskip
 
As a first step, one can use the following change of variables:
\begin{equation}\label{eq:newvar}
v_2\ \equiv\ (1-\epsilon a_2\partial_x^2)^{-1}u_2, \quad \text{and} \quad v_1\ \equiv\ (1-\epsilon b_1\partial_x^2)^{-1} (u_1+\epsilon a_1\partial_x^2 v_2),
\end{equation}
with $a_1\in\RR$ and $a_2,b_1\geq 0$, and one recovers the three-parameter family of Boussinesq/Boussinesq systems introduced in~\cite{Duchene10}.

Then, as it has been achieved in~\cite{BonaChenSaut02}, one can also use the classical BBM-trick~\cite{BenjaminBonaMahony72}, and inherit new choices as parameters. This trick is based on the following calculations: since we have from~\eqref{Bouss} at first order
\begin{equation}\label{eq:BBMtrick}
 \begin{array}{rlrl}\partial_t \eta_1 &= -\partial_x v_1 +\O(\epsilon), & \qquad \partial_t \eta_2 &= -\frac1\delta\partial_x v_2 +\O(\epsilon), \\
 \partial_t v_1 &= -\partial_x (\eta_1+\eta_2) +\O(\epsilon), & \qquad \partial_t v_2 &= -\partial_x (\gamma\eta_1+\eta_2)+\O(\epsilon) , \end{array}\end{equation}
we can deduce the following, with the parameters $\lambda_i\in[0,1]$ ($i=1\dots 4$),
\begin{align*}
 \partial_x^3 v_1&=\lambda_1\partial_x^3 v_1-(1-\lambda_1)\partial_x^2\partial_t \eta_1 +\O(\epsilon),\\
 \partial_x^3 v_2&=\lambda_2\partial_x^3 v_2-\delta(1-\lambda_2)\partial_x^2\partial_t \eta_2 +\O(\epsilon),\\
 \partial_x^2\partial_t v_1&=(1-\lambda_3)\partial_x^2\partial_t v_1-\lambda_3\partial_x^3 (\eta_1+\eta_2) +\O(\epsilon),\\
 \partial_x^2\partial_t v_2&=(1-\lambda_4)\partial_x^2\partial_t v_2-\lambda_4\partial_x^3  (\gamma\eta_1+\eta_2) +\O(\epsilon).
\end{align*}

In the end, one obtains the following system, formally equivalent to~\eqref{Bouss} at order $\O(\epsilon^2)$:
\begin{equation} \label{BOUSS}
\partial_t U+\A_0\partial_x U+\epsilon\big(\A(U)\partial_x U-\A_1\partial_x^3 U-\A_2\partial_x^2\partial_t U \big)=0,
\end{equation}
denoting $U=(\eta_1,\eta_2,v_1,v_2)$ and
\[\begin{array}{c}
	    \A_0=\begin{pmatrix}
  0 & 0 & 1 & 0 \\
  0 & 0 & 0 &1/\delta \\
  1 & 1 & 0 & 0\\
  \gamma & 1 & 0 & 0
            \end{pmatrix}, \ \ 
	    \A_1=\begin{pmatrix}
  0 & 0 & -\lambda_1\beta_1 & -\lambda_2\alpha_1 \\
  0  & 0 & 0 & -\lambda_2\frac{\alpha_2}{\delta} \\
  -\lambda_3b_1 -\lambda_4 \gamma a_1& -\lambda_3b_1 -\lambda_4 a_1& 0 & 0\\ -\lambda_4\gamma\beta_2 -\lambda_3\frac{\gamma}{2}& -\lambda_4\beta_2-\lambda_3\frac{\gamma}{2} & 0 & 0
            \end{pmatrix},\\   
	    \A_2=\begin{pmatrix}
  (1-\lambda_1)\beta_1 & \delta(1-\lambda_2)\alpha_1 &0 & 0 \\
  0 & (1-\lambda_2)\alpha_2  &0 & 0 \\
  0 & 0 & (1-\lambda_3)b_1 &  (1-\lambda_4)a_1 \\
  0 &0 & (1-\lambda_3)\frac{\gamma}{2} & (1-\lambda_4)\beta_2 
             \end{pmatrix}, \ \
	        \A(U)=\begin{pmatrix}
  v_1 & 0 & \eta_1 & 0 \\
  0  & v_2 & 0 & \eta_2 \\
  0 & 0 & v_1 & 0\\
  0 & 0 & 0 & v_2
              \end{pmatrix},	 
\end{array}\]
with the parameters $\beta_1=\frac{1}{3}-b_1$, $\alpha_1=\frac{1}{2\delta}-a_1$, $\alpha_2=\frac{1}{3\delta^2}-a_2$ and $\beta_2=a_2+\frac{\gamma}{\delta}$, so that the system depends on the real parameters $a_1$, $a_2$, $b_1$, and $\lambda_i$ ($i=1\dots 4$), that can be chosen freely.

\bigskip

In order to derive a system of the form~\eqref{SBOUSS}, we exhibit a good symmetrizer of~\eqref{BOUSS}, that is to say \[S(U)\equiv S_0+\epsilon S_1(U)-\epsilon \t S_2 \partial_x^2,\] such that when we multiply~\eqref{BOUSS} on the left by $S(U)$, and withdrawing the $\O(\epsilon^2)$ terms, we obtain a system~\eqref{SBOUSS} satisfying Assumption~\ref{Hyp}.

The symmetrization of the one-layer shallow water system is well known, and consists in multiplying the velocity equation by the water depth. An adaptation of this to our two-layer model leads to 
\begin{equation}\label{S0S1}S_0\equiv\left( \begin {array}{cccc} \gamma&\gamma&0&0\\
\gamma&1&0&0 \\
0&0&\gamma&0\\
0&0&0&1/\delta
\end {array} \right) \qquad \text{and} \qquad S_1(U)\equiv \left( \begin {array}{cccc} 0&0&\gamma v_1&0\\
0&0&0&v_2 \\
\gamma v_1&0&\gamma \eta_1&0\\
0&v_2&0&\eta_2
\end {array} \right),\end{equation}
so that $S_0$ and $\Sigma_0\equiv S_0\A_0$ are symmetric, and for all $U\in\RR^4$, $S_1(U)$ and ${\Sigma_1(U)\equiv S_1(U)\A_0+S_0\A(U)}$ are symmetric. Moreover, $S_0$ is definite positive for $\delta>0$ and $\gamma\in(0,1)$: its eigenvalues are \[\gamma,\ 1/\delta,\ \text{and} \ \ 1/2(1+\gamma \pm \sqrt{(1-\gamma)^2+4\gamma^2}).\]

Then, one can check that when we set
\[\t S_2\equiv\left( \begin {array}{cccc} a+(\delta-1+\gamma) b  &0 &0&0\\
a+ \delta b&0&0&0 \\
0&0&a+(\delta-1)b+\gamma(b_1\lambda_3-\lambda_1\beta_1)&b+\gamma \lambda_4 a_1\\
0&0&\gamma(\frac{\lambda_3}{2\delta}-\lambda_2\alpha_1)&\frac{\lambda_4\beta_2-\lambda_2\alpha_2}{\delta}
\end {array} \right)\]
with $a=\gamma((1-\lambda_2)\alpha_2-(1-\lambda_1)\beta_1)$ and $b= \gamma(1-\lambda_2)\alpha_1$, then 
$S_2\equiv S_0\A_2+\t S_2$ and $\Sigma_2\equiv S_0\A_1+\t S_2\A_0$ are symmetric. Then for any $K\in\RR$, we can substitute $\t S_2+K S_0$ for $\t S_2$, and one has again $S_2$ and $\Sigma_2$ are symmetric. Moreover, since $S_0$ is definite positive, one can choose $K$ big enough for $S_2$ to be definite positive.

Therefore, when we multiply~\eqref{BOUSS} by $S(U)=S_0+\epsilon S_1(U)-\epsilon (\t S_2+K S_0)\partial_x^2$, and withdrawing the $\O(\epsilon^2)$ terms, we obtain the perfectly symmetric system~\eqref{SBOUSS}.

\bigskip

Using the above calculations, it is now straightforward to obtain the following consistency result:
\begin{prpstn}\label{Prop:ConsistencyBoussSBOUSS}
 Let $U=(\eta_1,\eta_2,u_1,u_2)$ be a strong solution of system~\eqref{Bouss} such that ${V=(\eta_1,\eta_2,v_1,v_2)}$, given by the change of variables~\eqref{eq:newvar}, is uniformly bounded in $L^{1,\infty}_t([0,T];H^{s+5})$ with $s>1/2$. Then $V$ satisfies~\eqref{BOUSS} and~\eqref{SBOUSS} up to a residuals bounded by $\epsilon^2 C_0$ in the sense of $L^\infty H^s$ norm, with 
 \[C_0=C_0 \big(a_1,a_2,b_1,K,\delta+\frac1\delta\big)\big|V\big|_{L^{1,\infty}_t([0,T];H^{s+5})}.\] 
\end{prpstn}
\begin{proof}
The first step in order to obtain~\eqref{BOUSS} from~\eqref{Bouss} is to use the change of variables~\eqref{eq:newvar}. Thus, when replacing in $u_1$ by $v_1-\epsilon b_1 \partial_x^2 v_1-\epsilon a_1 \partial_x^2 v_2$ and $u_2$ by $v_2-\epsilon a_2 \partial_x^2 v_2$, we obtain straightforwardly that $V$ satisfies~\eqref{BOUSS} in the case $1-\lambda_1=1-\lambda_2=\lambda_3=\lambda_4=0$, up to terms of the form  \[\epsilon^2 \partial_x^4\partial_t f, \quad \epsilon^2 \partial_x (f\partial_x^2 g)\quad \text{and}  \quad
\epsilon^3 \partial_x \big((\partial_x^2 f)^2\big),\] where $f$ and $g$ are components of $V$. Using the fact that $V$ is bounded in $L^{1,\infty}_t([0,T);H^{s+5})$, and $H^s(\RR)$ is an algebra for $s>1/2$, the remaining terms are clearly bounded by $\epsilon^2 C_0(a_1,a_2,b_1,K,\delta+\frac1\delta,) \big|V\big|_{L^{1,\infty}_tH^{s+5}}$.

In the same way, when we substitute the relations of the BBM trick~\eqref{eq:BBMtrick} into the third-derivative terms of the equations, we obtain~\eqref{BOUSS} up to extra terms bounded by $\epsilon^2 C_0(a_1,a_2,b_1,\frac1\delta)\big|V\big|_{L^{1,\infty}_tH^{s+5}}$.

 Finally, in order to obtain~\eqref{SBOUSS}, we multiply~\eqref{BOUSS} by $S_0+\epsilon S_1(V)-\epsilon(\t S_2+K S_0)\partial_x^2$, and withdraw the terms \[\epsilon^2 \big(S_1(V)-(\t S_2+K S_0)\partial_x^2\big)\big(\A(V)\partial_x V-\A_1\partial_x^3 V-\A_2\partial_t\partial_x^2 V\big) .\] 
 Each of these terms are clearly bounded by $\epsilon^2 C_0 $, since $H^s(\RR)$ is an algebra for $s>1/2$, and
 \[\big|A(V)\partial_x V\big|_{H^{s}}\leq C_0\big|V\big|_{L^\infty H^{s+3}} \text{, }\  \big|S_1(V)\big|_{L^\infty}\leq C_0\big|V\big|_{L^\infty H^s},
  \text{ and }\ \big|\A_1 \partial_x^3V\big|_{H^{s}}+\big|\A_2 \partial_t\partial_x^2 V\big|_{H^{s}}\leq C_0\big|V\big|_{L^{1,\infty}_t H^{s+1}}
 \]
  Therefore, the system~\eqref{Bouss} is consistent with the system~\eqref{SBOUSS} at the precision $\epsilon^2 C_0 $.
\end{proof}

\subsection{Well-posedness and convergence results}
\label{sSec:WPAndCVResults}
The system~\eqref{SBOUSS} is a symmetric hyperbolic system, and can be studied using classical energy methods. We first prove in Proposition~\ref{Prop:WPSBOUSS} that such a system is well-posed over times of order $\O(1/\epsilon)$. The proof uses in particular an a priori estimate of the solution, in an adapted norm:
\[\big| U\big|_{H^{s+1}_\epsilon}^2 \ \equiv \ \big| U\big|_{H^{s}}^2 \ + \ \epsilon \big| U\big|_{H^{s+1}}^2.\]
Then, using a consistency result with energy estimates, we show that the solutions of our models converge towards bounded solutions of the full Euler system, assuming that such solutions exist (see Remark~\ref{rem:SurfTension}).

\medskip

Here and thereafter, we fix $0<\gamma_{\text{min}}\leq\gamma\leq\gamma_{\text{max}}<1$ and $0<\delta_{\text{min}}\leq\delta\leq\delta_{\text{max}}<+\infty$. The limit cases ($\delta \to 0,\infty$ and $\gamma\to0,1$) demand other scalings in the nondimensionalization than the ones presented in~\cite{Duchene10} (see Section A of~\cite{Lannes10} for example) and correspond to different regimes, such as the deep-water theory (from Benjamin~\cite{Benjamin67} and Ono~\cite{Ono75}), and lead to different models (see for example~\cite{CraigGuyenneKalisch05,BonaLannesSaut08} in the rigid lid configuration, and~\cite{SegurHammack82,KoopButler81,OstrovskyStepanyants05} in the free surface case). In all of these cases, the calculations of our justification break: the dependence of the constants ${C_0(\frac{1}{\gamma(1-\gamma)},\delta+\frac1\delta)=C_0(\big\|S_0\big\|,\big\|S_0^{-1}\big\|)}$ in the following theorems prevents the parameters to approach the limits.

In the same way, the following constants $C_0$ also depend on the set of parameters $(a_1,a_2,b_1,K,\lambda_i)$ $(i=1\cdots4)$. We decide to fix these to constants once for all, and do not write explicitly this dependency, in order to simplify the notations.

\bigskip

We now state the well-posedness of our symmetric Boussinesq/Boussinesq model.
\begin{prpstn}
\label{Prop:WPSBOUSS} 
 Let $U^0\in H^{s+1}$, with $s>3/2$. Then there exists a constant ${C_0=C_0(\frac{1}{\gamma(1-\gamma)},\delta+\frac1\delta)>0}$ such that for $\epsilon \leq \epsilon_0=(C_0 \big\vert U^0\big\vert_{H^{s+1}_\epsilon})^{-1}$, there exists a time ${T>0}$ independent of $\epsilon$, and a unique solution ${U\in C^0([0,T/\epsilon);H^{s+1}_\epsilon)\cap C^1([0,T/\epsilon);H^{s}_\epsilon)}$ of the Cauchy problem~\eqref{SBOUSS} with ${U\id{t=0}=U^0}$. 
 
 Moreover, one has the following estimate for $t\in[0,T/\epsilon]$:
\begin{equation}\label{EstWP}
\big| U\big|_{L^\infty([0,t];H^{s+1}_\epsilon)}+ \big| \partial_t U\big|_{L^\infty([0,t] ; H^{s}_\epsilon)} \leq C_0  \frac{\big| U^0\big|_{H^{s+1}_\epsilon}}{1-C_0  \big| U^0\big|_{H^{s+1}_\epsilon} \epsilon t}.
\end{equation}
\end{prpstn}
We postpone the somewhat technical proof to Appendix~\ref{Sec:ProofPropWPSBOUSS}. The key ingredients of the proof are quickly presented in the following remark.
\begin{rmrk}
The condition $s>3/2$ is necessary in order to obtain estimate~\eqref{EstWP}, and thus the well-posedness over times of order $\O(1/\epsilon)$. One could obtain, using the exact same method as in the proof, the same result over times of order $\O(1)$, with the sharper assumption $s>1/2$ (which is the standard regularity for one-dimensional hyperbolic systems). 

The smallness condition ($\epsilon \big\vert U^0\big\vert_{H^{s+1}_\epsilon} \leq {C_0}^{-1}$) is also necessary, for the nonlinear terms in $S_1(U)$ to remain negligible when compared with $S_0$. Indeed, under this condition, the energy of the system~\eqref{SBOUSS}, defined by
\[E_s(U)\ \equiv\  1/2(S_0 \Lambda^s U,\Lambda^s U)\ +\ \epsilon/2(S_1(U) \Lambda^s U,\Lambda^s U)\ +\ \epsilon/2(S_2 \Lambda^s \partial_x U,\Lambda^s \partial_x U),\]
is uniformly equivalent to the $\big\vert\cdot\big\vert_{H^{s+1}_\epsilon}$ norm, that is to say there exists $\alpha>0$ such that
\[
 \frac{1}{\alpha}\big\vert U\big\vert_{H^{s+1}_\epsilon}^2\ \leq\  E_s(U)\ \leq\ \alpha \big\vert U\big\vert_{H^{s+1}_\epsilon}^2.
\]
Moreover, the smallness condition on $\epsilon\big\vert U^0\big\vert_{H^{s+1}_\epsilon}$ allows us to prove that the operator \[P_\epsilon(U,\partial)\ =\ S_0+\epsilon \big(S_1(U)-S_2\partial^2\big):H^{s+1}\to H^{s-1}\] is one-to-one and onto, and that $P_\epsilon(U,\partial)^{-1}Q_\epsilon(U,\partial)$ is uniformly bounded $H^{s}_\epsilon\to H^{s}_\epsilon$. This leads to
\[\big\vert\partial_t U\big\vert_{H^{s}_\epsilon} = \big\vert P_\epsilon(U,\partial)^{-1}Q_\epsilon(U,\partial)\partial_x U \big\vert_{H^{s}_\epsilon}\leq C_0 \big\vert U\big\vert_{H^{s+1}_\epsilon}.\]
Both of these properties are crucial in order to prove estimate~\eqref{EstWP}.

The existence and uniqueness of the solution of the Cauchy problem follow from the a priori estimate, using classical methods.
\end{rmrk}

\bigskip

Using the previous Proposition, and the consistency of the full Euler system~\eqref{FullEuler} with our symmetric Boussinesq/Boussinesq model~\eqref{SBOUSS}, one can now easily deduce the following convergence Proposition:
\begin{prpstn} \label{Prop:ConvBoussEuler} Let $s>3/2$, $\epsilon>0$ and $U=(\zeta_1,\zeta_2,\psi_1,\psi_2)$ be a solution of the full Euler system~\eqref{FullEuler} such that $V=(\eta_1,\eta_2,v_1,v_2) \in C^0([0;T/\epsilon);H^{s+1}_\epsilon)\cap C^1([0;T/\epsilon);H^{s}_\epsilon)$ defined by \[
V \equiv(\zeta_1-\zeta_2,\zeta_2, (1-\epsilon b_1\partial_x^2)^{-1} (\partial_x\psi_1+\epsilon a_1\partial_x^2 v_2),(1-\epsilon a_2\partial_x^2)^{-1}\partial_x\psi_2 )\] is uniformly bounded in $L^{1,\infty}_t([0,T/\epsilon];H^{s+5})$ and~\eqref{eq:h} is satisfied. We assume that $\epsilon$ satisfies the smallness condition of Proposition~\ref{Prop:WPSBOUSS}, and denote by $V_B$ the solution of the symmetric Boussinesq/Boussinesq system~\eqref{SBOUSS}, with the same initial value ${V_B}\id{t=0}={V}\id{t=0}=V^0$ and the same domain of existence. 
Then one has for all $t\in[0,T/\epsilon]$, \[\big|V-V_B\big|_{L^\infty([0,t] ; {H}^{s+1}_\epsilon)}\leq \epsilon^2 t C_0,\]
with $C_0=C_0 (\frac{1}{h_{\min}},\frac{1}{\gamma(1-\gamma)},\delta+\frac1\delta,\big|V\big|_{L^{1,\infty}_tH^{s+5}},T)$. In particular, one has 
\[\big|V-V_B\big|_{L^\infty([0,T/\epsilon] ; {H}^{s+1}_\epsilon)}\leq \epsilon C_0,\]
with $C_0$ independent of $\epsilon$.
\end{prpstn}
\begin{proof}
Let us first point out that the full Euler system~\eqref{FullEuler} is consistent with~\eqref{SBOUSS} at the precision $\epsilon^2 C_0 $, that is to say that for any solution $U$ of the full Euler system such that $V\in C^0([0;T/\epsilon);H^{s+1})\cap C^1([0;T/\epsilon);H^{s})$ is uniformly bounded in $L^{1,\infty}_tH^{s+5}([0,T/\epsilon])$, then $V$ satisfies~\eqref{SBOUSS} up to a residual bounded by $\epsilon^2 C_0$ in the sense of $L^\infty H^s$ norm, with $C_0 =C_0(\frac{1}{h_{\min}},\delta+\frac1\delta,\big\vert V\big\vert_{L^{1,\infty}_tH^{s+5}})$. Proposition~\ref{Prop:ConsistencyBouss} states that the full Euler system is consistent with~\eqref{Bouss} at the precision $\epsilon^2 C_0 $, and Proposition~\ref{Prop:ConsistencyBoussSBOUSS} achieves the result.

Therefore, we know that $V$ satisfies~\eqref{SBOUSS} up to $\epsilon^2f$, with $f\in L^\infty([0,T/\epsilon];H^s)$, so that $R_s\equiv \Lambda^sV-\Lambda^sV_B$, $\Lambda^s$ being the Fourier multiplier defined by $\widehat{\Lambda^s u}(\zeta)\equiv (1+|\zeta|^2)^{s/2}\widehat u(\zeta)$, satisfies the system
\begin{equation} \label{eqRs} \big(S_0-\epsilon S_2\partial_x^2\big)\partial_t R_s+\epsilon \Lambda^s(S_1(V)\partial_t R_0+S_1(R_0)\partial_t V_B)+\big(\Sigma_0-\epsilon \Sigma_2\partial_x^2\big)\partial_x R_s 
+\epsilon \Lambda^s(\Sigma_1(V) \partial_x R_0+\Sigma_1(R_0) \partial_x V_B) =\epsilon^2 \Lambda^s f,\end{equation}
with $\big\vert f\big\vert_{L^\infty([0,T/\epsilon];H^s)}\leq C_0 (\frac{1}{h_{\min}},\delta+\frac1\delta,\big|V\big|_{L^{1,\infty}_tH^{s+5}}) $. 

\medskip

We can then carry on the calculations of Section~\ref{ProofUnicity}, with the extra term $\epsilon^2 \Lambda^s f$. We obtain that there exists $C_0(\frac{1}{\gamma(1-\gamma)},\delta+\frac1\delta)$ such that as long as \begin{equation}\label{smallR}\epsilon \big|R_0\big|_{H^{s+1}_\epsilon} \leq 1/C_0,\end{equation} one has the estimate
\begin{equation}
 \frac{d}{dt}E(R_s)\ \leq\ \epsilon C_0 (\big|V_B\big|_{H^s}+\big| V\big|_{H^s})\big| R_0\big|_{H^{s}}^2\ +\ \epsilon^2(\Lambda^{s} f, \Lambda^{s}R_s) ,
\end{equation}
with the energy $E(R_s)$ defined by
 \[E(R_s)\ \equiv\ 1/2(S_0 R_s,R_s)\ +\ \epsilon/2(S_1(R_0) R_s,R_s)\ +\ \epsilon/2(S_2 \partial_x R_s, \partial_x R_s).\]
 
 Now, since $S_0$ and $S_2$ are definite positive, the condition~\eqref{smallR} implies in particular
\begin{equation}
 \frac{1}{C_0}\big\vert R_0\big\vert_{H^{s+1}_\epsilon}^2\ \leq\  E(R_s)\ \leq\ C_0 \big\vert R_0\big\vert_{H^{s+1}_\epsilon}^2.
\end{equation}

Thus, under this condition, and since $V_B$ is uniformly bounded in $L^\infty([0,T/\epsilon];H^s)$ with respect to $\epsilon$ (from Proposition~\ref{Prop:WPSBOUSS}), one has
\[ \frac{d}{dt}E(R_s)\leq \epsilon C_0  E(R_s)+\epsilon^2 C_0  \big\vert f \big\vert_{H^s} E(R_s)^{1/2},\]
and Gronwall-Bihari's Lemma leads to
\[ E(R_s)^{1/2} \ \leq\  C_0 \epsilon\big\vert f \big\vert_{H^s}(e^{\epsilon C_0  t}-1).\]

Finally, since ${R_0}\id{t=0}=0$, and thanks to a continuity argument, there exists $T(C_0,\big|f\big|_{H^s})>0$ such that~\eqref{smallR} holds for $0\leq t \leq T/\epsilon$, and the estimate of the Proposition follows: 
\[\big\vert R_0\big\vert_{H^{s+1}_\epsilon}\ \leq\ C_0 E(R_s)^{1/2}\ \leq\ C_0 \epsilon^2 \big\vert f \big\vert_{H^s} t.\] 
\end{proof}

\section{The KdV approximation}
\label{Sec:TheKdVApproximation}
In this section, we offer a rigorous justification of the so-called KdV approximation, from a class of symmetric systems that contains the symmetric Boussinesq/Boussinesq system~\eqref{SBOUSS}. The KdV approximation consists in a decomposition of the flow into four parts, each of the components being lead by a Korteweg-de Vries equation.
The construction of the KdV approximation is precisely explained in Section~\ref{sSec:FormalDerivation}. Then, in Section~\ref{sSec:RigorousDemonstration}, we obtain the convergence rate between the solutions of the coupled systems and the solutions defined by the KdV approximation. 
As a consequence, when we combine this result with the convergence Proposition~\ref{Prop:ConvBoussEuler}, it follows immediately that any strong solution of the full Euler system existing over times $\O(1/\epsilon)$ and bounded in a sufficiently high Sobolev norm, is well approximated by the KdV approximation. More precisely, we state the following:
\begin{thrm}\label{Th:CVKdVEuler} 
Let $s>3/2$ and $U=(\zeta_1,\zeta_2,\psi_1,\psi_2)$ be a solution of the full Euler system~\eqref{FullEuler} such that $V=(\eta_1,\eta_2,v_1,v_2) \in C^0([0;T/\epsilon);H^{s+1}_\epsilon)\cap C^1([0;T/\epsilon);H^{s}_\epsilon)$ defined by \[V\ \equiv\ (\ \zeta_1-\zeta_2\ ,\ \zeta_2\ ,\ (1-\epsilon b_1\partial_x^2)^{-1} (\partial_x\psi_1+\epsilon a_1\partial_x^2 v_2)\ ,\ (1-\epsilon a_2\partial_x^2)^{-1}\partial_x\psi_2\ )\] is uniformly bounded in $L^{1,\infty}_t([0,T/\epsilon];H^{s+5})$ and~\eqref{eq:h} is satisfied. Then there exists\footnote{one has explicit expressions for the basis $(\e_1,\dots,\e_4)$ and the coefficients $c_i,\lambda_i,\mu_i$, that we display in Remark~\ref{rem:coefsKdV} page \pageref{rem:coefsKdV}.} $(\e_1,\dots,\e_4)$ a basis of $\RR^4$, and coefficients $c_i,\lambda_i,\mu_i$ ($i=1\cdots4$) such that, denoting by $u_i$ the solution of the KdV equation
\begin{equation}
 \partial_t u_i \ +\ c_i\partial_x u_i \ +\ \epsilon \lambda_i u_i\partial_x u_i\ +\ \epsilon\mu_i \partial_x^3 u_i\ =\ 0
\end{equation}
with ${u_i}\id{t=0}=\e_i\cdot S_0 V\id{t=0}$ ($S_0$ defined in~\eqref{S0S1}), one has for all $t\in[0,T/\epsilon]$,
\[\big|V-\sum_{i=1}^4 u_i \e_i \big|_{L^\infty([0,t] ; {H}^{s+1}_\epsilon)}\leq \epsilon \sqrt t C_0,\]
with $C_0=C_0 (\frac{1}{h_{\min}},\frac{1}{\gamma(1-\gamma)},\delta+\frac1\delta,\big|V\big|_{L^{1,\infty}_tH^{s+5}})$.

\medskip
 
Moreover, if $V$ satisfies ${(1+x^2) V\id{t=0} \in H^{s+4}}$, then one has the better estimate
\[\big|V-\sum_{i=1}^4 u_i \e_i\big|_{L^\infty([0,T/\epsilon]; H^{s+1}_\epsilon)}\leq \epsilon C_0',\]
with $C_0'=C_0 (\frac{1}{h_{\min}},\frac{1}{\gamma(1-\gamma)},\delta+\frac1\delta,\big|V\big|_{L^{1,\infty}_tH^{s+5}},\big\vert(1+x^2) V\id{t=0}\big\vert_{H^{s+4}})$.
\end{thrm}
\begin{proof}
The proof proceeds from different results of the paper; the completion is as follows. In Proposition~\ref{convergenceKdVBouss}, we prove the convergence between the solutions of systems of the form~\eqref{eq:BoussForKdV} (thus containing symmetric Boussinesq/Boussinesq systems~\eqref{SBOUSS}) and the approximate solution ${U_{\text{app}}=\sum_{i=1}^4 u_i \e_i+\epsilon U_1}$, defined in Definition~\ref{Def:Uapp}. The residual $U_1$ is then estimated in Proposition~\ref{EstU1}. Finally, since we have from Proposition~\ref{Prop:ConvBoussEuler} the convergence between the solutions of the full Euler system~\eqref{FullEuler} and the solutions of the symmetric Boussinesq/Boussinesq system~\eqref{SBOUSS} with a better rate, the Theorem follows; see also Remark~\ref{rem:proof} below.
\end{proof}
\begin{rmrk}
 The difference on the convergence rate for different sets of initial values is not simply a technical issue. Indeed, one can see in Figure~\ref{fig:ErrorInTime}, page~\pageref{fig:ErrorInTime}, that if the condition of sufficient decreasing in space is not satisfied, the convergence will be worse than $\O(\epsilon)$. Requiring the initial data (and thus the solutions of the KdV equations) to lie in weighted Sobolev spaces guarantees that the nonlinear interaction between the four traveling waves can be neglected. As a matter of fact, this condition on the sufficient decay in space of the initial data appears also naturally for the KdV approximation of the one-layer problem, as we see in~\cite{SchneiderWayne00,BonaColinLannes05}.
\end{rmrk}

\subsection{Formal derivation} \label{sSec:FormalDerivation}
The class of system that we now study is the following:
\begin{equation}\label{eq:BoussForKdV}\Big(S_0+\epsilon \big(S_1(U)- S_2\partial_x^2\big)\Big)\partial_t U+\Big(\Sigma_0+\epsilon\big(\Sigma_1(U)-\Sigma_2 \partial_x^2\big)\Big)\partial_x U=0,\end{equation}
with the following hypothesis:
\begin{assmptm}\label{AssForKdV}
\begin{enumerate}
 \item The matrices $S_0$, $\Sigma_0$, $S_2$, $\Sigma_2 \in \M_4(\RR)$ are symmetric.
 \item $S_1(\cdot)$ and $\Sigma_1(\cdot)$ are linear mappings with values in $\M_4(\RR)$, and for all $U\in \RR^4$, $S_1(U)$ and $\Sigma_1(U)$ are symmetric.
 \item $S_0$ is definite positive, and $S_0^{-1}\Sigma_0$ has four different non zero eigenvalues $c_i$ $(i=1\dots 4)$.
\end{enumerate}
\end{assmptm}

\begin{rmrk}
The symmetric Boussinesq/Boussinesq systems~\eqref{SBOUSS} derived in Section~\ref{sSec:ANewSymmetricModel} immediately satisfy Assumption~\ref{AssForKdV}, with \[c_i=\pm\sqrt{\frac{1+\delta\pm\sqrt{(1-\delta)^2+4\gamma\delta}}{2\delta}}.\]
\end{rmrk}

Following the classical WKB method, we look for an approximate solution of the Cauchy problem~\eqref{eq:BoussForKdV} with initial data $U^0$ under the form
\[U_{\text{app}}(t,x)=U_0(\epsilon t,t,x)+\epsilon U_1(\epsilon t,t,x),\]
with the profiles $U_0(\tau,t,x)$ and $\epsilon U_1(\tau,t,x)$ satisfying ${U_0}\id{t=\tau=0}=U^0$ et ${U_1}\id{t=\tau=0}=0$.

We plug the Ansatz into~\eqref{eq:BoussForKdV}, and obtain
\begin{align}\label{eq:AnsatzInBouss}
 (S_0\partial_t +\Sigma_0\partial_x) U_0+\epsilon S_0\partial_\tau U_0+\epsilon\big(S_1(U_0)\partial_t U_0+\Sigma_1(U_0)\partial_x U_0-S_2\partial_x^2\partial_t U_0 -\Sigma_2\partial_x^3 U_0\big)  \nonumber \\ +\epsilon(S_0\partial_t +\Sigma_0\partial_x) U_1+\epsilon^2 R =0 .
\end{align}
We now deduce the equations satisfied by $U_0(\tau,t,x)$ and $U_1(\tau,t,x)$, solving~\eqref{eq:AnsatzInBouss} at each order.

\bigskip

\noindent{\bf At order $\O(1)$:} We solve
\begin{equation}\label{eq:TRANSPORT}(S_0\partial_t +\Sigma_0\partial_x) U_0=0.\end{equation}
We assumed that $S_0$ is symmetric definite positive, and hence induce a scalar product: \[\langle u,v \rangle\equiv u\cdot S_0 v=u^T S_0 v.\]
Since $S_0^{-1}\Sigma_0$ is real and symmetric for the scalar product $\langle\cdot,\cdot\rangle$, it is diagonalizable in an orthonormal basis. We denote by $\e_i$ $(i=1\dots 4)$ the basis vectors, which are the unitary eigenvectors of $S_0^{-1}\Sigma_0$. By definition, they satisfy for $1\leq i,j\leq 4$,\[\langle\e_i,S_0^{-1}\Sigma_0 \e_j\rangle\ \equiv\ \e_i\cdot \Sigma_0 \e_j\ =\ c_i\delta_{i,j}, \quad \text{ and  \quad }\langle\e_i,\e_j\rangle\ \equiv\ \e_i\cdot S_0 \e_j\ =\ \delta_{i,j},\]
with $\delta_{i,j}$ the classical Kronecker delta symbol. Therefore, when we define $u_i\equiv \e_i\cdot S_0 U_0$ (and hence $U_0=\sum\limits_{i=1}^4 u_i \e_i$), multiplying~\eqref{eq:TRANSPORT} on the left by $\e_i$, we obtain \[(\partial_t +c_i\partial_x) u_i \ =\ 0\]
for all $i=1,\dots,4$. Finally, since $u_i$ satisfies a transport equation, we use the notation
\begin{equation}
 \label{eq:transport} u_i(\tau,t,x)\ =\ u_i(\tau,x-c_i t)\ =\ u_i(\tau,x_i),
\end{equation}
with initial data ${u_i}(0,x_i)=\e_i\cdot S_0 U^0(x_i)$.

\bigskip

\noindent{\bf At order $\O(\epsilon)$:} We solve
\begin{equation}\label{eq:Order2} S_0\partial_\tau U_0+\Sigma_1(U_0)\partial_x U_0+S_1(U_0)\partial_t U_0-\Sigma_2\partial_x^3 U_0-S_2\partial_x^2\partial_t U_0 +(S_0\partial_t +\Sigma_0\partial_x) U_1=0,\end{equation}
that we can split\footnote{This splitting is in fact necessary. Indeed, the multiscale WKB expansion can be justified for times of order $\O(1/\epsilon)$ only if the growth of the corrector term $U_1$ is sublinear. As we see in Proposition~\ref{EstU1}, thanks to Lemma~\ref{LannesLemma}, the particular form of~\eqref{eqnrest} allows to obtain a square-root growth, and even better if $U^0$ is sufficiently decreasing in space.} in
\begin{equation}\label{KdV}\partial_\tau u_i\ +\ \lambda_i u_i\partial_{x_i} u_i\ +\ \mu_i \partial_{x_i}^3 u_i\ =\ 0,\end{equation}
with $\lambda_i\equiv \e_i\cdot\big(\Sigma_1(\e_i)-c_i S_1(\e_i)\big)\e_i$ and $\mu_i\equiv \e_i\cdot \big(-\Sigma_2+c_i S_2\big)\e_i$; and in the other hand,
\begin{equation}\label{eqnrest} (\partial_t +c_i\partial_x) \langle\e_i, U_1\rangle + \sum_{(j,k)\neq(i,i)}\alpha_{ijk} u_k(\tau,x- c_k t)\partial_x u_j(\tau,x- c_j t)=\sum_{j\neq i} \beta_{ij}\partial_x^3 u_j(\tau,x- c_j t),\end{equation}
with $\alpha_{ijk}\equiv\e_i\cdot (\Sigma_1(\textbf{e} _k) -c_j S_1(\textbf{e} _k)) \e_j$ and $\beta_{ij}\equiv\e_i \cdot (\Sigma_2-c_j S_2)\e_j$.

\medskip

It is clear that $u_i$ satisfies~\eqref{eq:transport} and~\eqref{KdV}, if and only if $u_i(\epsilon t,t,x)$ satisfies the Korteweg-de Vries equation of Theorem~\ref{Th:CVKdVEuler}:
\begin{equation}\label{eq:KdV}
 \partial_t u_i\ +\ c_i\partial_x u_i\ +\ \epsilon \lambda_i u_i\partial_x u_i\ +\ \epsilon\mu_i \partial_x^3 u_i\ =\ 0.
\end{equation}

\begin{rmrk}\label{rem:coefsKdV} 
In the specific case of symmetric Boussinesq/Boussinesq systems~\eqref{SBOUSS}, one obtains the following values for the coefficients:
\begin{align*}
 (c_i)_{i\in\{1\dots4\}}&=(c_+,-c_+,c_-,-c_-), & &\text{ with } c_\pm^2=\frac{1+\delta\pm\sqrt{(1-\delta)^2+4\gamma\delta}}{2\delta}, \\
 (\lambda_i)_{i\in\{1\dots4\}}&= (\lambda_+,\lambda_+,\lambda_-,\lambda_-), & &\text{ with } \lambda_\pm=\frac3{2}\frac{(2-\delta)c_\pm^2 +\delta-\frac1\delta-(1-\gamma)}{\Theta_\pm(c_+^2-c_-^2)} ,\\
 (\mu_i)_{i\in\{1\dots4\}}&= (\mu_+,-\mu_+,\mu_-,-\mu_-), & &\text{ with } \mu_\pm=\frac{c\pm}{6}\frac{(1+\frac{3\gamma}\delta+\frac1{\delta^2})(c_\pm^2-\frac{1-\gamma}{\delta+1})-\frac1\delta c_\pm^2 }{c_\pm^2-2\frac{1-\gamma}{\delta+1}} ,
\end{align*}
with $\Theta_\pm\equiv\sqrt{2\delta(c_+^2-c_-^2)|c_\pm^2-1|}$.
 The unitary eigenvectors of $S_0^{-1}\Sigma_0$ are given by
\[\e_1=\frac{1}{\Theta_+}\begin{pmatrix} \frac{1}{c_+} \\ c_+-\frac{1}{c_+} \\ 1 \\ \delta c_+^2-\delta \end{pmatrix} ,\ 
\e_2=\frac{1}{\Theta_+}\begin{pmatrix} -\frac{1}{c_+} \\ \frac{1}{c_+}-c_+ \\ 1 \\ \delta c_+^2-\delta \end{pmatrix},\ 
\e_3=\frac{1}{\Theta_-}\begin{pmatrix} \frac{-1}{c_-} \\ \frac{1}{c_-}-c_- \\ -1 \\ \delta-\delta c_-^2 \end{pmatrix} ,\ 
\e_4=\frac{1}{\Theta_-}\begin{pmatrix} \frac{1}{c_-} \\ c_--\frac{1}{c_-} \\ -1 \\ \delta-\delta c_-^2 \end{pmatrix}.\]

Let us first remark that these coefficients only depend on $\gamma$ and $\delta$, so that neither the change of variables~\eqref{eq:newvar} nor the BBM-trick~\eqref{eq:BBMtrick} affect the coefficients of the KdV approximation.

\medskip

 It is also worth pointing out that the dispersion coefficients $\mu_\pm$ cannot be zeros for $\gamma\in(0,1)$ and $\delta>0$. Therefore, the KdV approximation cannot degenerate into Burgers-type equations.
\end{rmrk}
The coefficients of the KdV equation as an approximation to describe solitary internal waves in the free surface case have already been introduced in~\cite{KakutaniYamasaki78}, and then in~\cite{LeoneSegurHammack82,Matsuno93,MichalletBarth'elemy98}, and correspond to the ones we present here. In a slightly different regime, Craig, Guyenne and Kalisch~\cite{CraigGuyenneKalisch05} obtained a model that consists in four independent propagating waves, two of them satisfying the KdV equations of our slow mode waves (the other two being solutions of the first order transport equation $\partial_t u \pm c_+\partial_x u=0$). We state here that in the long wave regime, \emph{any} bounded solution of the full Euler system can be decomposed into four propagating waves, each of them being well approximated by independent solutions of four KdV equations. This simultaneous decomposition of the flow is new, and allows to compare with other models, such as Boussinesq-type models. %In particular, we are able to deduce the behavior (polarity, magnitude, thickness) of the solitary waves predicted by the KdV approximation in Section~\ref{sSec:AnalysisCoefs}, depending on the parameters $\gamma$ and $\delta$, for both rigid lid and free surface configurations. We offer in the following section a rigorous justification of the KdV approximation.

\subsection{Rigorous demonstration}\label{sSec:RigorousDemonstration}
The strategy is the following. Using the previous calculations, we define the approximate solution:
\begin{dfntn}\label{Def:Uapp}
 Let $U^0\in H^{s+2}(\RR)$, with $s>1/2$. We call {\em approximate solution} of the system~\eqref{eq:BoussForKdV} any function ${U_{\text{app}}(t,x)\in C^0([0,T/\epsilon);H^s)^4}$ such that
\begin{equation}  
  U_{\text{app}}(t,x)=\sum_{i=1}^4 u_i(t,x) \e_i +\epsilon U_1(\epsilon t,t,x),
\end{equation}
where $(u_i)_{i=1\dots4}$ is the solution of the four uncoupled KdV equations~\eqref{eq:KdV}, with ${u_i}\id{t=0}=\e_i\cdot S_0 U^0$, and the correcting term $U_1$ is a solution of~\eqref{eqnrest} with ${U_1}\id{\tau=t=0}=0$.
\end{dfntn}
 We first prove that such solutions exist in the strong sense for sufficiently smooth initial data, over times of order $\O(1/\epsilon)$. Then we use estimates on $U_0$ and $U_1$, to obtain a consistency result. This result allows us to show that $U_{\text{app}}(t,x)$ indeed approximates the solution of~\eqref{eq:BoussForKdV} with the same initial data, at least at order $\O(\epsilon^{3/2} t)$.

\begin{prpstn}\label{Prop:ExistenceUapp} 
 Let $U^0\in H^{s+2}(\RR)$, with $s>1/2$. Then, one has: \begin{enumerate}
                                           \item                             
For all $i=1\dots4$, there exists a unique strong solution of the Cauchy problem~\eqref{KdV}, with initial data ${{u_i}\id{t=0}=\e_i\cdot S_0 U^0}$. Moreover, there exists $T>0$, such that one has the estimate
\[ \sum_{i=1}^4 \big\vert u_i\big\vert_{L^\infty([0,T];H^{s+2})} \leq C_0\left(\frac{1}{\gamma(1-\gamma)},\delta+\frac1\delta\right) \big\vert U^0 \big\vert_{H^{s+2}}.\]
\item There exists a function $U_1\in L^\infty([0,T]\times\RR;H^{s})$, strong solution of~\eqref{eqnrest}, with ${U_1}\id{\tau=t=0}=0$.
    \end{enumerate}
\end{prpstn}
\begin{proof}
\begin{enumerate}\item         
 The existence and uniqueness of the solutions of~\eqref{KdV} is classical: see~\cite{KenigPonceVega93} for the local well-posedness of the KdV equation, and~\cite{CollianderKeelStaffilaniEtAl03} for the global well-posedness.
 One obtains the estimate as usual: as we multiply the equation~\eqref{KdV} by $\Lambda^{2k} u_i$ (with $3/2<k\leq s+2$) and integrate with respect to the space variable, one obtains 
 \[\frac12\frac{d}{dt}\int_\RR(\Lambda^k u_i)^2 dx= \left|\lambda_i \int_\RR \Lambda^k (u_i\partial_x u_i)\Lambda^k u_i dx\right|.
 \]
Thanks to the Kato-Ponce Lemma, we manage to estimate the right-hand side as follows: 
 \[\left|\int_\RR \Lambda^k (u_i\partial_x u_i)\Lambda^k u_i dx \right|\leq \left|\frac12\int_\RR  \partial_x u_i(\Lambda^k u_i )^2\right| + \left| \int_\RR  [\Lambda^k,u_i]\partial_x u_i(\Lambda^k u_i )dx\right| \leq C_0 \big\vert u_i\big\vert_{H^k}^3,\]
and we conclude by applying Gronwall-Bihari's lemma, that reads
 \[\big\vert u_i\big\vert_{H^k}\leq C_0 \frac{\big\vert{u_i}\id{t=0}\big\vert_{H^k}}{1-C_0 t\big\vert{u_i}\id{t=0}\big\vert_{H^k}}, \]
so that the estimate of the proof follows for $T$ sufficiently small, since \[\big\vert{u_i}\id{t=0}\big\vert_{H^{s+2}}=\big\vert \e_i\cdot S_0 U^0\big\vert_{H^{s+2}}\leq C_0\left(\frac{1}{\gamma(1-\gamma)},\delta+\frac1\delta\right) \big\vert U^0 \big\vert_{H^{s+2}}.\]

Let us recall the notation~\eqref{eq:transport}: we have thus proved that $u_i(\tau,t,x)\in L^\infty([0,T]\times\RR;H^{s+2})$, with
\[ \sum_{i=1}^4 \big\vert u_i\big\vert_{L^\infty([0,T]\times\RR;H^{s+2})} \leq C_0\left(\frac{1}{\gamma(1-\gamma)},\delta+\frac1\delta\right) \big\vert U^0 \big\vert_{H^{s+2}}.\]
\medskip

\item We can then exhibit $U_1$: let us write~\eqref{eqnrest} under the form \[(\partial_t+c_i\partial_x)\langle\e_i, U_1\rangle=\sum_{(j,k)\neq (i,i)}f_{ijk}(\tau,t,x)+\sum_{j\neq i}\partial_x g_{ij}(\tau,x-c_j t).\] From the above estimate on $u_i$, one has $f_{ijk}\in L^\infty([0,T]\times\RR;H^{s+1})$, and $g_{ij}\in L^\infty([0,T];H^{s})$. Hence, for $s>1/2$, one can set
\[\langle\e_i, U_1\rangle(\tau,t,x)=\sum_{(j,k)\neq (i,i)}\int_0^t f_{ijk}(\tau,s,x+c_i(s-t)) ds+\sum_{j\neq i} \frac1{c_i-c_j}\big(g_{ij}(\tau,x-c_j t)-g_{ij}(\tau,x-c_i t)\big),\]
and $U_1$ satisfies the hypotheses of the Proposition: $U_1(\tau,t,x) \in L^\infty([0,T]\times\RR;H^{s})$ and ${U_1}\id{\tau=t=0}=0$.
\end{enumerate}\end{proof}

 We now prove that $U_1$, which is the corrector term defined by~\eqref{eqnrest}, and that contains all the coupling effects between the different components, obeys to a sublinear secular growth. The key point is given by the following Lemma, that proceeds from Propositions 3.2 and 3.5 of~\cite{Lannes03}: 
\begin{lmm} \label{LannesLemma}Let $u$ be the solution of 
\begin{equation}
 \left\{\begin{array}{l}
       (\partial_t+c\partial_x) u=g(v_1,v_2)\\
       u\id{t=0}=0
        \end{array}\right. \ \text{ with }\ \quad\forall i\in\{1,2\},  \quad
 \left\{\begin{array}{l}
       (\partial_t+c_i\partial_x) v_i=0\\
       {v_i}\id{t=0}=v^0_i
        \end{array}\right.
\end{equation}
with $c_1\neq c_2$, $v^0_1$, $v^0_2 \in H^{s}(\RR)$, ${s}>1/2$, and $g$ is a bilinear mapping defined on $\RR^2$ and with values in $\RR$. Then one has the following estimates:
\begin{enumerate}
\item If $c=c_1$ , then $\lim\limits_{t \to \infty }\frac{1}{\sqrt t}\big\vert u(t,\cdot) \big\vert_{H^{s}(\RR)}=0$.
\item If $c\neq c_1\neq c_2$, then $\frac{1}{\sqrt t}\big\vert u(t,\cdot) \big\vert_{H^{s}(\RR)}=\O(1)$.
\end{enumerate}
Moreover, if there exists $\alpha>1/2$ such that $v^0_1 (1+x^2)^\alpha$,  and $v^0_2 (1+x^2)^\alpha\in H^{s}(\RR)$, then one has the better estimate 
\[ \big\vert u \big\vert_{L^\infty H^{s}(\RR)}\ \leq\ C_0 \big\vert v^0_1 (1+x^2)^\alpha \big\vert_{H^{s}(\RR)}\big\vert v^0_2 (1+x^2)^\alpha \big\vert_{H^{s}(\RR)},\]
with $C_0=C_0(c,c_1,c_2) $.
\end{lmm}
We are now able to give the following crucial estimate on $U_1$:
\begin{prpstn}
\label{EstU1}
 Let $s>1/2$ and $U^0\in H^{s+2}$. Then with $U_1\in L^\infty([0,T]\times\RR;H^{s})$ a strong solution of~\eqref{eqnrest} with ${U_1}\id{\tau=t=0}=0$, one has the estimate:
 \[\big\vert U_1 \big\vert_{L^\infty([0,T]\times[0,t];H^{s})} \leq  C_0  \sqrt t ,\]
 with $C_0=C_0 \big(\frac1{\gamma(1-\gamma)},\delta+\frac1\delta,\big\vert U^0\big\vert_{H^{s+2}}\big)$.
 
 Moreover, if $U^0$ satisfies $U^0 (1+x^2) \in H^{s+1}(\RR)$, then one has the uniform estimate 
 \[\big\vert U_1 \big\vert_{L^\infty([0,T]\times\RR;H^{s})} \leq  C_0 \big\vert U^0(1+x^2)\big\vert_{H^{s+1}}^2,\]
  with $C_0=C_0 \big(\frac1{\gamma(1-\gamma)},\delta+\frac1\delta,\big\vert U^0\big\vert_{H^{s+2}}\big)$.
 
\end{prpstn}
\begin{proof}
Let us decompose $U_1$ as a sum of functions as in Proposition~\ref{Prop:ExistenceUapp}: 
\begin{align*}\langle\e_i, U_1\rangle(\tau,t,x)&=\sum_{(j,k)\neq (i,i)}\int_0^t f_{ijk}(\tau,s,x+c_i(s-t)) ds+\sum_{j\neq i} \frac1{c_i-c_j}\big(g_{ij}(\tau,x-c_j t)-(g_{ij}(\tau,x-c_i t)\big)\\
 &= \sum_{(j,k)\neq (i,i)} U^{ijk} \ + \ \sum_{j\neq i} V^{ij}.
\end{align*}
with the functions $f_{ijk}$ and $g_{ij}$ coming from~\eqref{eqnrest} written in a simplified form:
\[(\partial_t+c_i\partial_x)\langle\e_i, U_1\rangle=\sum_{(j,k)\neq (i,i)}f_{ijk}(\tau,t,x) \ + \ \sum_{j\neq i}\partial_x g_{ij}(\tau,x-c_j t).\]
From Proposition~\ref{Prop:ExistenceUapp}, we know that the following bounds hold: 
\[\forall \tau \in [0,T], \ \sum_{(j,k)\neq (i,i)} \big\vert f_{ijk}\big\vert_{L^\infty([0,T]\times\RR;H^{s+1})} \ + \ \sum_{j\neq i} \big\vert g_{ij}\big\vert_{L^\infty([0,T];H^{s})} \leq C_0 \big\vert U^0\big\vert_{H^{s+2}},\]
with $C_0=C_0(\frac1{\gamma(1-\gamma)},\delta+\frac1\delta)$. Therefore, one has
\[\forall j\neq i,\qquad  \big\vert  V^{ij}\big\vert_{L^\infty([0,T]\times\RR;H^{s})}\ \leq \ C_0(\frac1{\gamma(1-\gamma)},\delta+\frac1\delta)\big\vert U^0\big\vert_{H^{s+2}}.\] Moreover, again for $j\neq i$, we remark that $f_{ijj}$ can be written as 
\[f_{ijj}(\tau,t,x) \ \equiv \ \alpha_{ijj}u_j(\tau,x-c_j t)\partial_x u_j(\tau,x-c_j t) \ \equiv \ \partial_x h_{ij}(\tau,x-c_j t), \]
so that $U^{ijj}$ has the same form as $V^{ij}$, and can be treated in the same way. And since $f_{ijj}\in L^\infty([0,T]\times\RR;H^{s+1})$, $U^{ijj}$ is uniformly bounded in $H^{s}$. The last terms that have to be bounded are $U^{ijk}$, for $(j,k)\neq (i,i)$ with $j\neq k$, $U^{ijk}$. One can easily check that $U^{ijk}$ satisfies the hypothesis of Lemma~\ref{LannesLemma}, with $f_{ijk}=g(u_j,\partial_x u_k)$, for any $\tau\in[0,T]$. 
We then immediately deduce
\[\big\vert U_1\big\vert_{L^\infty([0,T]\times[0,t];H^{s})} \leq\sqrt{t} C_0\Big(\frac1{\gamma(1-\gamma)},\delta+\frac1\delta,\big\vert U^0\big\vert_{H^{s+2}}\Big).\]

 As for the second estimate of the proposition, let us first remark that the estimates of $V^{ij}$ and $U^{ijj}$ are time-independant, and in agreement with the improved estimate. Therefore, the only remaining terms we have to control are $U^{ijk}$ with $j\neq k$. Of course, we will use the second case of Lemma~\ref{LannesLemma}, but we have to check first that for every $\tau\in[0,T]$, the initial data $u_j(\tau,0,x)$ and $\partial_x u_k(\tau,0,x)$ are localized in space, that is   
 \[\forall \tau\in[0,T],\qquad \big\vert  (1+x^2) {u_j}(\tau,0,x) \big\vert_{H^{s}} \ + \ \big\vert  (1+x^2) \partial_x{u_k}(\tau,0,x)\big\vert_{H^{s}} \ < \ \infty.\]
 This property is true at $\tau=0$ (by hypothesis of the proposition), and is propagated to $\tau>0$, using the fact that $u_i(\tau,x_i)$ satisfies the KdV equation~\eqref{KdV}. This propagation of the localization in space has been proved by Schneider and Wayne in~\cite[Lemma 6.4]{SchneiderWayne00}. We do not recall the proof here, and use directly the statement: 
 \begin{lmm}
  If $(1+x^2) U^0\id{t=0} \in H^{s+1}$, then there exists $C_1,\t C_1>0$ such that
 \[ \big\vert  (1+x^2) {u_j}(\tau,0,x) \big\vert_{H^{s+1}} \ \leq \ C_1\big\vert  (1+x^2) {u_j}\id{\tau=t=0} \big\vert_{H^{s+1}} \ \leq \  \t C_1\big\vert  (1+x^2) {U^0}\id{t=0} \big\vert_{H^{s+1}} .\]
 \end{lmm}
This Lemma, together with the second estimate of Lemma~\ref{LannesLemma}, allows to control $U^{ijk}$, uniformly in time.
Every term of the decomposition of $U_1$ has been controlled, and one has the following estimate:
\begin{equation}\label{eqn:estU12}
\big\vert U_1\big\vert_{L^\infty([0,T]\times\RR;H^{s})} \leq\ C_0\Big(\frac1{\gamma(1-\gamma)},\delta+\frac1\delta,\big\vert U^0\big\vert_{H^{s+2}}\Big)\big\vert (1+x^2) U^0\big\vert_{H^{s+1}}^2.
 \end{equation}
 This concludes the proof. \end{proof}

The next step consists in proving the consistency of our approximation with the symmetric system~\eqref{eq:BoussForKdV}.
\begin{prpstn}
\label{Prop:ConsistencyKdVBouss}
 If $U^0\in H^{s+5}$ with $s>1/2$, then $U_{\text{app}}(t,x)$ defined in Definition~\ref{Def:Uapp} satisfies the symmetric system~\eqref{eq:BoussForKdV} up to a residual of order $\O(\epsilon^{3/2})$ in $L^\infty([0,T/\epsilon];H^{s})$. 
 
 Moreover, $U_0$ satisfies $U^0 (1+x^2) \in H^{s+4}(\RR)$, then the residual is uniformly bounded in $L^\infty([0,T/\epsilon];H^{s})$ by $\epsilon^{2} C_0$, with $C_0=C_0\left(\frac{1}{\gamma(1-\gamma)},\delta+\frac1\delta,\big\vert (1+x^2) U^0\big\vert_{H^{s+4}},T\right)$. 
\end{prpstn}
\begin{proof}
 Plugging $U_{\text{app}}(t,x)$ into~\eqref{eq:AnsatzInBouss}, we see from the calculations of Section~\ref{sSec:FormalDerivation} that the only remaining term we have to control is $\epsilon^2 R(\epsilon t,t,x)$, with 
\[\begin{array}{r}R\equiv \partial_\tau U_1+\Sigma_1(U_0)\partial_x U_1+\Sigma_1(U_1)\partial_x U_0+S_1(U_0)\partial_t U_1+S_1(U_1)\partial_t U_0-\Sigma_2\partial_x^3 U_1-S_2\partial_x^2\partial_t U_1 \\  +\epsilon \Sigma_1(U_1)\partial_x U_1+\epsilon S_1(U_1)\partial_x U_1,\end{array}\]
where $U_0(\epsilon t,t,x)=\sum_{i=1}^4 u_i(t,x)\e_i$. 

\medskip

Each term of the right hand side is suitably bounded in the Sobolev $H^{s}$-norm, as we show in the following. Indeed, from Proposition~\ref{EstU1}, one has immediately
\[\big\vert \Sigma_2 \partial_x^3 U_1(\epsilon t,t,\cdot) \big\vert_{H^{s}} \leq C_0 \big\vert U_1(\epsilon t,t,\cdot) \big\vert_{H^{s+3}} \leq \sqrt{t} C_0 \left(\frac{1}{\gamma(1-\gamma)},\delta+\frac1\delta,\big\vert U^0\big\vert_{H^{s+5}},T\right).\]
Then, from~\eqref{eqnrest}, we deduce \[\langle\e_i,\partial_t U_1\rangle=-c_i\langle\e_i,\partial_x U_1\rangle+f_i\] with ${f_i\in L^\infty([0,T]\times\RR;H^{s+2})}$, and ${\big\vert f_i\big\vert_{H^{s+2}}\leq C_0 \big\vert U^0\big\vert_{H^{s+5}}}$, so that one has identically
\[\big\vert S_2 \partial_x^2\partial_t U_1(\epsilon t,t,\cdot) \big\vert_{H^{s}} \leq C_0 \big\vert \partial_t U_1(\epsilon t,t,\cdot) \big\vert_{H^{s+2}} \leq \sqrt{t} C_0 \left(\frac{1}{\gamma(1-\gamma)},\delta+\frac1\delta,\big\vert U^0\big\vert_{H^{s+5}}\right).\]

One obtains in the same way the desired estimates for $\Sigma_1(U_0)\partial_x U_1$, $\Sigma_1(U_1)\partial_x U_0$, $S_1(U_0)\partial_t U_1$, $S_1(U_1)\partial_t U_0$, $\Sigma_1(U_1)\partial_x U_1$ and $S_1(U_1)\partial_x U_1$. 

Finally, in order to estimate $\partial_\tau U_1$, we differentiate~\eqref{eqnrest} with respect to $\tau$. Since $u_i$ satisfies~\eqref{eq:KdV}, one has $\partial_\tau u_i \in L^\infty([0,T);H^{s+2})$. We are on the frame of the Lemma~\ref{LannesLemma}, so that we can obtain as in Proposition~\ref{EstU1} that $\partial_\tau U_1 \in L^\infty([0,T]\times[0,t];H^{s})$, with 
\[\big\vert \partial_\tau U_1(\epsilon t,t,\cdot)\big\vert_{H^{s}} \leq C_0  \sqrt{t}\big\vert  \partial_t U^0 \big\vert_{H^{s+2}}\leq C_0  \sqrt{t}\big\vert U^0 \big\vert_{H^{s+5}}.\]
Hence, $R(\epsilon t,t,\cdot)\in L^\infty ([0,T/\epsilon] ;H^{s})$, and
\[\big\vert R \big\vert_{H^{s}}\leq C_0  \sqrt{T/\epsilon}\big\vert U^0 \big\vert_{H^{s+5}} ,\] which concludes the first part of the proof.

\medskip

The second part follows in the exact same way, using the second estimate of Proposition~\ref{EstU1}.
\end{proof}

Finally, thanks to the consistency result and the estimate on $U_1$, we are able to set the following convergence Proposition:
\begin{prpstn}\label{convergenceKdVBouss}
Let $U^0\in H^{s+5}$, $s>1/2$, $U_B\in L^\infty([0,T/\epsilon];H^{s+5})$ be a family of solutions of~\eqref{eq:BoussForKdV} with ${U_B}\id{t=0}=U^0$ and $U_{\text{app}}$ be defined by Definition~\ref{Def:Uapp}, with the same initial value. Then one has 
\[\big\vert U_{\text{app}}-U_B\big\vert_{L^\infty([0,t];H^{s+1}_\epsilon)}\leq C_0  \epsilon^{3/2} t,\]
with $C_0 =C_0 \big(\frac{1}{\gamma(1-\gamma)},\delta+\frac1\delta,\big|U^0\big|_{H^{s+5}},T\big)$.

Moreover, $U^0$ satisfies $U^0 (1+x^2) \in H^{s+4}(\RR)$, then one has the better estimate 
 \[\big\vert U_{\text{app}}-U_B\big\vert_{L^\infty([0,t];H^{s+1}_\epsilon)}\leq C_0  \epsilon^{2} t ,\]
  with $C_0 =C_0 \big(\frac{1}{\gamma(1-\gamma)},\delta+\frac1\delta,\big\vert U^0(1+x^2)\big\vert_{H^{s+4}},T\big)$.
\end{prpstn}
\begin{proof}
 Let us set $R^\epsilon\equiv  U_{\text{app}}-U_B$. Thanks to Proposition~\ref{Prop:ConsistencyKdVBouss}, we know that 
\begin{equation}\label{RBOUSS}\big(S_0-\epsilon S_2\partial_x^2+\epsilon S_1(U_{\text{app}} \big) \partial_t R^\epsilon +\big(\Sigma_0-\epsilon \Sigma_2\partial_x^2+\epsilon \Sigma_1(U_{\text{app}})\big)\partial_x R^\epsilon=\epsilon^{3/2} f+\epsilon\A+\epsilon\B,\end{equation}
with $\A=\partial_t S_1(U_{\text{app}})R^\epsilon-S_1(R^\epsilon)\partial_t U_B$, $\B=\partial_x \Sigma_1(U_{\text{app}}) R^\epsilon-\Sigma_1(R^\epsilon)\partial_x U_B$ and a function $f\in L^\infty H^{s}$.

Then, we can follow the same path as for the proof of Proposition~\ref{Prop:WPSBOUSS}, in Section~\ref{ProofUnicity} (see also the proof of Proposition~\ref{Prop:ConvBoussEuler}). We define the energy as
\[E_s(R^\epsilon)\equiv\frac12(S_0 \Lambda^s R^\epsilon, \Lambda^s R^\epsilon)+\frac{\epsilon}{2}(S_2 \partial_x \Lambda^s R^\epsilon,\partial_x \Lambda^s R^\epsilon)+\frac{\epsilon}{2}(S_1(U_{\text{app}}) \Lambda^s R^\epsilon, \Lambda^s R^\epsilon),\]
and the exact same calculations lead to the following inequality: 
\[\frac{d}{dt}E_s(R^\epsilon)\leq C_0  \epsilon E_s(R^\epsilon)+ C_0 \epsilon^{3/2}(E_s(R^\epsilon))^{1/2},\]
with $C_0=C_0(\frac{1}{\gamma(1-\gamma)},\delta+\frac1\delta,\big|U^0\big|_{H^{s+5}})$.

From Gronwall-Bihari's theorem, we get ${E_s(R^\epsilon) \leq C_0  \epsilon^{1/2} (e^{C_0 \epsilon t}-1)}$, and finally for $\epsilon t\leq T$ ,
\[\big\vert U_{\text{app}}-U_B\big\vert_{H^{s+1}_\epsilon}\leq C_0 E_s(R^\epsilon) \leq C_0(\frac{1}{\gamma(1-\gamma)},\delta+\frac1\delta,\big|U^0\big|_{H^{s+5}},T)  \epsilon^{3/2} t.\]

The second part of the proof follows in the same way, using the consistency at order $\O(\epsilon^2)$ of Proposition~\ref{Prop:ConsistencyKdVBouss}.
\end{proof}
\begin{rmrk}\label{rem:proof}
 From this Proposition, one can immediately deduce the convergence rate between the solution of the symmetric system~\eqref{eq:BoussForKdV}, and the KdV approximation~\eqref{eq:KdV}. Indeed, since we know from Proposition~\ref{EstU1} the growth of the correcting term $U_1$, and since we restrict ourselves to times $0\leq  t\leq T/\epsilon$, the convergence rate is of order $\O(\epsilon\sqrt t)$ in general, and of order $\O(\epsilon)$ if the initial data is sufficiently decreasing in space.
 
 In the same way, we obtain the convergence rate between bounded solutions of the full Euler system~\eqref{FullEuler}, and the KdV approximation~\eqref{eq:KdV}, using Proposition~\ref{Prop:ConvBoussEuler}. This result is stated rigorously in Theorem~\ref{Th:CVKdVEuler}. One sees that even in the case where the initial data is rapidly decreasing in space, the convergence rate of the symmetric Boussinesq/Boussinesq model~\eqref{SBOUSS} (namely $\O(\epsilon^2 t)$) is better than the one of the KdV approximation (namely $\O(\epsilon)$). This is due to the interaction between the traveling waves of different wave modes, that is captured by the Boussinesq/Boussinesq system, and not by the uncoupled KdV approximation, and which is of order $\O(\epsilon)$ for times of order $\O(1)$. The decreasing in space of the initial data allows this error to remain of order $\O(\epsilon)$ for times of order $\O(1/\epsilon)$.
 
 Numerical simulations for both the Boussinesq/Boussinesq models and the KdV approximation are computed in Section~\ref{sSec:NumericalResults}. In particular, the relationship between the convergence rate and the decreasing in space of the initial data is discussed and enhanced in Figure~\ref{fig:ErrorInTime}.
\end{rmrk}

\subsection{The models under the rigid lid assumption}
\label{sSec:RigidLid}
In this section, we formally recover models existing in the literature in the rigid lid configuration. Starting from our Boussinesq/Boussinesq model~\eqref{Bouss}, we recover the three-parameter family of rigid lid Boussinesq/Boussinesq systems presented in~\cite{BonaLannesSaut08}. One can then apply the method presented in the previous section, in order to obtain the KdV approximation in this case.

\medskip

The rigid lid models use the variables $(\zeta,v)$, where $\zeta$ is the interface deviation ($-\eta_1=\eta_2\equiv\zeta$), and $v$ is the shear velocity defined by
\[v\equiv(\partial_x\phi_2-\gamma\partial_x\phi_1)\id{z=\epsilon\zeta}.\]
Using the calculations in~\cite{Duchene10}, one has
\[v=u_2-\gamma u_1 -\epsilon \left(\frac\gamma6\partial_x^2 u_1 +\left(\frac1{3\delta^2}+\frac{\gamma}{2\delta}\right)\partial_x^2 u_2\right)+\O(\epsilon^2).\]
Then, adding the first two equations of~\eqref{Bouss} leads to $\partial_x(h_1 u_1+h_2 u_2)=0$, so that one has $u_1+\frac1\delta u_2=\O(\epsilon)$, and 
\[v=\frac{\delta+\gamma}\delta u_2+\O(\epsilon)=-(\delta+\gamma) u_1+\O(\epsilon).\]
Therefore, using a straightforward combination the equations of~\eqref{Bouss}, one checks that the system becomes
\begin{equation}\label{eq:BoussTS}
\left\{ \begin{array}{l}
\partial_t \zeta + \frac{1}{\delta+\gamma} \partial_x v + \epsilon\frac{\delta^2-\gamma}{(\gamma+\delta)^2}\partial_x (\zeta v)+\epsilon\frac{1+\gamma\delta}{3\delta(\gamma+\delta)^2}\partial_x^3 v=\O(\epsilon^2), \\
\partial_t v +(1-\gamma)\partial_x \zeta+\epsilon\frac{\delta^2-\gamma}{(\delta+\gamma)^2}v\partial_x v =\O(\epsilon^2).
\end{array} 
\right. 
\end{equation}
Finally, using BBM-tricks as in~\eqref{eq:BBMtrick}, and the change of variable $v_\beta=(1-\epsilon \beta \partial_x^2)^{-1}v$ (with $\beta\geq 0$), one obtains eventually the three-parameter family of rigid lid Boussinesq/Boussinesq systems presented in~\cite{BonaLannesSaut08}:
\begin{equation}\label{eq:BoussTS3}
\left\{ \begin{array}{l}
(1-\epsilon b\partial_x^2)\partial_t \zeta + \frac{1}{\delta+\gamma} \partial_x v_\beta +\epsilon \frac{\delta^2-\gamma}{(\gamma+\delta)^2}\partial_x (\zeta v_\beta)+\epsilon a\partial_x^3 v_\beta=\O(\epsilon^2), \\
(1-\epsilon d\partial_x^2)\partial_t v_\beta +(1-\gamma)\partial_x \zeta+\epsilon\frac{\delta^2-\gamma}{(\delta+\gamma)^2}v_\beta\partial_x v_\beta +\epsilon c(1-\gamma)\partial_x^3 \zeta=\O(\epsilon^2),
\end{array} 
\right. 
\end{equation}
with $a$, $b$, $c$ and $d$ set (with $\theta_1\geq 0$, $\theta_2\leq 1$, $\beta\geq0$) as
\[(\gamma+\delta)a=(1-\theta_1)\frac{1+\gamma\delta}{3\delta(\gamma+\delta)}-\beta,\quad b=\theta_1\frac{1+\gamma\delta}{3\delta(\gamma+\delta)}, \quad c=\beta\theta_2,\quad d=\beta(1-\theta_2).
\]

From this system, one can easily follow the path of Section~\ref{sSec:FormalDerivation}, and deduce the KdV approximation related to system~\eqref{eq:BoussTS3}. One would obtain a similar result as in Theorem~\ref{Th:CVKdVEuler}. Eventually, the KdV approximation consists in decomposing the approximate solution $(\zeta_{\text{KdV}},v_{\text{KdV}})$ as
\[(\zeta_{\text{KdV}},v_{\text{KdV}})=u_+\e_+ + u_-\e_-,\]
with $u_+$ and $u_-$ two solutions of independent Korteweg-de Vries equations, namely
\begin{equation}\label{eq:KdV2}
 \partial_t u_\pm\ \pm\ c\partial_x u_\pm\ +\ \epsilon \lambda_\pm u_\pm\partial_x u_\pm\ \pm\ \epsilon\mu \partial_x^3 u_\pm=0,
\end{equation}
with the following vectors and coefficients:
\[\e_\pm=\frac1{\sqrt2}\begin{pmatrix} \pm\frac1{\sqrt{1-\gamma}} \\ \sqrt{\gamma+\delta} \end{pmatrix},\ c= \sqrt\frac{1-\gamma}{\gamma+\delta},\ \lambda=\frac{1}{\sqrt{2(1-\gamma)}}\frac{3c}{2}\frac{\delta^2-\gamma}{\gamma+\delta},\ \mu= \frac{c}6 \frac{(1+\gamma\delta)}{\delta(\gamma+\delta)}.\]

In that way, when looking at the decomposition of deformation of the interface, the KdV approximation leads to two counter-propagating waves, that is to say that one can write $\eta=\eta_++\eta_-$, with $\eta_\pm$ solution of
\[
 \partial_t \eta_\pm \ \pm\ \sqrt{\dfrac{1-\gamma}{\gamma+\delta}}\partial_x \eta_\pm \ +\ \epsilon \dfrac{3c}{2}\dfrac{\delta^2-\gamma}{\gamma+\delta} \eta_\pm\partial_x \eta_\pm \ \pm\ \epsilon\dfrac{c}6\dfrac{(1+\gamma\delta)}{\delta(\gamma+\delta)} \partial_x^3 \eta_\pm\ =\ 0.
\]
We recover the classical KdV equations in the rigid lid configuration (see for example~\cite{DjordjevicRedekopp78,KoopButler81,MichalletBarth'elemy97,CraigGuyenneKalisch05}). Following Section~\ref{sSec:RigorousDemonstration}, one would obtain in the same way a rigorous justification for the KdV approximation, under the rigid lid assumption.

One sees that whereas the KdV approximation in the rigid lid case leads to a decomposition into two waves, the free surface configuration predicts the decomposition into four waves, each of them solution of a KdV equation with different velocities. This striking fact leads to think that the rigid lid assumption may induce a significant alteration of the behavior of the solutions, and thus cannot be considered as a harmless statement in all configurations.
This remark has already been addressed in~\cite{Keulegan53,MichalletBarth'elemy98}, but to our knowledge, in the absence of the exhaustive decomposition given in Theorem~\ref{Th:CVKdVEuler}, the analysis of the difference between the two configurations has never been extensively discussed.

The following Section is devoted to a detailed study of the differences between the rigid lid and free surface configurations, depending on the values of the density ratio $\gamma$ and the depth ratio $\delta$.

\subsection{Discussion}
 \label{sSec:AnalysisCoefs}
Let us recall here that the KdV approximation in the free surface configuration consists in decomposing $U=(\eta_1,\eta_2,u_1,u_2)$ as ${U\sim \sum_{i=1}^4 u_i\e_i}$, with $u_i$ satisfying the KdV equation
\[\partial_t u_i \ +\ c_i \partial_x u_i \ +\ \epsilon\lambda u_i\partial_x u_i\ +\ \epsilon\mu_i\partial_x^3 u_i \ = \ 0 .\]
The coefficients $c_i,\lambda_i$ and $\mu_i$, as well as the vectors $\e_i$ are given in Remark~\ref{rem:coefsKdV}, page \pageref{rem:coefsKdV}. This leads to the following decomposition for the respective deformations of the interface and the surface:
\begin{enumerate}
 \item For the interface: $\eta_2=\displaystyle\sum_{i=1}^4 u_i\e_{i,2}\equiv\sum_{(j,k)=(\pm,\pm)} \eta_{j,k}$, and $\eta_{\pm,\pm}$ satisfies the KdV equation
 \[\partial_t \eta_{\pm,k} \ \pm\ c_k \partial_x \eta_{\pm,k} \ +\ \epsilon\lambda^i_k \eta_{\pm,k}\partial_x \eta_{\pm,k}\ \pm\ \epsilon\mu_k\partial_x^3 \eta_{\pm,k}\ = \ 0,\]
 where $\e_{i,j}$ denotes the $j^{th}$ component of the vector $\e_i$ and with the following coefficients\footnote{These are the coefficients that are displayed in~\cite{KakutaniYamasaki78,LeoneSegurHammack82,MichalletBarth'elemy98}.}: \[\begin{array}{lr}
 \displaystyle c_\pm^2=\dfrac{1+\delta\pm\sqrt{(1-\delta)^2+4\gamma\delta}}{2\delta}, &\mu_\pm=\dfrac{c\pm}{6}\dfrac{(1+\frac{3\gamma}\delta+\frac1{\delta^2})(c_\pm^2-\frac{1-\gamma}{\delta+1})-\frac1\delta c_\pm^2 }{c_\pm^2-2\frac{1-\gamma}{\delta+1}},  \\                                                                                                                                                                            
          \multicolumn{2}{c}{\displaystyle \lambda^i_\pm=\dfrac{3c_\pm}{2}\frac{(2-\delta)c_\pm^2+\delta-\frac1\delta-(1-\gamma)}{|(c_+^2-c_-^2)(1-c_\pm^2)|}.}                                                                                                                                                                                  \end{array}\]
 \item For the surface: $\zeta_1=\eta_1+\eta_2=\displaystyle\sum_{i=1}^4 u_i(\e_{i,1}+\e_{i,2})\equiv\sum_{(j,k)=(\pm,\pm)} \zeta_{j,k}$, and $\zeta_{\pm,\pm}$ satisfies the KdV equation
 \[\partial_t \zeta_{\pm,k} \ \pm\ c_\pm \partial_x \zeta_{\pm,k}\ +\ \epsilon\lambda^s_k \zeta_{\pm,k}\partial_x \zeta_{\pm,k}\ \pm\ \epsilon\mu_k\partial_x^3 \zeta_{\pm,k}\ = \ 0,\]
 with the same values as previously for $c_\pm$ and $\mu_\pm$, and \[\lambda^s_\pm=\dfrac{3c_\pm}{2}\dfrac{(2-\delta)c_\pm^2 +\delta-\frac1\delta-(1-\gamma)}{(c_+^2-c_-^2)c_\pm^2}.\]
 
\end{enumerate}
\begin{figure}[ht]
 \centering
\subfigure[$\delta=1/2$]{\includegraphics [width=0.45\textwidth,keepaspectratio=true]{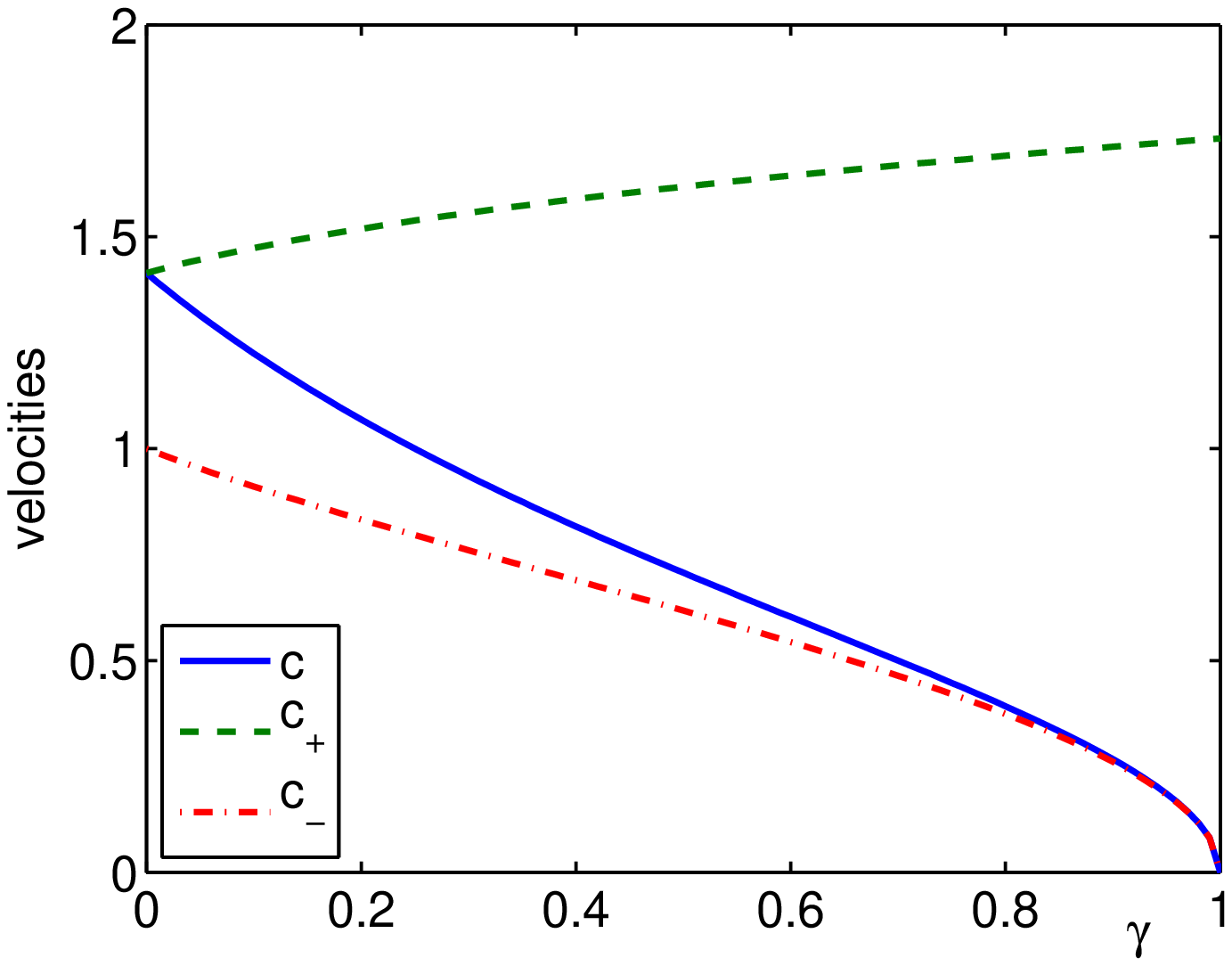}} 
\subfigure[$\delta=2$]{\includegraphics [width=0.45\textwidth,keepaspectratio=true]{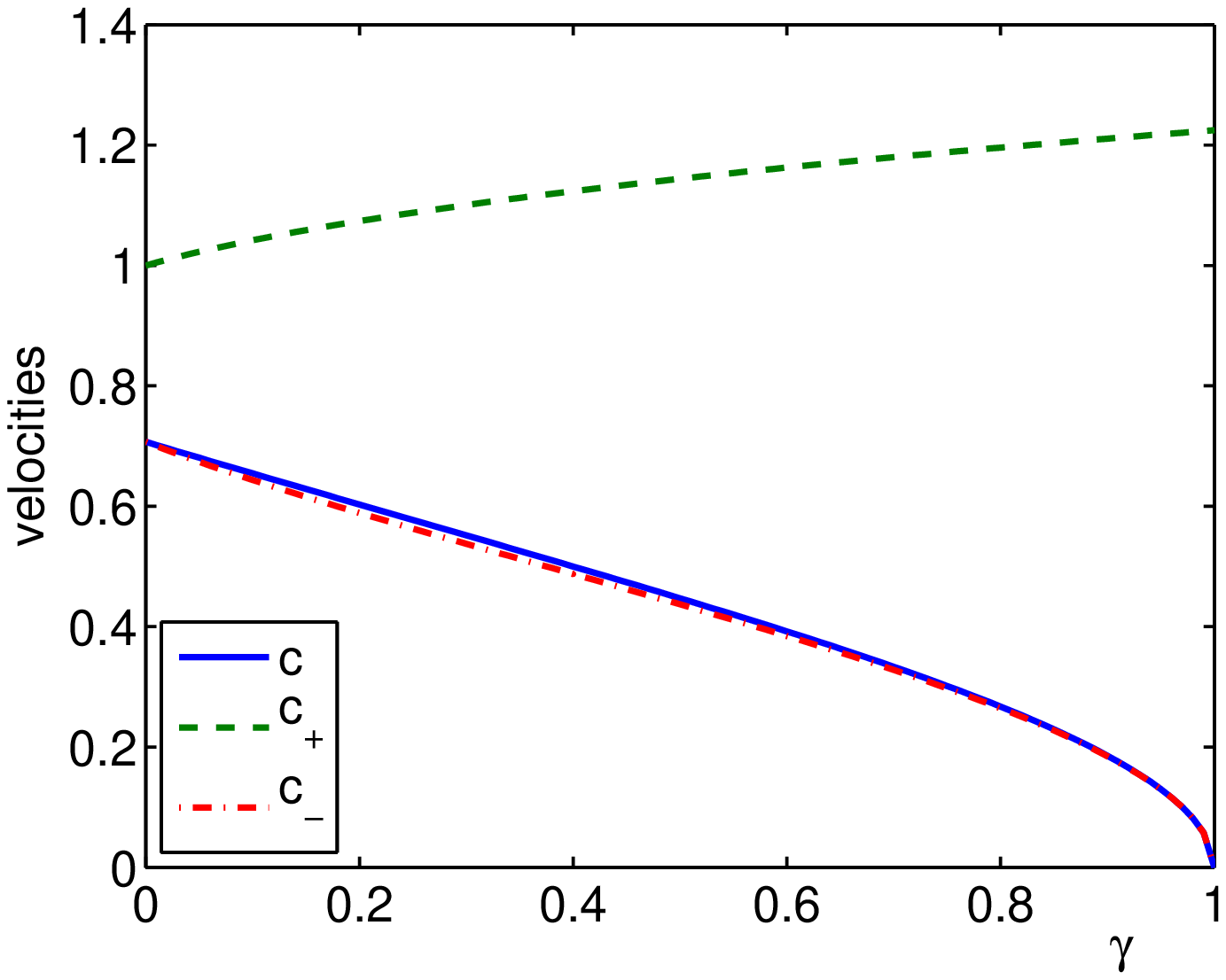}} 
 \caption{The different velocities, in the rigid lid configuration (${c= \sqrt\frac{1-\gamma}{\gamma+\delta}}$) and free surface  case (${c_\pm=\sqrt{\frac{1+\delta\pm\sqrt{(1-\delta)^2+4\gamma\delta}}{2\delta}}}$), for $\gamma\in(0,1)$ and (a) $\delta=1/2$,  (b) $\delta=2$.}
 \label{fig:velocity}
\end{figure}
Therefore, one sees that for both the surface and the interface elevations, the KdV approximation predicts the evolution of four different waves, two of them corresponding to the velocities $\pm c_+$ (we call them fast mode waves), and the other two corresponding to the velocities $\pm c_-$ (we call them slow mode waves), with $c_+>c_->0$. The fact that such different modes exist is characteristic of the free surface configuration, as only two counter-propagating waves appear in the rigid lid case. Moreover, as we see in Figure~\ref{fig:velocity}, the two velocities corresponding to the free surface configuration can be very different from the velocity in the rigid lid case (namely $c= \sqrt\frac{1-\gamma}{\gamma+\delta}$), depending on the values of $\delta$ and $\gamma$. In these cases, one expect the solutions in the two configurations to behave very differently.

The aim of this section is to study more in depth the behavior of the KdV approximation in the two different configurations, with respect to the parameters $\gamma$ and $\delta$. The first part is devoted to the study of solitary waves, as known solutions of the KdV equations. Then, we study the case where the initial data satisfy the rigid lid hypothesis, and explore the evolution the surface in that case. The results we obtain are summarized in Section~\ref{sSec:sum}; we let the reader refer to Figures~\ref{fig:solitonA}--\ref{fig:rigidlidD} for a numerical illustration of our statements.

\begin{rmrk}
 In the following study, and especially in Table~\ref{tab:rigidlidlimits}, we allow ourselves to look at the behavior of the system in the limit cases of the parameters ($\delta\to0,\infty$ and $\gamma\to0,1$), despite the fact that the rigorous justification of the KdV approximation, as well as the Boussinesq/Boussinesq models, break in these limits. However, the KdV approximation when $\gamma\to1$ has been widely used in the literature (see for example~\cite{GrimshawPelinovskyTalipova97,HelfrichMelville06,SakaiRedekopp07,SegurHammack82}), and such limits offer striking
 illustrations of our discussion.
\end{rmrk}

\subsubsection{Solitary waves}
 It is well known that the solitary wave solutions of the generic KdV equation
 \[ \partial_t u \ +\ c \partial_x u \ +\ \epsilon\lambda u\partial_x u\ +\ \epsilon\mu\partial_x^3 u \ = \ 0\]
can be expressed as follows:
  \begin{equation}\label{soliton}
 u(t,x)=\frac{M}{\cosh(k(x-x_0-c' t))^2},
\end{equation}
with $c'=c+\epsilon\frac{\lambda M}{3}$, $k=\sqrt\frac{\lambda M}{12 \mu}$, and $M$ and $x_0$ arbitrary. 

We discuss in the following the polarity, magnitude and thickness of such waves, as the parameters $\gamma$ and $\delta$ specify the coefficients of the KdV approximation.

\bigskip
\noindent{\bf Polarity.}
It is obvious that for $k=\sqrt\frac{\lambda M}{12 \mu}$ to be real valued, the sign of the ratio $\lambda/\mu$ determines the sign of the acceptable values of $M$.  Therefore, we are able to predict, depending on the parameters $\gamma$ and $\delta$, the polarity of the solitary waves predicted by the KdV approximation (elevation or depression). We give here the result, first for the free surface case, and then in the rigid lid configuration.

First, one can check that for every value of $\delta>0$ and $\gamma\in(0,1)$, the three coefficients $\lambda^i_+$, $\lambda^s_+$ and $\mu_+>0$ are positive. Hence, the fast mode solitary waves will always be of elevation/elevation type (both surface wave and interface wave are convex upward). The qualitative nature of the fast mode is thus similar to that of the one-layer water-wave problem, which is always of elevation type. In particular, when we set $\gamma\to0^+$ and $\delta\to1^-$, one recovers the classical KdV equation for a single layer at the interface ($c_+\to1$, $\lambda^i_+\to 3/2$, $\mu_+\to1/6$), and when we set $\delta\to\infty$, one recovers the classical KdV equation for a single layer at the surface ($c_+\to1$, $\lambda^s_+\to 3/2$, $\mu_+\to1/6$).

The behavior of the slow mode is more peculiar, as the nonlinear coefficients $\lambda^{s,i}_-$ can have both signs, depending on the size of the thickness ratio $\delta$. Indeed, as it has been pointed out in~\cite{Walker73,KakutaniYamasaki78} and then in~\cite{MichalletBarth'elemy98,OstrovskyStepanyants05}, for every $\gamma\in(0,1)$, there exists a critical ratio $\delta_c(\gamma)$ such that if $\delta>\delta_c(\gamma)$, then $\lambda^{i}_->0$ and $\lambda^{s}_-<0$ (and conversely if $\delta<\delta_c$). Since one has $\mu_->0$ for every value of $\delta>0$ and $\gamma\in(0,1)$, we know that the slow mode solitary wave will be of elevation/depression type (surface wave convex upward, and interface wave concave) if $\delta<\delta_c$, and of depression/elevation type if $\delta>\delta_c$. There is no solitary waves in the case $\delta=\delta_c$.

\begin{wrapfigure}{r}{0.5\textwidth}
 \centering 
             \includegraphics[width=0.5\textwidth,keepaspectratio=true]{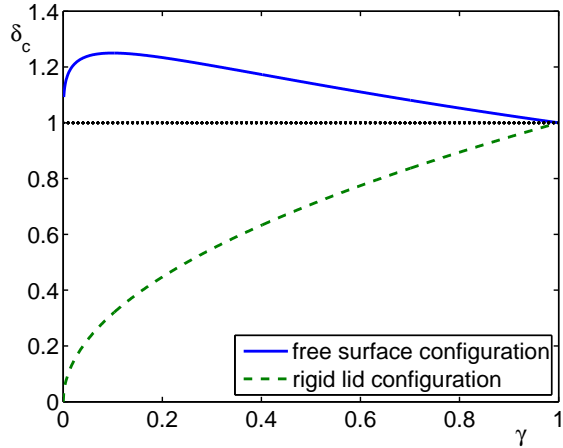} 
 \caption{Dependence of the critical ratio $\delta_c$ on the density ratio $\gamma$, for both the rigid lid and free surface configurations.}
 \label{fig:criticalratio}
\end{wrapfigure}

More precisely, the critical ratio is the unique real solution of the equation
\[X^3+(\gamma^2+3\gamma-3)X^2+(3-4\gamma)X-1=0,\]
and takes values in $\delta_c\in(1,5/4]$ for $\gamma\in(0,1)$ (see Figure~\ref{fig:criticalratio}).

In that way, the behavior of the slow mode waves resembles that of the interface waves with a rigid lid. Indeed, such a critical ratio appears straightforwardly in the rigid lid configuration: if $\delta<\delta_c'\equiv\sqrt\gamma$, then the interface wave is of depression type, and if $\delta>\delta_c'$, then the interface wave is of elevation type.
This critical ratio is well known in the literature, and has led to many extended models, where a cubic nonlinear term becomes the major source of nonlinearity for $\delta\sim\delta_c$ (see for example~\cite{DjordjevicRedekopp78,KoopButler81,FunakoshiOikawa86,Guyenne06,HelfrichMelville06}).

It is interesting to see that even if the respective polarity of the interface slow mode waves in the free surface configuration, and the interface waves in the rigid lid configuration follow qualitatively the same behavior, the value of the critical ratio is notably different. Indeed $\delta_c$ is always located above $1$ in the free surface case, on the contrary to the rigid lid case. Moreover, as we can see in Figure~\ref{fig:criticalratio}, they have considerably different values, except when $\gamma\sim 1$. In the area between the two curves, the polarities of the interface waves in the free surface and in the rigid lid configurations are reversed; the two models therefore lead to considerably different results.

\bigskip
\noindent{\bf Magnitude of the deformations.}
Depending of the parameters $\gamma$ and $\delta$, we are able to compare the magnitudes of the respective amplitudes of the surface and the interface waves. Indeed, since $\eta_2=\sum_{i=1}^4 u_i\e_{i,2}=\sum \eta_{\pm,\pm}$, and $\zeta_1=\eta_1+\eta_2=\sum_{i=1}^4 u_i(\e_{i,1}+\e_{i,2})=\sum \zeta_{\pm,\pm}$, one sees immediately that the surface and interface deformations, for each mode, are proportional, and satisfy 
\[\frac{\eta_{\pm,\pm}}{\zeta_{\pm,\pm}} = \frac{\e_{i,2}}{\e_{i,1}+\e_{i,2}}= \frac{c_\pm^2-1}{c_\pm^2}.\]
In that way, one deduces that for the fast mode, the surface deformation is always bigger than the interface deformation, and the ratio tends to zero when $\gamma\to0$ with $\delta\geq1$, or when $\delta\to\infty$.

Meanwhile, as remarked in~\cite{KakutaniYamasaki78}, the amplitude of the surface deformation is bigger than the one of the interface for the slow mode if $0<\delta\leq 2(1-2\gamma)$, and conversely if $\delta>2(1-2\gamma)$. Moreover, the ratio tends to zero when $\gamma\to0$ with $\delta\leq 1$, and tends to $\infty$ when $\gamma \to 1$ or $\delta\to \infty$.

\bigskip
\noindent{\bf Thickness.} In addition to forcing the polarity of the solitary wave, the ratio $\frac{\lambda M}{12 \mu}$ is also related to the thickness, or the wavelength of the wave. Indeed, defining the thickness of a wave as in~\cite{MichalletBarth'elemy98} by
\[l(u)\equiv\frac{1}{2M}\int_{-\infty}^{+\infty} u(x)dx,\]
one obtains for the function~\eqref{soliton}: $l(u)=\frac1k=\sqrt\frac{12 \mu}{\lambda M}$.

As an immediate result, and since we know the ratio between the magnitude of the deformations at the surface and the interface, the thickness of the deformations are identical at the surface and at the interface:
\[\frac{l_\pm^s}{l_\pm^i}=\sqrt{\frac{\lambda^{i}_\pm}{\lambda^{s}_\pm}\left| \frac{c_\pm^2-1}{c_\pm^2}\right|}=1.\]

\begin{figure}[ht]
 \centering
\subfigure[$\delta=1/2$]{\includegraphics [width=0.3\textwidth,keepaspectratio=true]{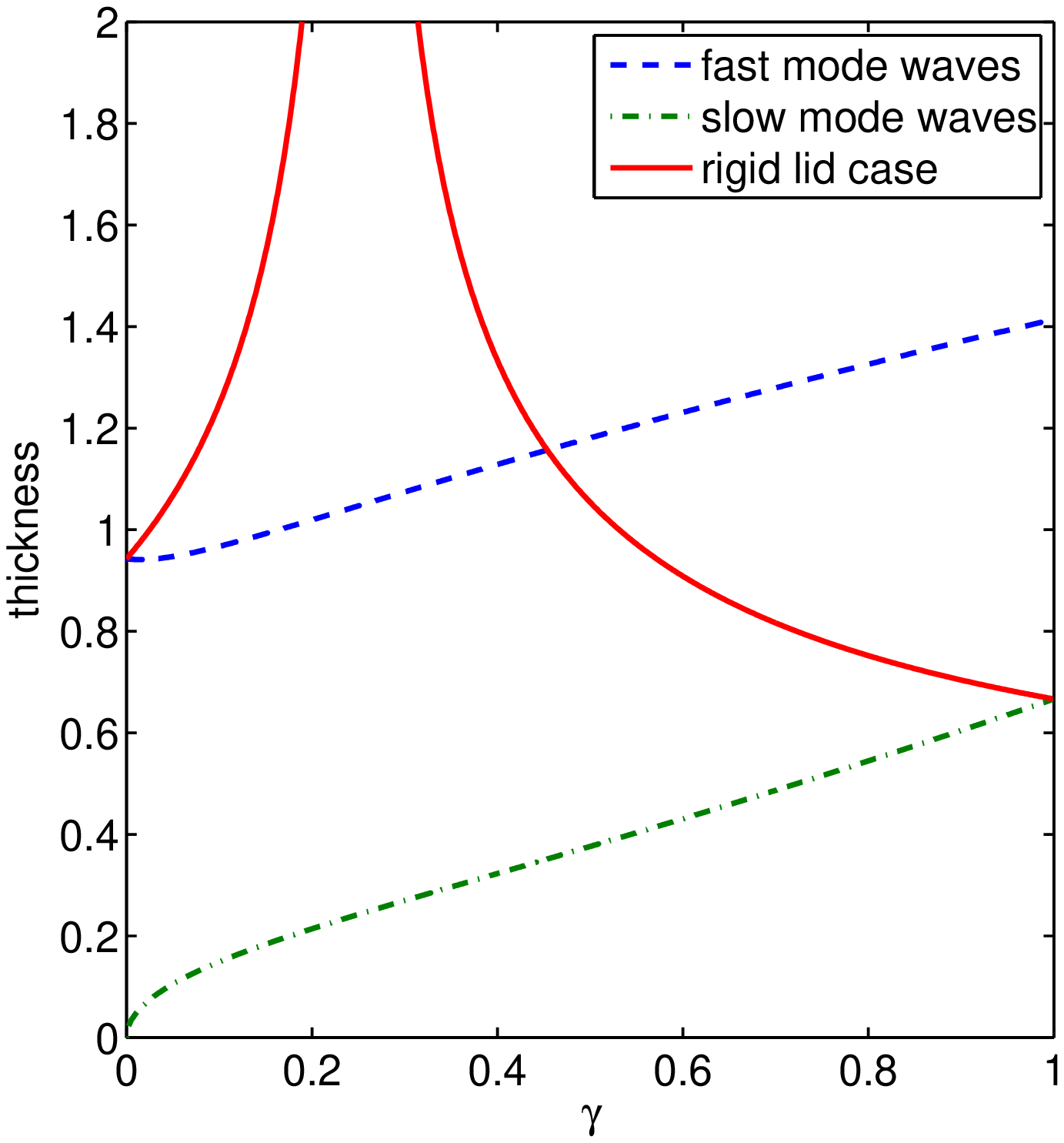}} 
\subfigure[$\delta=1$]{\includegraphics [width=0.3\textwidth,keepaspectratio=true]{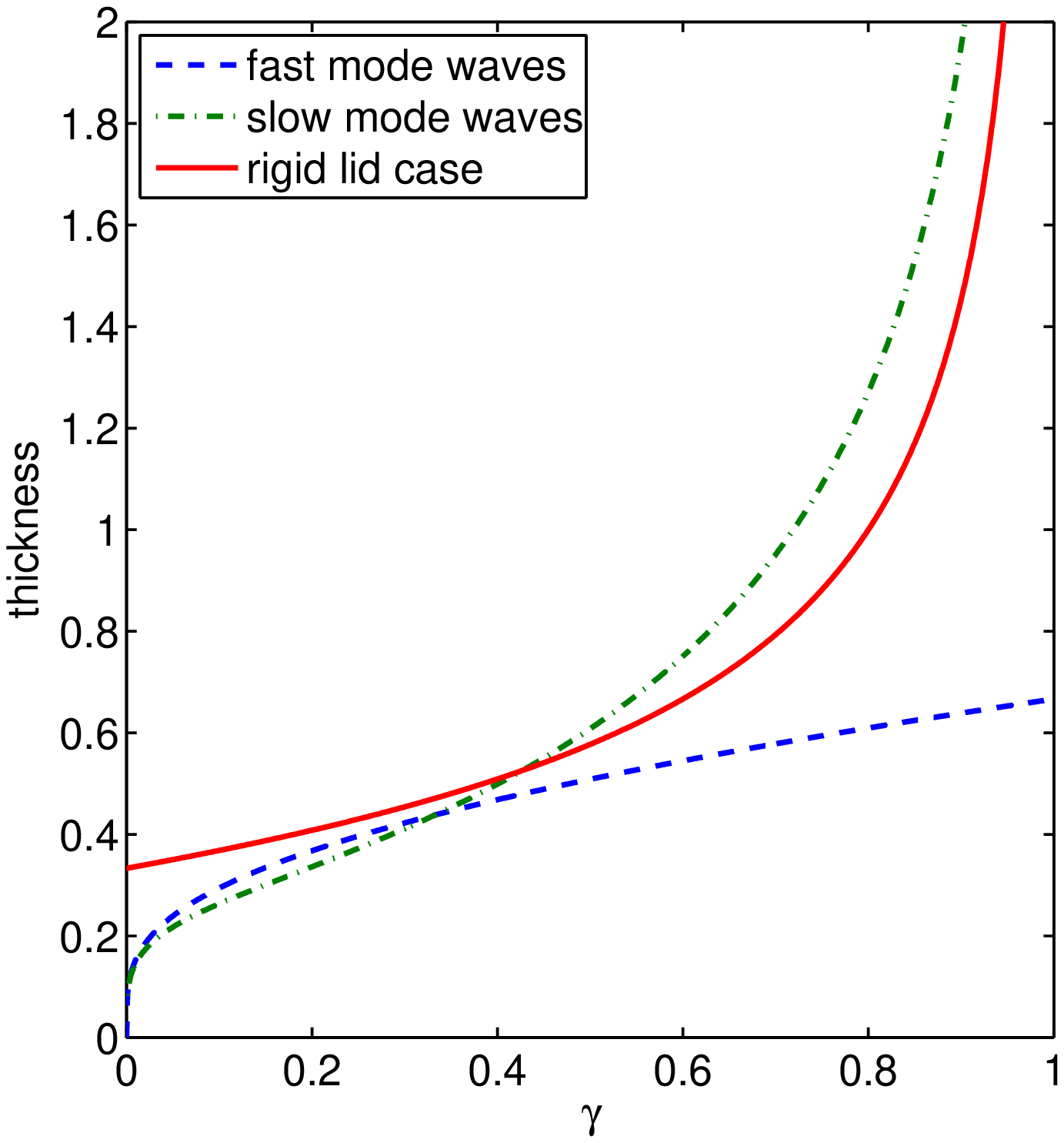}} 
\subfigure[$\delta=2$]{\includegraphics [width=0.3\textwidth,keepaspectratio=true]{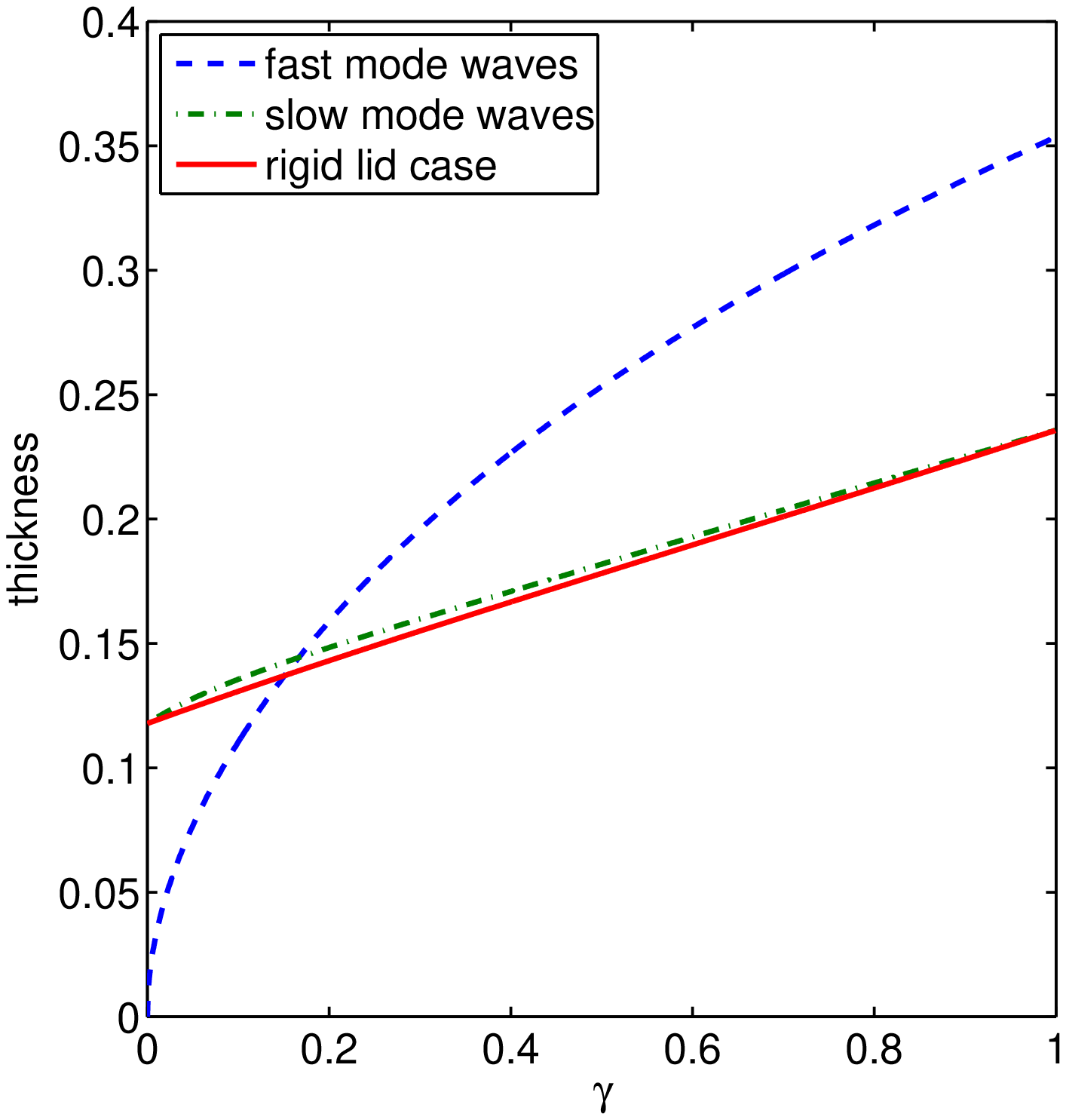}} 
 \caption{The different thicknesses, in the rigid lid and free surface configuration, for the parameters $\gamma\in(0,1)$ and (a) $\delta=1/2$,  (b) $\delta=1$,  (c) $\delta=2$.}
 \label{fig:thickness}
\end{figure}

One can now compare the thickness of the different wave modes, and the ones in the rigid lid configuration, for waves of same heights. One computes in Figure~\ref{fig:thickness} the ratio $\sqrt\frac{\mu_+}{\lambda^i_+}$, $\sqrt\frac{\mu_-}{\lambda^i_-}$ and $\sqrt\frac{\mu}{\lambda^i}$ for $\gamma \in (0,1)$, and $\delta=1/2,1,2$. 

Of course, the thickness of the solitary wave tends to infinity when the depth ratio approaches its critical ratio $\delta=\delta_c$, so that the waves predicted in the rigid lid and in the free surface configurations are excessively unlike. Additionally, one sees that when $\delta$ is small, then the thickness of the internal waves predicted in the free surface configuration is largely different from both the slow mode and fast mode thicknesses, except in the limit cases $\gamma\to 0$ and $\gamma\to 1$. 

When $\gamma\ll 1$ and $\delta\leq 1$, one has at the same time the similitude of the thickness and the velocities of the waves in the rigid lid configuration, and the fast mode waves in the free surface configuration. However, the amplitude of the surface deformation is bigger than the one of the interface for the fast mode waves, so that the rigid lid is not valid hypothesis. This resemblance has to do with the fact that both the fast mode waves and the waves in the rigid lid configuration converge to waves of the one-layer problem, when $\gamma\to 0$ and $\delta \to 1^-$. 
Conversely, when $\gamma\sim 1$ and when $\delta \gg 1$, the characteristics of the waves in the rigid lid configuration resemble the slow mode waves ones. We explore thereafter the validity of this hypothesis for different situations.

\subsubsection{Evolution of the rigid lid hypothesis}
Now, we restrict ourselves to initial data that are compatible for both the free surface and rigid lid configurations, and compare the different evolutions of the solutions. Following Section~\ref{sSec:RigidLid}, if the initial data with the rigid lid assumption is $(\eta^0,v^0)$, then the corresponding initial data in the free surface case  (in the limit $\epsilon\to0$) is $U^0=(-\eta^0,\eta^0,\frac{-1}{\gamma+\delta}v^0,\frac{\delta}{\gamma+\delta}v^0)$). Therefore, one has initially ${u_i}\id{t=0}=\e_i\cdot S_0 U^0$, so
\begin{align*}
\sum_{(j,k)=(\pm,\pm)} {\eta_{j,k}}\id{t=0}&=\sum_{i=1}^4\e_{i,2} {u_i}\id{t=0} =\e_{i,2}\left((1-\gamma)\e_{i,2}\eta^0+ \frac{(\e_{i,3}-\gamma\e_{i,4})v^0}{\gamma+\delta}\right),\\
\sum_{(j,k)=(\pm,\pm)} {\zeta_{j,k}}\id{t=0}&=(\e_{i,1}+\e_{i,2})\left((1-\gamma)\e_{i,2}\eta^0+ \frac{(\e_{i,3}-\gamma\e_{i,4})v^0}{\gamma+\delta}\right).
\end{align*}
Using the values given in Remark~\ref{rem:coefsKdV} page~\pageref{rem:coefsKdV}, this reads
\begin{align*}
\eta_{j,k}^0={\eta_{j,k}}\id{t=0}&=\frac{k}{2\delta(c_+^2-c_-^2)}\left((1-\gamma)\frac{c_k^2-1}{c_k^2}\eta^0+j \frac{1+\gamma\delta-\gamma\delta c_k^2}{(\gamma+\delta)c_k}v^0 \right),\\
\zeta_{j,k}^0={\zeta_{j,k}}\id{t=0}&=\frac{k}{2\delta(c_+^2-c_-^2)}\left((1-\gamma)\eta^0+j \frac{c_k^2}{c_k^2-1}\frac{1+\gamma\delta-\gamma\delta c_k^2}{(\gamma+\delta)c_k}v^0 \right).
\end{align*}
Since the KdV equation preserves mass, knowing the size of the initial data allows to predict the significance of the waves. In particular, the rigid lid hypothesis will be valid for long times only if $\vert\zeta_{\pm,\pm}^0\vert\ll \vert\eta_{\pm,\pm}^0\vert$. We give in Table~\ref{tab:rigidlidlimits} the different behavior of these variables, in the limits $\gamma\to 1$, $\gamma\to 0$, $\delta\to \infty$ and $\delta\to 0$. As we can see, the rigid lid hypothesis will be valid for long times when $\gamma \sim 1$, or when $\delta \gg 1$. For each of these cases, one sees that the main deformation comes from the slow mode waves, which correspond to the waves predicted by the models in the rigid lid configuration. 
\begin{table}[!ht]
 
{  \begin{tabular}{|c|c|c|c|c|c|} \hline
                   & $\gamma\to 1$ & $\gamma\to 0,\ \delta>1$ & $\gamma\to 0,\ \delta<1$ &  $\delta\to \infty$ & $\delta\to 0$ \\ \hline 
 $\eta_{\pm,+}^0$  & $0$ & $\dfrac{\pm\nu}{2\sqrt\delta} v^0$ & $\frac12\left(\eta^0 \pm \nu v^0 \right)$ & $0$ & $\frac12(1-\gamma)\eta^0$  \\ \hline 
 $\eta_{\pm,-}^0$ & $\frac12\eta^0 \pm \varsigma v^0$ & $\frac12\left(\eta^0 \mp \nu v^0 \right)$ &$\mp\frac12 \nu v^0 $  &  $\frac12\eta^0 $ & $\frac12\left(\gamma\eta^0 \mp \frac{1}{\gamma\sqrt{1-\gamma}}v^0\right) $  \\ \hline 
 $\zeta_{\pm,+}^0$ & $0$ & $\dfrac{1}{2(\delta-1)}\left(\eta^0 \pm \tau v^0\right)$ &$\dfrac{1}{2(1-\delta)}\left(\eta^0 \pm \nu v^0\right)$ & $0$ & $\frac12(1-\gamma)\eta^0$ \\ \hline 
 $\zeta_{\pm,-}^0$ & $0$ & $\dfrac{1}{2(\delta-1)}\left(\eta^0 \mp \nu v^0\right)$ &$\dfrac{1}{2(1-\delta)}\left(-\eta^0 \pm \tau v^0\right)$ & $0$ & $\dfrac{-1}{2}\left((1-\gamma)\eta^0\mp\frac{\sqrt{1-\gamma}}{\gamma^2}v^0\right)$\\ \hline
 
\end{tabular}} \\ with $\varsigma\ \mathop{\sim}\limits_{\gamma\to 1}\ \frac{1}{2\sqrt{(\delta+1)(1-\gamma})}$, $\tau\ \mathop{\sim}\limits_{\gamma\to0}\ \frac{\delta-1}{\delta\gamma}$ and $\nu=\frac1{\sqrt\delta (\delta-1)}$.
 \caption{Initial magnitudes of the different waves at the surface and at the interface, for an initial with a flat surface, in the limit cases.}
  \label{tab:rigidlidlimits}
 \end{table}

\medskip

As a specific example, when the initial data has zero velocities (that is to say $v^0=0$), then we are able to compare straightforwardly the different magnitudes of the four waves. Indeed, one deduces from the previous calculations that when $v^0=0$, one has
\[
\eta_{j,k}^0=\frac{k(1-\gamma)}{2\delta(c_+^2-c_-^2)}\frac{c_k^2-1}{c_k^2}\eta^0 ,\ \ {\rm and}\ \
\zeta_{j,k}^0=\frac{k(1-\gamma)}{2\delta(c_+^2-c_-^2)} \eta^0.
\]
Consequently, the four different waves have the same weight at the surface. The situation is more sophisticated at the interface, and one has eventually
\[\begin{array}{ll}
   |\eta_{\pm,+}|_{L^2}\ \geq\  |\eta_{\pm,-}|_{L^2}&{\rm if }\ \delta\ \leq\ 1-2\gamma, \\
   |\eta_{\pm,+}|_{L^2}\ <\ |\eta_{\pm,-}|_{L^2}&{\rm if }\ \delta\ >\ 1-2\gamma ,
  \end{array} 
\]
so that the fast mode wave is more significant than the slow mode wave when $\delta \leq 1-2\gamma$, and conversely otherwise. In the limits $\delta\to\infty$ and $\gamma \to 1$, the magnitudes of the fast mode waves tend to 0, so that the energy is only shared by the slow mode waves. Meanwhile, in the limit $\delta\to 0$, the significance of the fast mode and the slow modes respectively tend to $\frac{1-\gamma}2|\eta^0|_{L^2} $ and $\frac\gamma2|\eta^0|_{L^2}$. Finally, in the limit $\gamma \to 0$, then the magnitudes of the fast mode waves tend to 0 when $\delta>1$, and in the contrary carry all the energy when $\delta<1$.

\medskip

% Finally, in the case of initial data leading to solitary waves, when comparing the magnitudes of the deformation at the surface and at the interface, one has as 
% \[\begin{array}{ll}
%    |\zeta_{\pm,+}|_{L^2}\ >\  |\eta_{\pm,-}|_{L^2}&{\rm for\ all }\ \delta\ \in(0,+\infty), \\
%    |\zeta_{\pm,-}|_{L^2}\ \geq\  |\eta_{\pm,-}|_{L^2}&{\rm if }\ \delta\ \leq\ 2(1-2\gamma), \\
%    |\zeta_{\pm,-}|_{L^2}\ <\ |\eta_{\pm,-}|_{L^2}&{\rm if }\ \delta\ >\ 2(1-2\gamma) ,
%   \end{array}
% \]
% Thus for the fast mode, the surface elevation is always bigger than the interface elevation.
% As for the slow mode, the amplitude of the surface deformation is bigger than the one of the interface if $0<\delta\leq 2(1-2\gamma)$, and conversely if $\delta>2(1-2\gamma)$. 

\begin{figure}[ht]
 \centering
\subfigure[$\delta=1/2$]{\includegraphics [width=0.3\textwidth,keepaspectratio=true]{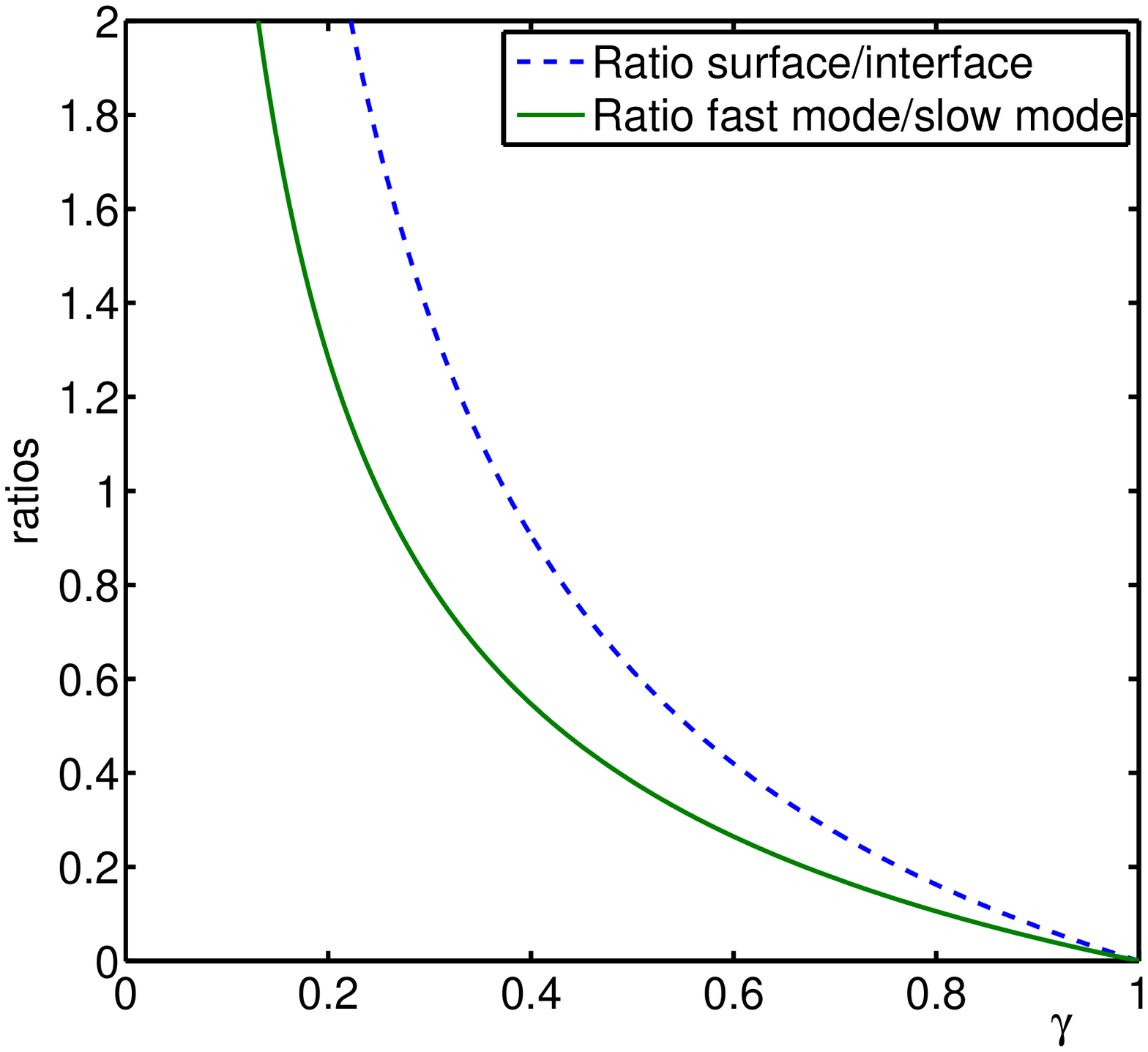}} 
\subfigure[$\delta=1$]{\includegraphics [width=0.3\textwidth,keepaspectratio=true]{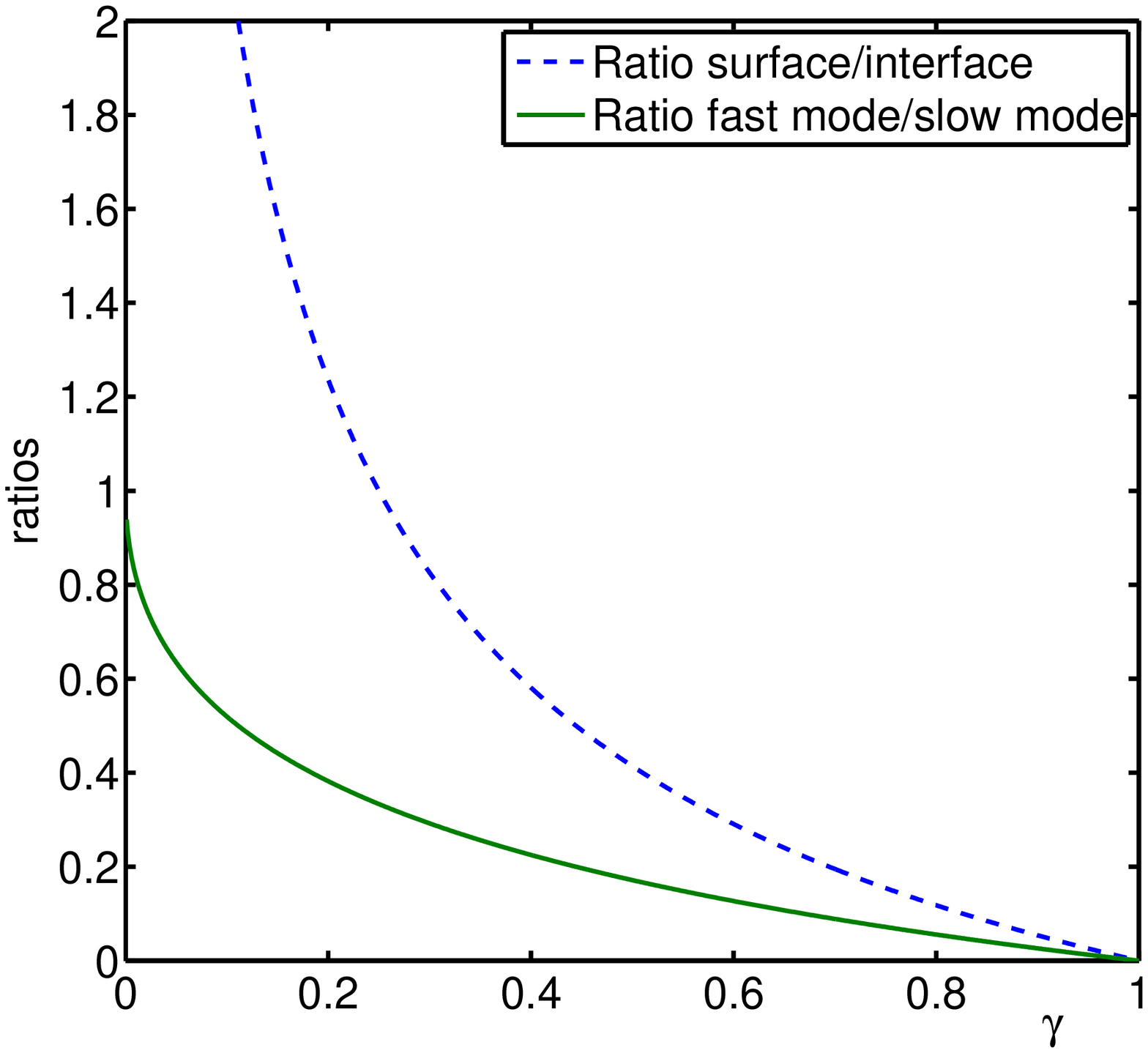}} 
\subfigure[$\delta=2$]{\includegraphics [width=0.3\textwidth,keepaspectratio=true]{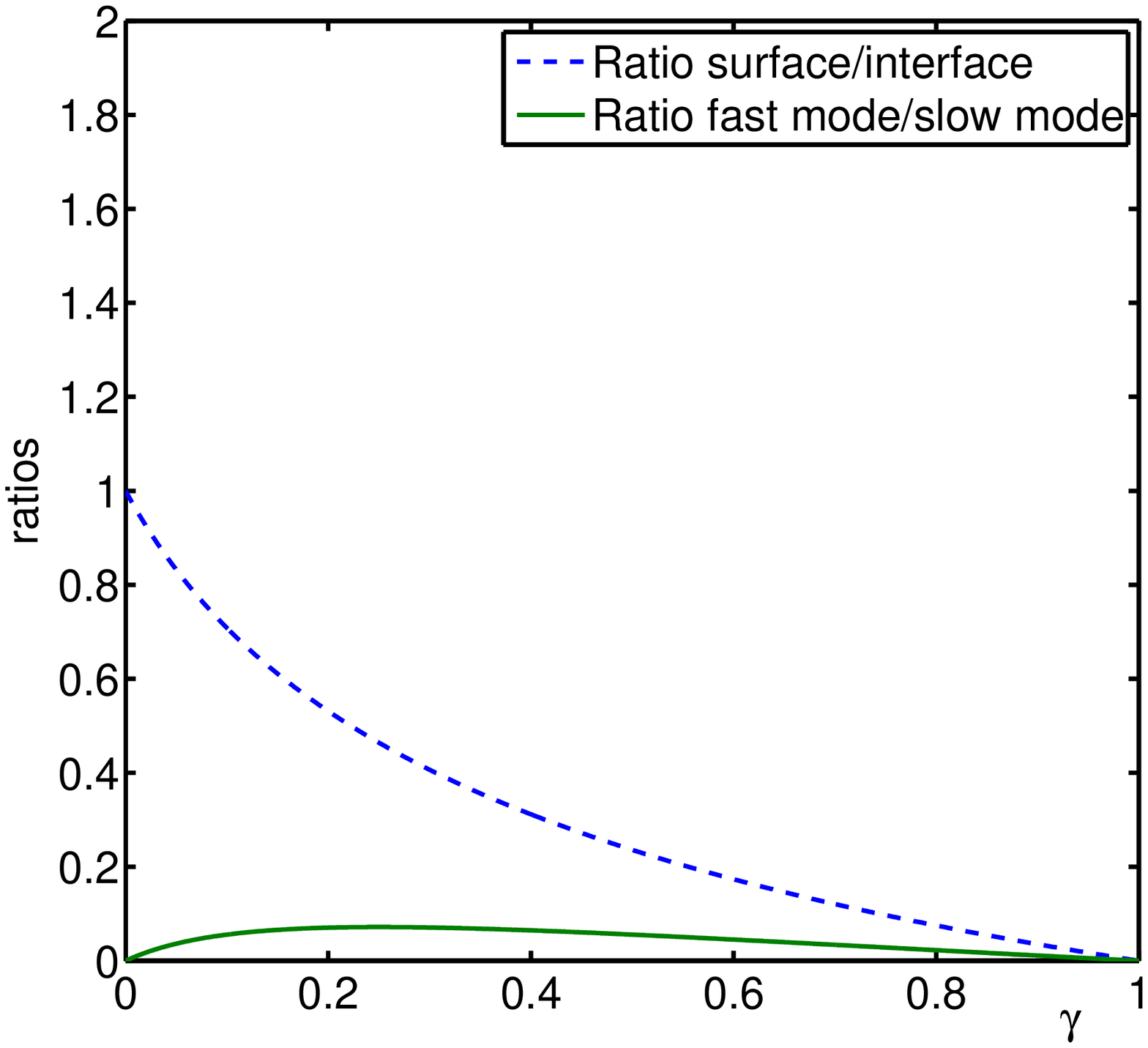}} 
 \caption{Magnitude of the deformations: the fast mode internal waves and the surface waves, when compared with the slow mode internal waves. Ratios for $\gamma\in(0,1)$ and (a) $\delta=1/2$,  (b) $\delta=1$,  (c) $\delta=2$.}
 \label{fig:ratiodeformations}
\end{figure}
We plot in Figure~\ref{fig:ratiodeformations} the different ratio of magnitudes, with the slow mode internal waves chosen as reference; that is to say
\[\frac{|\zeta_{\pm,+}|_{L^2}}{|\eta_{\pm,-}|_{L^2}} =\frac{|\zeta_{\pm,-}|_{L^2}}{|\eta_{\pm,-}|_{L^2}}   \ \ \text{ and }\ \ \ \frac{|\eta_{\pm,+}|_{L^2}}{|\eta_{\pm,-}|_{L^2}} \ \ .\]
The rigid lid hypothesis is valid for small values of these ratios. Again, one sees that it occurs only when $\gamma\sim1$, or when $\delta$ is big. This fact has already been addressed, for example in~\cite{Keulegan53,MichalletBarth'elemy98}, but its precise confirmation has never been exposed, to our knowledge. For example, if $\delta=1/2$ and $\gamma=0.8$, one has
\[\frac{|\zeta_{\pm,+}|_{L^2}}{|\eta_{\pm,-}|_{L^2}} =\frac{|\zeta_{\pm,-}|_{L^2}}{|\eta_{\pm,-}|_{L^2}}\sim \dfrac16   \ \ \text{ and }\ \ \ \frac{|\eta_{\pm,+}|_{L^2}}{|\eta_{\pm,-}|_{L^2}}\sim \dfrac1{10} \ \ ,\]
so that the rigid lid assumption is inaccurate.

\subsubsection{Summary}
\label{sSec:sum}
\begin{wrapfigure}{r}{0.5\textwidth}%\vspace*{-20pt}
             \includegraphics[width=0.5\textwidth,keepaspectratio=true]{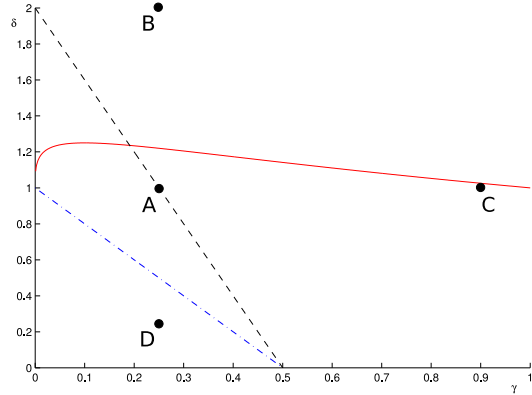} 
	     %\vspace*{-10pt}
 \caption{Properties of the KdV approximation, depending on the depths ratio $\delta$ and the density ratio $\gamma$.}
 %\vspace*{-10pt}
 \label{fig:recap}
\end{wrapfigure}
Let us summarize in Figure~\ref{fig:recap} the different results we obtained concerning the dependence of the behavior of the KdV approximation, depending on the parameters of the problem.

When $\delta$ is above the plain curve ($\delta>\delta_c$), then we know that the slow mode solitary waves will be of elevation type at the interface and of depression type at the surface, and conversely if $\delta<\delta_c$. The fast mode solitary waves are always of elevation type at the surface and at the interface.

Above the dashed line, the interface deformation is more important than the surface one for the slow mode solitary waves, and conversely if $\delta<2(1-2\gamma)$. As for the fast mode waves, the surface elevation is always bigger than the interface elevation.

The dash-dotted line concerns waves created by an initial data with zero velocities, and with a flat surface. In that case, the fast mode waves will be smaller than the slow mode waves above the line $\delta=1-2\gamma$, and conversely below.

\medskip

We see that for big values of $\gamma$ and/or $\delta$, the KdV approximation of the two-layer problem with a free surface gives a solution that resembles the interface problem with a rigid lid: the fast mode waves are smaller than the slow mode waves, and the magnitude of the deformation of the surface is of less importance than the deformation of the interface. On the contrary, when $\gamma$ and $\delta$ are small, then the solutions of our problem, when considered at the interface, are comparable to solutions of the one-layer water wave models.

The dots on Figure~\ref{fig:recap} represent the pair of parameters $(\gamma,\delta)$ for which numerical simulations have been computed in Section~\ref{sSec:NumericalResults}. 
 We have first computed the solutions of the symmetric Boussinesq/Boussinesq model and the KdV approximation, for both for the case of solitary waves, and zero velocities-flat surface initial data. The results are plotted for parameters corresponding to \textbf{A} ($\gamma=1/4,\delta=1$) in Figures~\ref{fig:solitonA} and~\ref{fig:clocheA}, and to \textbf{B} ($\gamma=1/4,\delta=2$) in Figures~\ref{fig:solitonB} and~\ref{fig:clocheB}. Then, we have compared the KdV approximation in the rigid lid and free surface configurations. Figure~\ref{fig:rigidlidA},~\ref{fig:rigidlidB},~\ref{fig:rigidlidC} and~\ref{fig:rigidlidD} corresponds respectively to the points \textbf{A}, \textbf{B}, \textbf{C} ($\gamma=9/10,\delta=1$) and \textbf{D} ($\gamma=\delta=1/4$).

\section{Numerical comparison}
\label{Sec:NumericalComparison}
\subsection{The numerical schemes}
\label{sSec:NumericalSchemes}
This section is devoted to the numerical comparison between the different models displayed in this article, namely the symmetric Boussinesq/Boussinesq model~\eqref{SBOUSS}, and the KdV equations~\eqref{eq:KdV}.
We first provide a numerical scheme for generic KdV equations that can easily be adjusted for the uncoupled KdV approximations~\eqref{eq:KdV} and~\eqref{eq:KdV2}, and its adaptation to the Boussinesq/Boussinesq system~\eqref{SBOUSS}. 

\medskip

Each time, we use a Crank Nicholson scheme and replace the costly numerical treatment of the nonlinear term by a predictive step. This method has been introduced by Besse and Bruneau in~\cite{BesseBruneau98}, justified in~\cite{Besse98}, and used in the water wave framework by Chazel in~\cite{Chazel09}, and more recently by Durufl\'e and Israwi in~\cite{DurufleIsrawi}. The method is formally of order two in space and time, which is confirmed by the simulations, and appears to be unconditionally stable.

\subsubsection{The KdV equation}
We present here the numerical scheme for the generic KdV equation:
\begin{equation}
 \label{KdVg} \partial_t u +c \partial_x u +\lambda u\partial_x u+\mu\partial_x^3 u.
\end{equation}

First, we use the following semi-discretized in time equation:
\[\begin{array}{r}\dfrac{u^{n+1}-u^n}{dt}+c\partial_x \left(\dfrac{u^{n+1}+u^n}{2}\right)+\lambda \left(\alpha u^{n+1/2} \partial_x \left(\frac{u^{n+1}+u^n}{2}\right) +(1-\alpha) \partial_x u^{n+1/2} \left(\frac{u^{n+1}+u^n}{2}\right)\right) \\ +\mu\partial_x^3\left(\frac{u^{n+1}+u^n}{2}\right)=0,
 \end{array}
\]
with $\alpha\in[0,1]$, and where $u^{n+1/2}$ is a predictive term defined by
\begin{equation}\label{un+1/2}u^n=\frac{u^{n+1/2}+u^{n-1/2}}{2}.\end{equation}

The scheme takes advantage of the two possible discretizations of the nonlinear term $u\partial_x u$, that is to say $u^{n+1/2} \partial_x \left(\frac{u^{n+1}+u^n}{2}\right)$ and  $\frac{u^{n+1}+u^n}{2} \partial_x\left(u^{n+1/2}\right)$, by introducing a parameter $\alpha\in[0,1]$ and taking a convex combination of these possibilities.

It is easy to check that, in order to preserve the semi-discrete $L^2$-norm, one has to choose $\alpha=2/3$. As for the spatial discretization, we use the Crank-Nicholson scheme, adjusted so that the discrete $L^2$-norm is preserved. This leads to the final discretization:
\begin{align}\label{KdVscheme}
\frac{u^{n+1}_i-u^n_i}{dt}+ c\left(D_1\frac{u^{n+1}+u^n}{2}\right)_i+\frac{\lambda}{3} \left( \left(\frac{u^{n+1}_{i+1}+u^n_{i+1}+u^{n+1}_{i-1}+u^n_{i-1}}{4}\right) \left(D_1 u^{n+1/2}\right)_i\right. &\nonumber\\ \left. +\left(u^{n+1/2}_i+\frac{u^{n+1/2}_{i+1}+u^{n+1/2}_{i-1}}{2} \right) \left(D_1 \frac{u^{n+1}+u^n}{2}\right)_i\right) +\mu \left(D_3 \frac{u^{n+1}+u^n}{2}\right)_i& \ = \ 0,
\end{align}
with $D_1$ and $D_3$ the classical centered discretizations of the derivatives $\partial_x$ and $\partial_x^3$ with periodic boundary conditions. The scheme is given at each step by~\eqref{KdVscheme} and~\eqref{un+1/2}, with a simple explicit scheme for the first half step. One can then check (see~\cite[Theorem 2]{DurufleIsrawi}):
 \[\forall\ n\in\NN,\quad \sum_i |u_i^n|^2 \ = \  \sum_i |u_i^0|^2. \]

\subsubsection{The Boussinesq/Boussinesq system} 
We use the same ideas as in the previous section for the discretization of the Boussinesq/Boussinesq system.
Even if the $L^2$ norm is not preserved by~\eqref{SBOUSS}, and no simple quantity either, we decide to use the same parameter $\alpha=2/3$ as in the previous section. This leads to the following discretization of the spatial nonlinear term:
\[ \Sigma_1(U)\partial_x U \sim \frac{2}{3}\Sigma_1(U^{n+1/2})\partial_x \left( \frac{U^{n+1}+U^n}{2}\right) + \frac13 \Sigma_1\left(\frac{U^{n+1}+U^n}{2}\right)\partial_x U^{n+1/2}.\]
As for the nonlinear term in time, we simply use
\[ S_1(U)\partial_t U \sim S_1(U^{n+1/2})\partial_t \left( \frac{U^{n+1}+U^n}{2}\right).\] 
We also need to construct the linear mappings $\t \Sigma_1$ with values on $\M_4(\RR)$, such that 
\[ \forall U,V\in \RR^4, \quad \Sigma_1(U)V\ =\ \t \Sigma_1(V)U.\]
This finally leads to the following scheme:
\begin{equation}\label{Bousscheme}
\begin{array}{r}\displaystyle\left(S_0+\epsilon S_1(U^{n+1/2}_i) -\epsilon S_2 D_2\right)\frac{U^{n+1}_i-U^n_i}{dt} \ + \left(\Sigma_0 D_1 \frac{U^{n+1}+U^n}{2}\right)_i-
\epsilon\left(\Sigma_2 D_3 \frac{U^{n+1}+U^n}{2}\right)_i \\ \displaystyle+\frac{\epsilon}{3} \left( \Sigma_1\left(U^{n+1/2}_i+\frac{U^{n+1/2}_{i+1}+U^{n+1/2}_{i-1}}{2} \right) \left(D_1 \frac{U^{n+1}+U^n}{2}\right)_i \right.\\
\displaystyle\left. + \left(\t \Sigma_1\left(D_1 U^{n+1/2}\right)\frac{U^{n+1}+U^n}{2}\right)_i \right) =0,
 \end{array}
\end{equation}
where $D_1$, $D_2$ and $D_3$ are block-diagonal, with the classical centered discretizations of the derivatives $\partial_x$, $\partial_x^2$ and $\partial_x^3$ (with periodic boundary conditions) as diagonal blocks.

The scheme is given at each step by~\eqref{Bousscheme} and~\eqref{un+1/2}, with a simple explicit scheme for the first half step.

\subsubsection{Validation of the numerical method}\label{SubsubsectionValidationOfTheNumericalMethod}
As said previously, the method is formally of order 2 in space and time, which is confirmed by the simulations, and appears to be unconditionally stable.
In order to validate the schemes, we use the known solitary wave solutions of~\eqref{KdVg}, expressed as follows:
\begin{equation}\label{solitonscheme}
 u(t,x)=\frac{M}{\cosh^{2}(k(x-c' t))},
\end{equation}
with $c'=c+\frac{\lambda M}{3}$, $k=\sqrt\frac{\lambda M}{12 \mu}$, and $M$ arbitrary.

Therefore, we are able to construct an initial data that will lead the solutions of the KdV approximation defined in Theorem~\ref{Th:CVKdVEuler} to be steady traveling waves. Indeed, if we set $\gamma\in(0,1)$ and $\delta>0$, the coefficients of the KdV equations solved by the KdV approximation are explicit and given in Remark~\ref{rem:coefsKdV} page \pageref{rem:coefsKdV}. 
Consequently, we choose $M_1,\dots,M_4\in\RR$, and set $U^0=\sum_{i=1}^4 u_i(x,0) \e_i$, with $u_i$ of the form given by~\eqref{solitonscheme}, with the given parameters, and the uncoupled KdV approximation will be given by $U(t,x)=\sum_{i=1}^4 u_i(\epsilon t,x-c_i t)\e_i$.

Unfortunately, we cannot exhibit such an exact solution for the symmetric Boussinesq/Boussinesq system~\eqref{SBOUSS}. In order to validate the scheme~\eqref{Bousscheme}, we plug the solution of the uncoupled KdV approximation in~\eqref{SBOUSS}, and obtain a forcing term $F(t,x)$, so that $U(t,x)$ is solution of the modified system
\[\Big(S_0+\epsilon \big(S_1(U)-S_2 \partial_x^2\big)\Big)\partial_t U+\Big(\Sigma_0+\epsilon\big(\Sigma_1(U)-\Sigma_2 \partial_x^2\big)\Big)\partial_x U=F.\]
It is trivial then to modify the scheme~\eqref{Bousscheme} by adding a forcing term $F^n_i=F(t^n,x_i)$.

We present in Table~\ref{tab:ErSchemas} the results that we obtain using the scheme~\eqref{KdVscheme} and~\eqref{Bousscheme}, for several values of spatial and time discretization steps $dx$, $dt$, and different values of $\epsilon$, and for times $T=1/\epsilon$. The relative errors are computed in the discrete $L^2$ norm. These results allow to validate the schemes proposed. 
\begin{table}[!ht]
\centering {\begin{tabular}{@{}ccccccc@{}} \Hline
 $dx$ & $dt$ & $L$ & $T$ & $\epsilon$ & KdV scheme & Boussinesq scheme \\ \hline
 0.01 & 0.01 & 120 & 5\hphantom{0} & 0.2\hphantom{0} & $9.6317.10^{-5}$ & $9.8719.10^{-5}$ \\  
 0.02 & 0.02 & 120 & 10 & 0.1\hphantom{0} & $7.7094.10^{-4}$ & $7.9861.10^{-4}$ \\  
 0.05 & 0.05 & 120 & 20 & 0.05 & $9.5663.10^{-3}$ & $9.8587.10^{-3}$\\ \Hline
 \end{tabular}}
 \caption{Numerical errors of the KdV and the Boussinesq schemes.}
\label{tab:ErSchemas}
\end{table}

In the following simulations, we always choose $dx=dt=0.01$.

\subsection{Numerical results}
\label{sSec:NumericalResults}
 We compute our schemes for two different different forms of initial data. The first leads to solitary waves as exact solutions of the KdV approximation (thus of the form~\eqref{soliton}), and is therefore exponentially decreasing in space. The other initial data consist in deformations of the interface of the form $\frac{M}{\sqrt{1+k x^2}}$ (therefore, they do not satisfy the spatial rapidly decreasing assumption for the better convergence rate in Proposition~\ref{convergenceKdVBouss}), as the surface is flat, and initial velocities are zeros. All the forthcoming results are expressed in non-dimensionalized variables.
 
 \medskip
 
 We first look at the behavior of the relative difference between the solutions of the different models for fixed time, and for different values of $\epsilon$. In Table~\ref{tab:ErrorInEpsilonBousSymBous}, we provide the difference between the original Boussinesq/Boussinesq model~\eqref{Bouss} and our symmetric Boussinesq/Boussinesq system~\eqref{SBOUSS}, at time $T=1$ and for different values of $\epsilon$. The same results are given for the comparison between the symmetric Boussinesq/Boussinesq system and the KdV approximation in Table~\ref{tab:ErrorInEpsilonBousKdV}. As predicted, when compared with the solution of the symmetric Boussinesq/Boussinesq model, the deviation of the original Boussinesq/Boussinesq system at fixed time is of order $\O(\epsilon^2)$, whereas the error of the KdV approximation is of order $\O(\epsilon)$. This supports our choice of the symmetric Boussinesq/Boussinesq system~\eqref{SBOUSS}, as an equivalent model for system~\eqref{Bouss}.
\begin{table}[!ht]
\centering    {\begin{tabular}{@{}ccccc|cc@{}} \Hline
 $dx$ & $dt$ & $L$ & $T$ & $\epsilon$ & \multicolumn{2}{c}{relative error of the solutions with the initial value:} \\
  & & & & & $U^0=\dfrac{M}{\cosh(kx)^2}$ & $U^0=\dfrac{M}{\sqrt{1+(kx)^2}}$ \\ \hline
 0.01 & 0.01 & 120 & 1 & 0.1\hphantom{0} & $1.4859.10^{-3}$ & $2.4785.10^{-3}$ \\ 
 0.01 & 0.01 & 120 & 1 & 0.05 & $5.2714.10^{-4}$ & $1.0882.10^{-3}$\\ 
 0.01 & 0.01 & 120 & 1 & 0.01 & $3.1656.10^{-5}$ & $7.1958.10^{-5}$ \\ \Hline
 \end{tabular}}
 \caption{Relative error between the solutions of the symmetric Boussinesq/Boussinesq system and the original Boussinesq/Boussinesq model.}
 \label{tab:ErrorInEpsilonBousSymBous}
\end{table}\begin{table}[!ht]
 \centering  
{\begin{tabular}{@{}ccccc|cc@{}} \Hline
 $dx$ & $dt$ & $L$ & $T$ & $\epsilon$ & \multicolumn{2}{c}{relative error of the solutions with the initial value: } \\
  & & & & & $U^0=\dfrac{M}{\cosh(kx)^2}$ & $U^0=\dfrac{M}{\sqrt{1+(kx)^2}}$ \\ \hline
 0.01 & 0.01 & 120 & 1 & 0.1\hphantom{0} & $4.0805.10^{-3}$ & $5.9295.10^{-3}$ \\ 
 0.01 & 0.01 & 120 & 1 & 0.05 & $2.2785.10^{-3}$ & $2.5506.10^{-3}$\\ 
 0.01 & 0.01 & 120 & 1 & 0.01 & $5.2496.10^{-4}$ & $3.5434.10^{-4}$ \\ \Hline
 \end{tabular}}
 \caption{Relative error between the solution of the symmetric Boussinesq/Boussinesq system and the solution of the KdV approximation.}
  \label{tab:ErrorInEpsilonBousKdV}
\end{table}

\bigskip

Then, we give a numerical confirmation of the convergence rate obtained in Proposition~\ref{convergenceKdVBouss}. It is stated that for any initial data, the difference between the solution of the Boussinesq/Boussinesq system and the KdV approximation is bounded by $\O(\epsilon)$ over times of order $\O(1/\epsilon)$ if the initial data is sufficiently decreasing in space, and $\O(\epsilon\sqrt t )$ otherwise. 

\begin{figure}[!ht]
 \centering
\subfigure[solitary waves]{\includegraphics [width=0.45\textwidth,keepaspectratio=true]{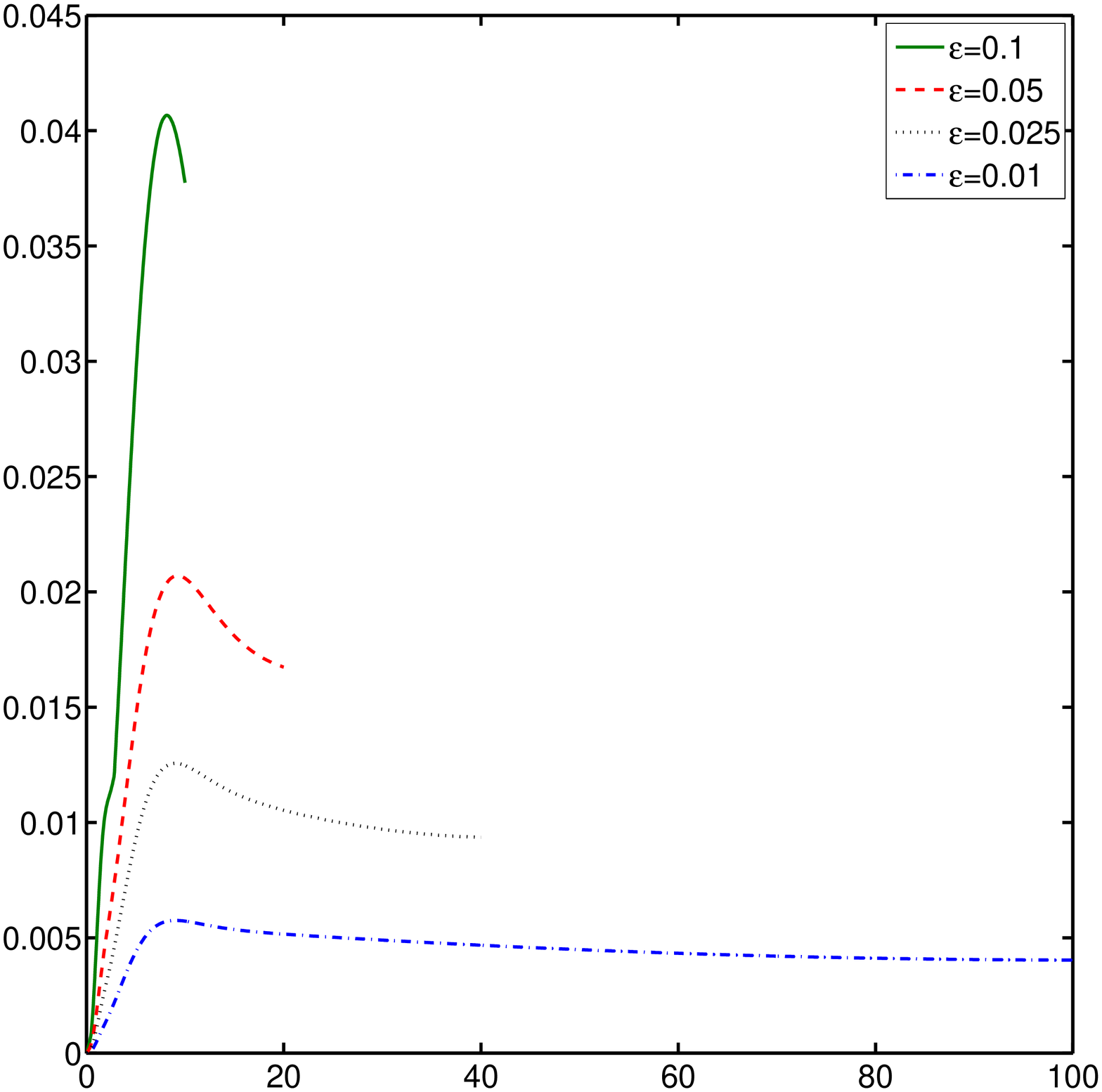}} 
\subfigure[zero velocities-flat surface initial value]{\includegraphics [width=0.45\textwidth,keepaspectratio=true]{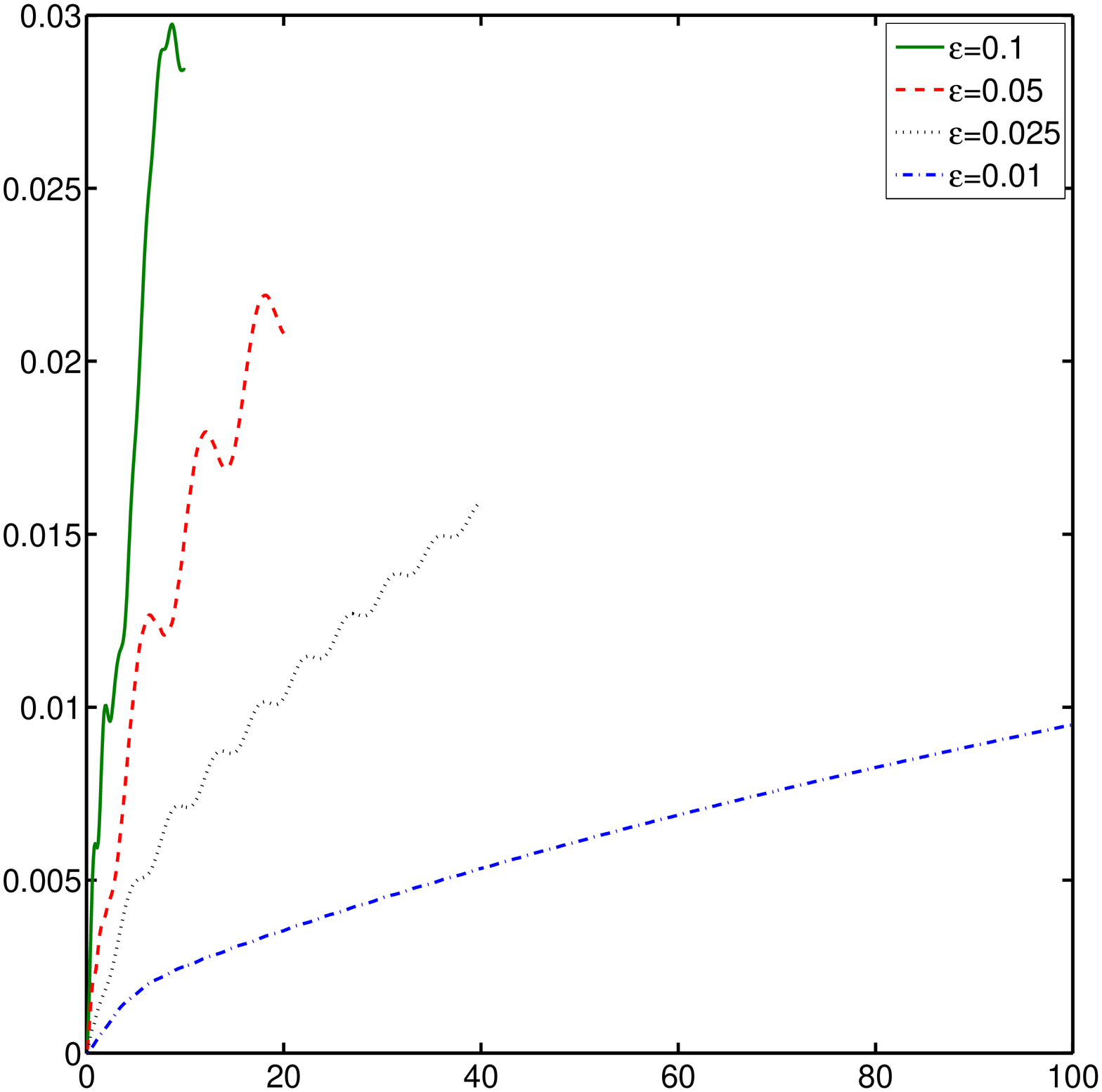}} 
 \caption{Relative error between the KdV approximation and the symmetric Boussinesq/Boussinesq system: (a) for solitary waves, exponentially decreasing in space (b) for zero velocities-flat surface initial data, slowly decreasing in space.}
 \label{fig:ErrorInTime}
\end{figure}

We plot in Figure~\ref{fig:ErrorInTime} the difference (in the discrete $L^2$ norm) of two solutions, obtained  respectively by the KdV scheme and the Boussinesq scheme. The first plot concerns solitary waves of the form~\eqref{soliton}, thus exponentially decreasing in space, and the second is given by a zero velocities-flat surface initial data, with the interface of the form $\frac{M}{\sqrt{1+(kx)^2}}$. 

For both configurations, we set $\delta=1$, $\gamma=1/4$, and simulate at the different values $\epsilon=0.1,0.05,0.025,0.01$, throughout time $T=1/\epsilon$. One sees that the results match the theory, and that the decreasing in space at infinity is indeed of great concern. 

 \begin{figure}[hp]
 \centering
\subfigure[$t=20$. $\epsilon=0.1,\delta=1,\gamma=1/4$.]{\includegraphics [width=0.45\textwidth,keepaspectratio=true]{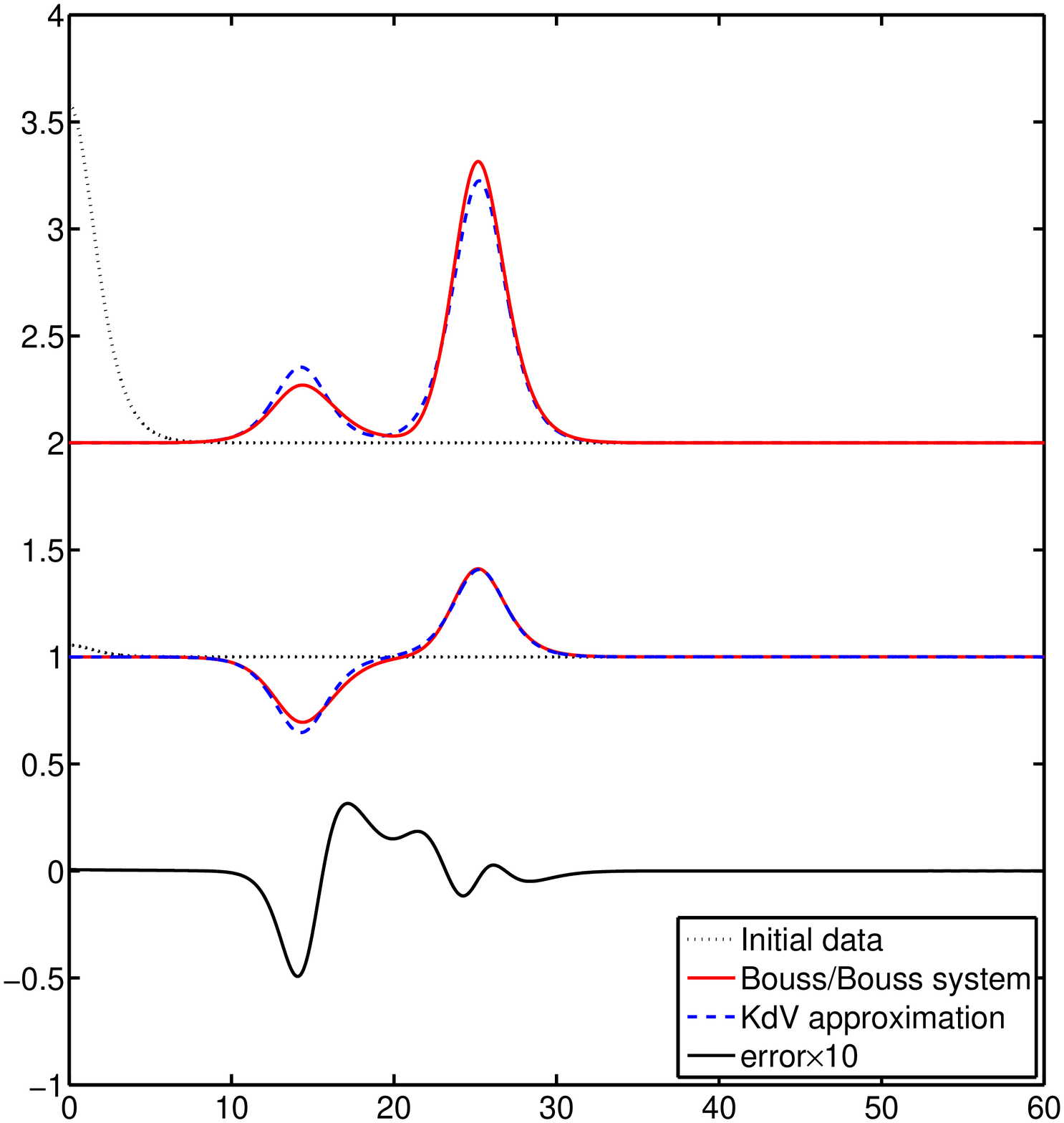}} 
\subfigure[$t=40$. $\epsilon=0.1,\delta=1,\gamma=1/4$.]{\includegraphics [width=0.45\textwidth,keepaspectratio=true]{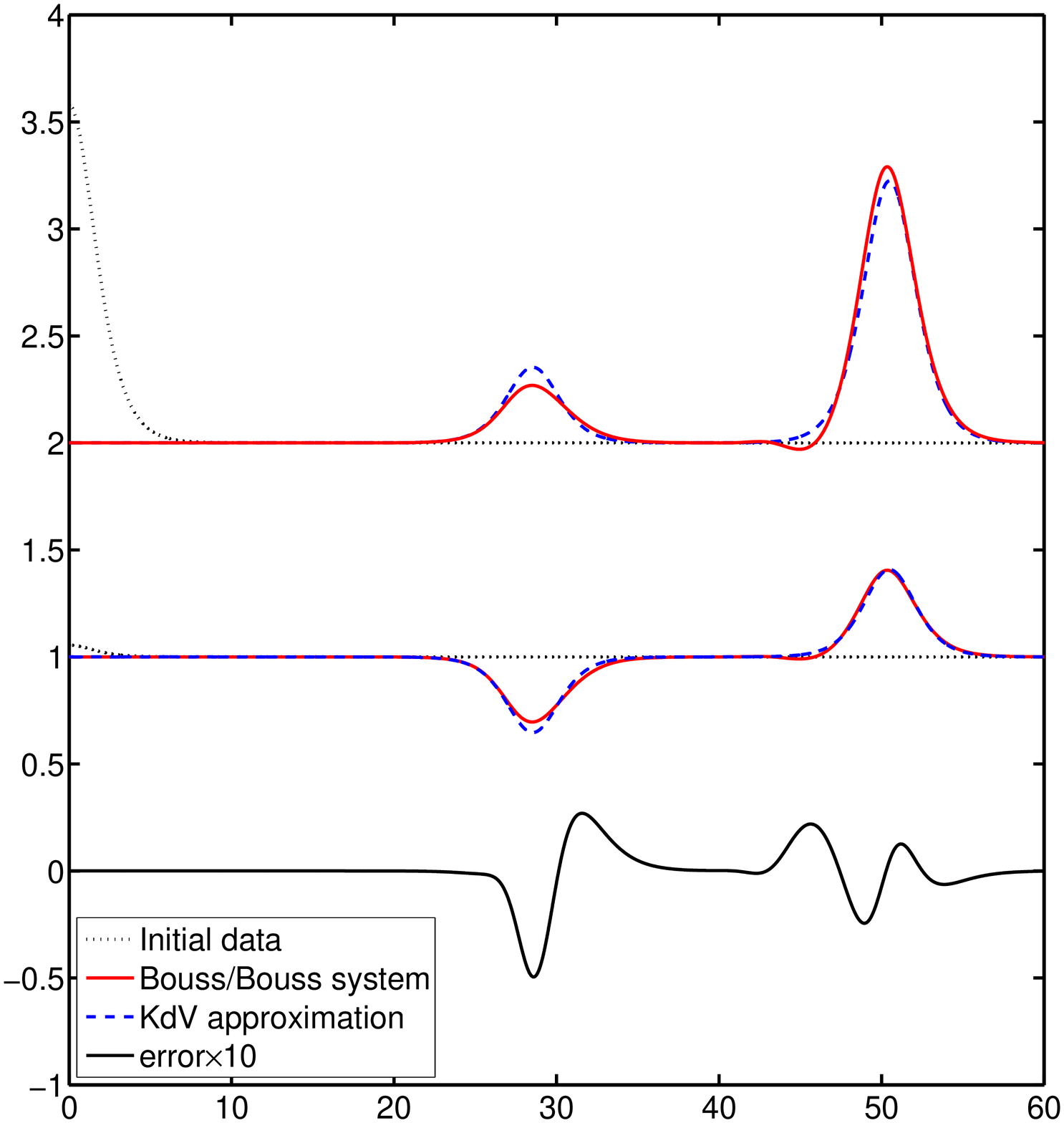}} 
 \caption{Solitary wave solution of the KdV approximation, and the symmetric Boussinesq/Boussinesq model, at times (a) t=20, (b) t=40.}
 \label{fig:solitonA}
\end{figure}

\begin{figure}[hp]
 \centering
\subfigure[$t=20$. $\epsilon=0.1,\delta=1,\gamma=1/4$.]{\includegraphics [width=0.45\textwidth,keepaspectratio=true]{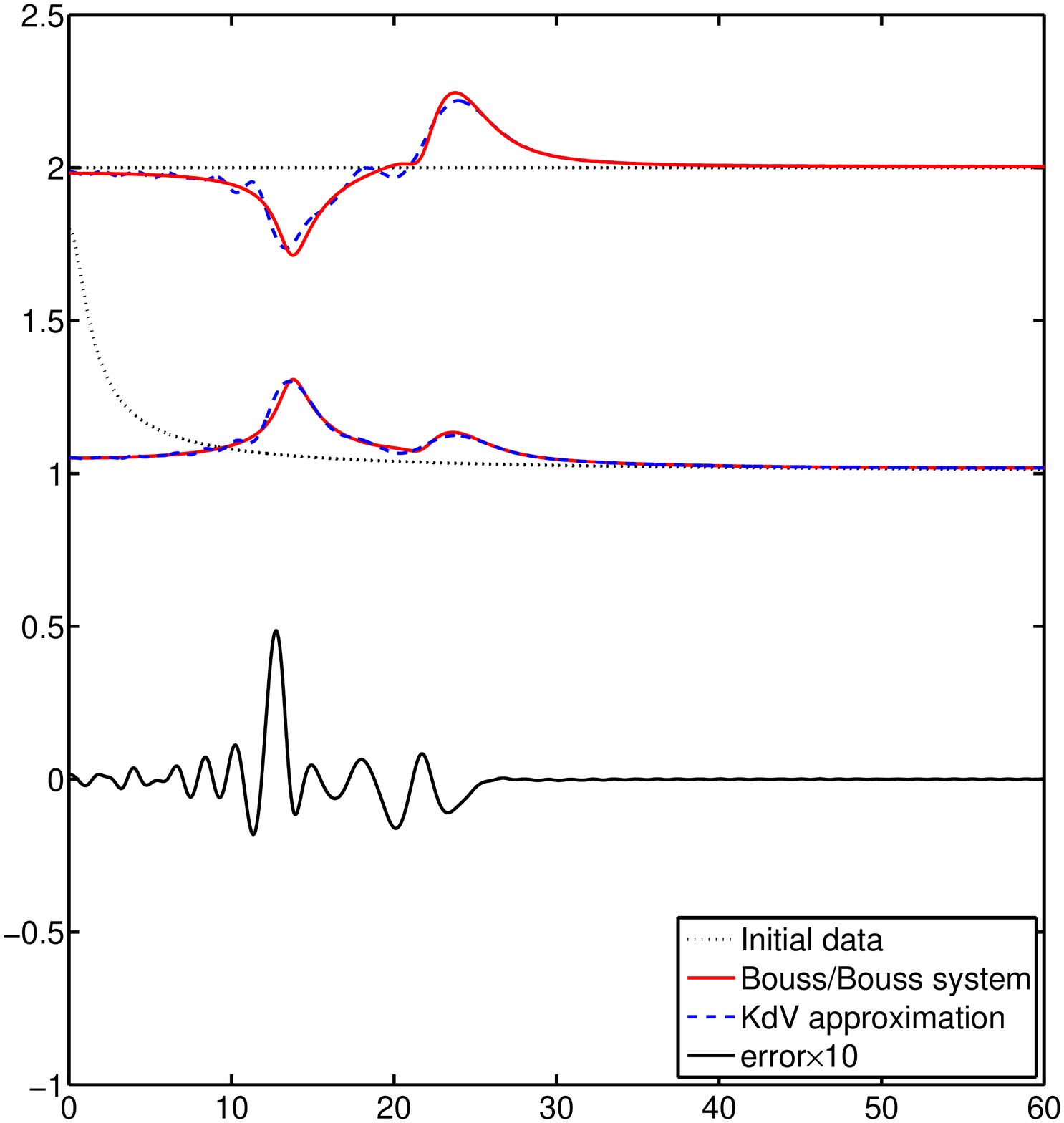}} 
\subfigure[$t=40$. $\epsilon=0.1,\delta=1,\gamma=1/4$.]{\includegraphics [width=0.45\textwidth,keepaspectratio=true]{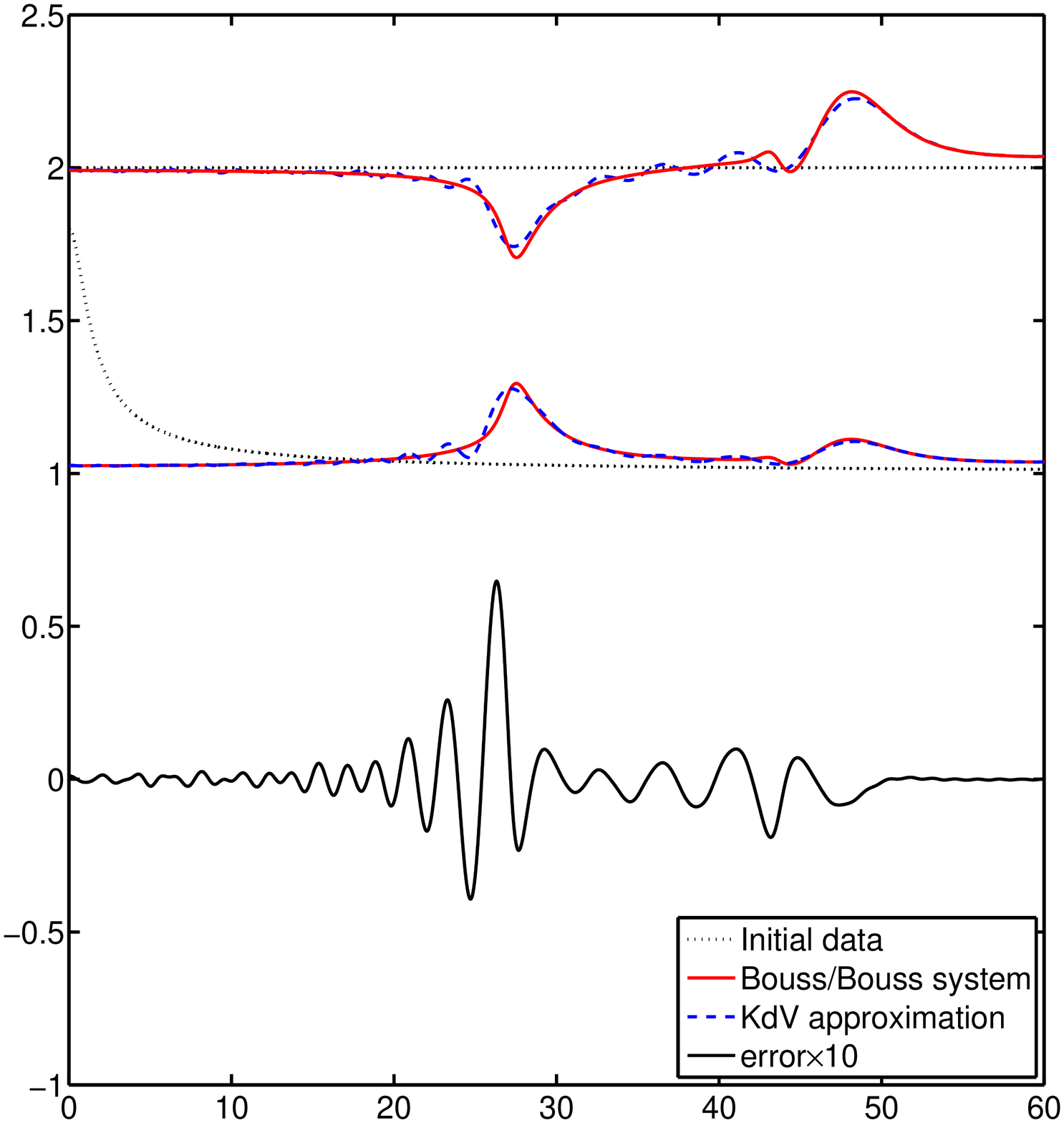}} 
 \caption{Solution of the KdV approximation, and symmetric Boussinesq/Boussinesq system for a zero velocities-flat surface initial value, at times (a) t=20, (b) t=40.}
 \label{fig:clocheA}
\end{figure}
 \begin{figure}[hp]
 \centering
\subfigure[$t=20$. $\epsilon=0.1,\delta=2,\gamma=1/4$.]{\includegraphics [width=0.45\textwidth,keepaspectratio=true]{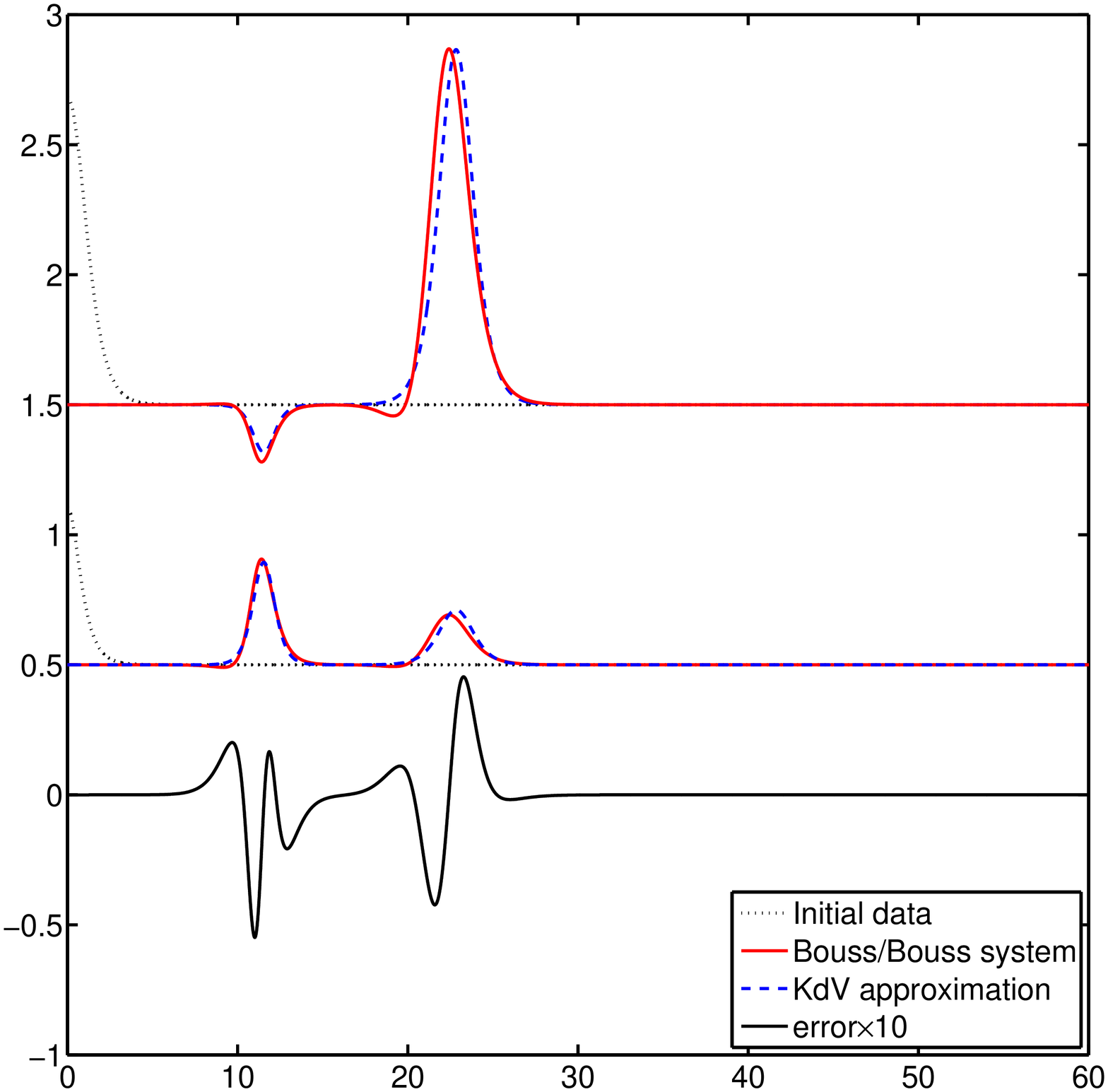}} 
\subfigure[$t=40$. $\epsilon=0.1,\delta=2,\gamma=1/4$.]{\includegraphics [width=0.45\textwidth,keepaspectratio=true]{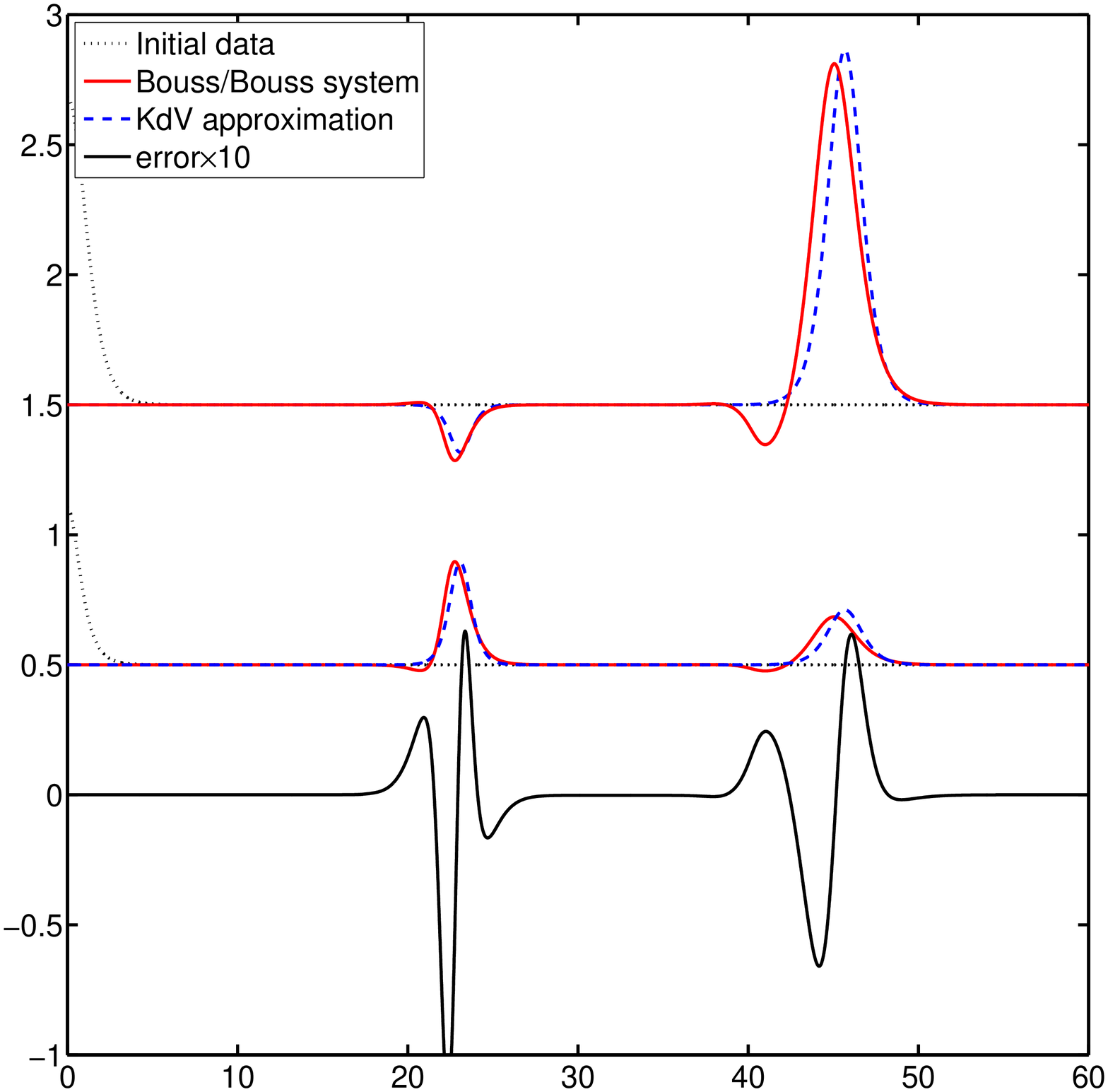}} 
 \caption{Solitary wave solution of the KdV approximation, and symmetric Boussinesq/Boussinesq model, at times (a) t=20, (b) t=40.}
 \label{fig:solitonB}
\end{figure}

\begin{figure}[hp]
 \centering
\subfigure[$t=20$. $\epsilon=0.1,\delta=2,\gamma=1/4$.]{\includegraphics [width=0.45\textwidth,keepaspectratio=true]{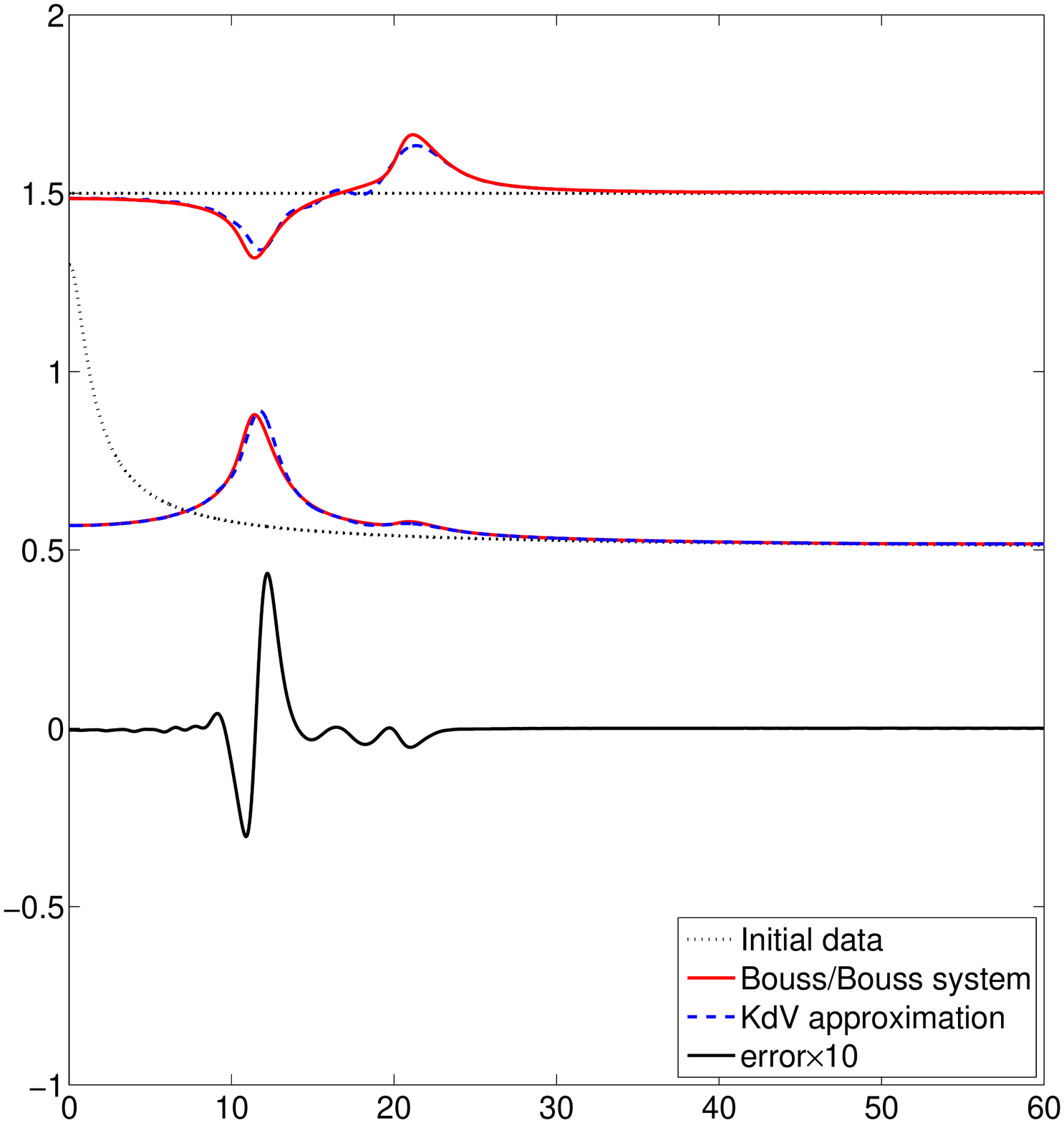}} 
\subfigure[$t=40$. $\epsilon=0.1,\delta=2,\gamma=1/4$.]{\includegraphics [width=0.45\textwidth,keepaspectratio=true]{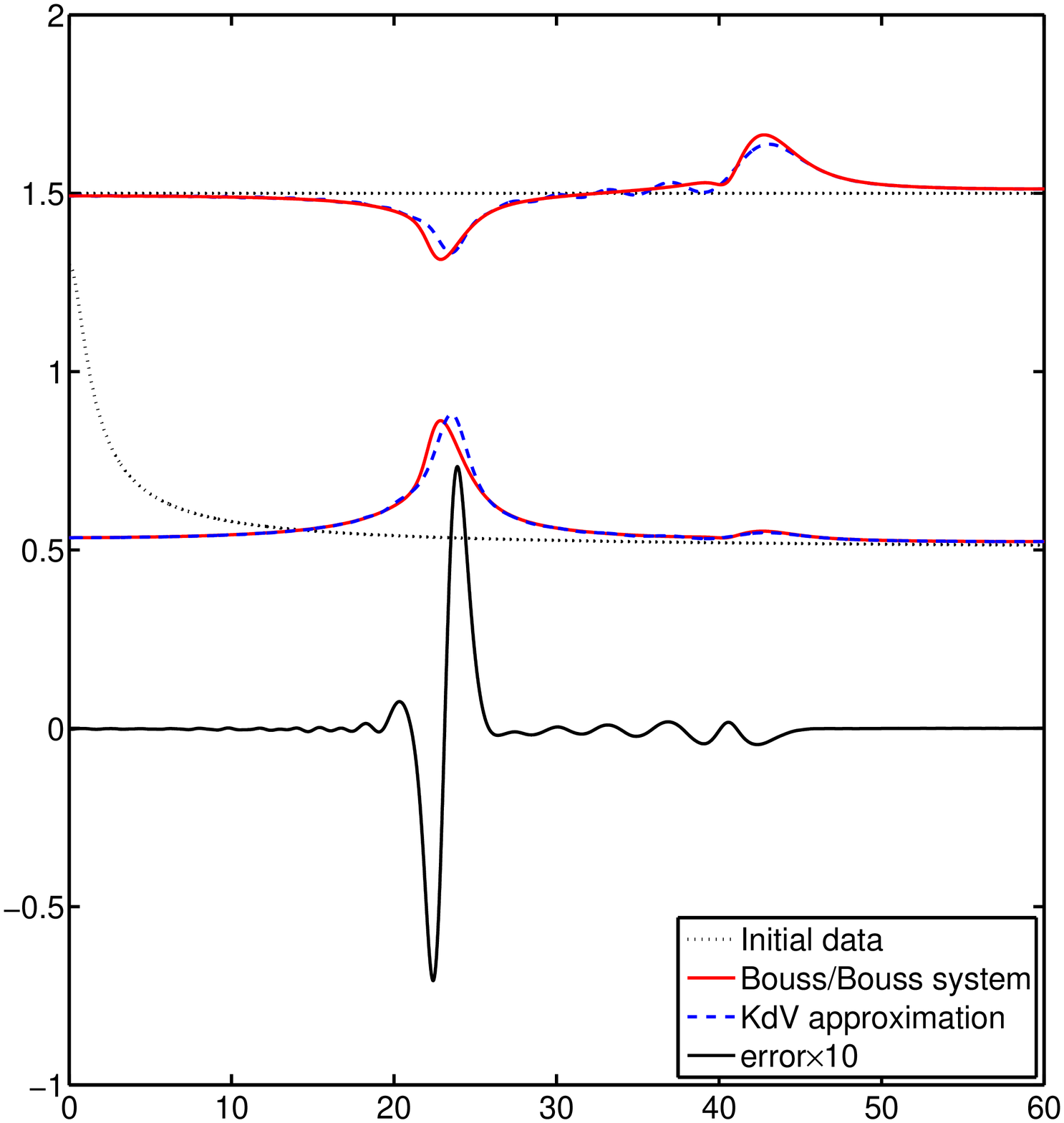}} 
 \caption{Solution of the KdV approximation, and symmetric Boussinesq/Boussinesq system for a zero velocities-flat surface initial value, at times (a) t=20, (b) t=40.}
 \label{fig:clocheB}
\end{figure}
 \begin{figure}[hp]
 \centering
\subfigure[$t=20$. $\epsilon=0.1,\delta=1,\gamma=1/4$.]{\includegraphics [width=0.45\textwidth,keepaspectratio=true]{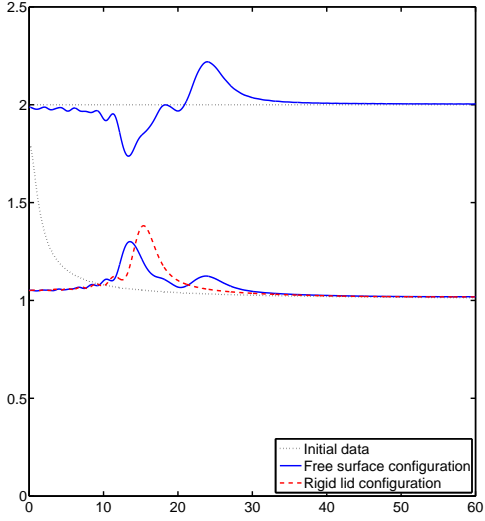}} 
\subfigure[$t=40$. $\epsilon=0.1,\delta=1,\gamma=1/4$.]{\includegraphics [width=0.45\textwidth,keepaspectratio=true]{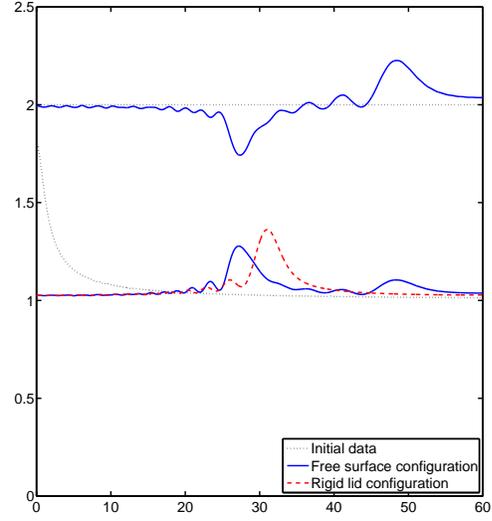}} 
 \caption{Solution of the KdV approximations, in both rigid lid and free surface configurations, with two fluids of highly different densities, at times (a) t=20, (b) t=40.}
 \label{fig:rigidlidA}
\end{figure}
 \begin{figure}[hp]
 \centering
\subfigure[$t=20$. $\epsilon=0.1,\delta=1,\gamma=9/10$.]{\includegraphics [width=0.45\textwidth,keepaspectratio=true]{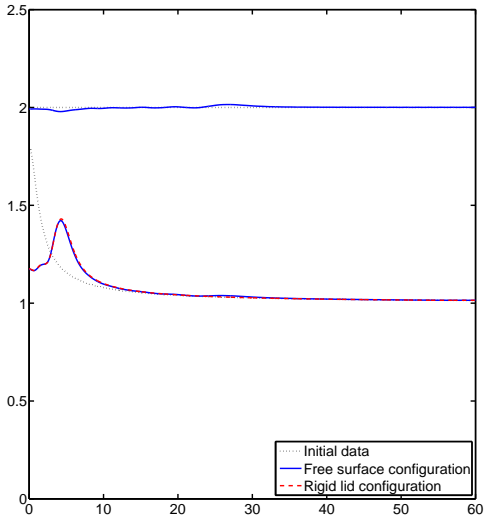}} 
\subfigure[$t=40$. $\epsilon=0.1,\delta=1,\gamma=9/10$.]{\includegraphics [width=0.45\textwidth,keepaspectratio=true]{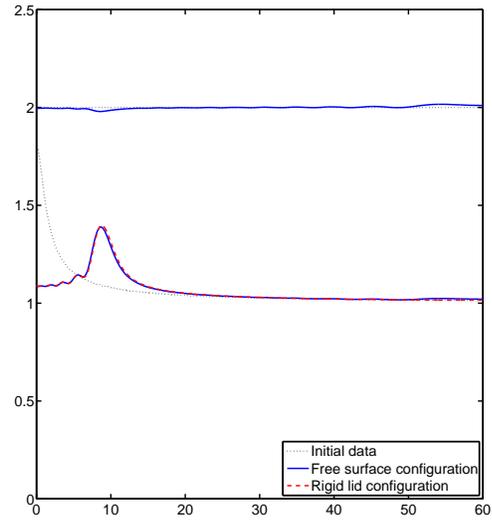}} 
 \caption{Solution of the KdV approximations, in both rigid lid and free surface configurations, with two fluids of near equal densities, at times (a) t=20, (b) t=40.}
 \label{fig:rigidlidC}
\end{figure}
 \begin{figure}[hp]
 \centering
\subfigure[$t=20$. $\epsilon=0.1,\delta=2,\gamma=1/4$.]{\includegraphics [width=0.45\textwidth,keepaspectratio=true]{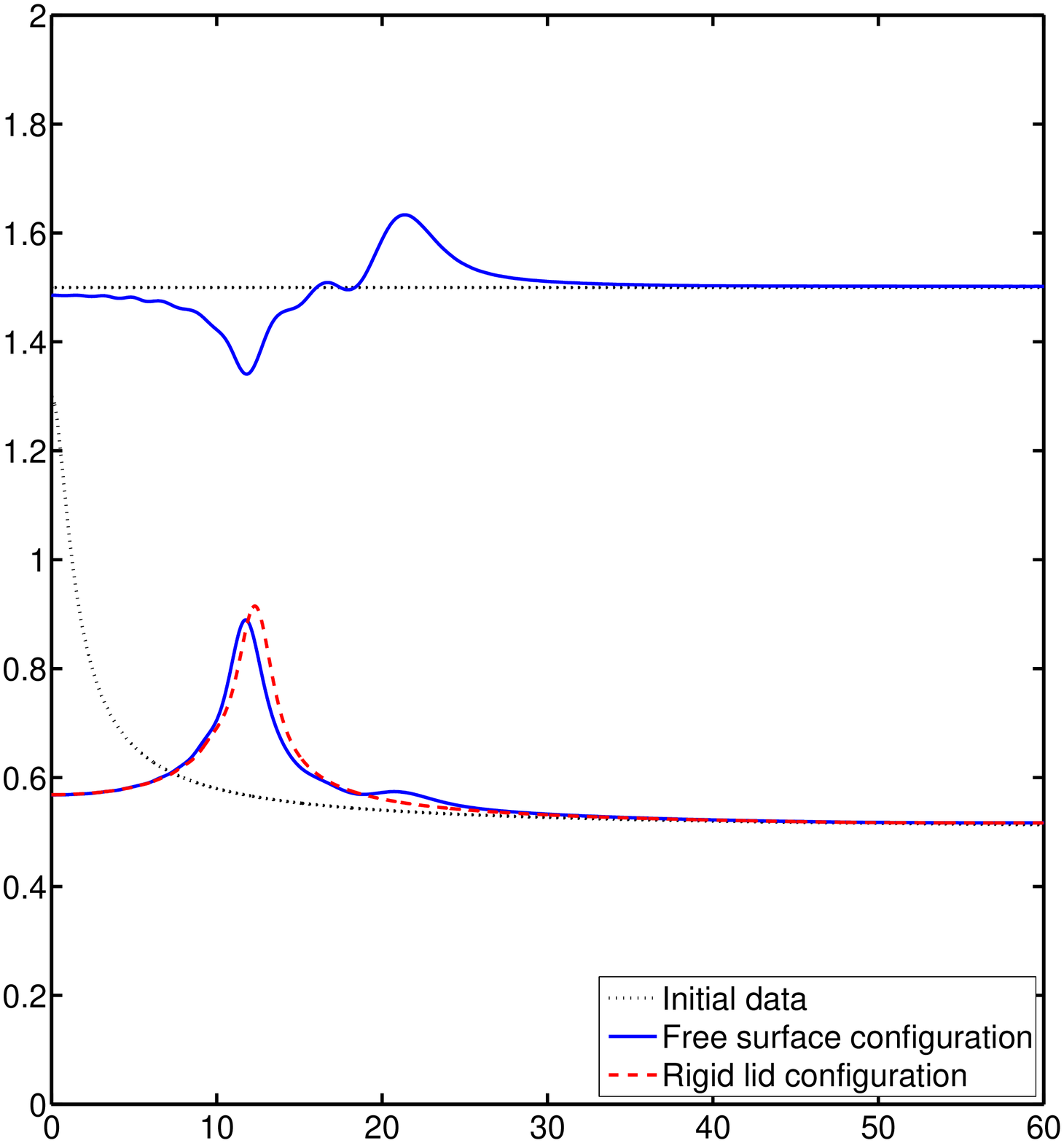}} 
\subfigure[$t=40$. $\epsilon=0.1,\delta=2,\gamma=1/4$.]{\includegraphics [width=0.45\textwidth,keepaspectratio=true]{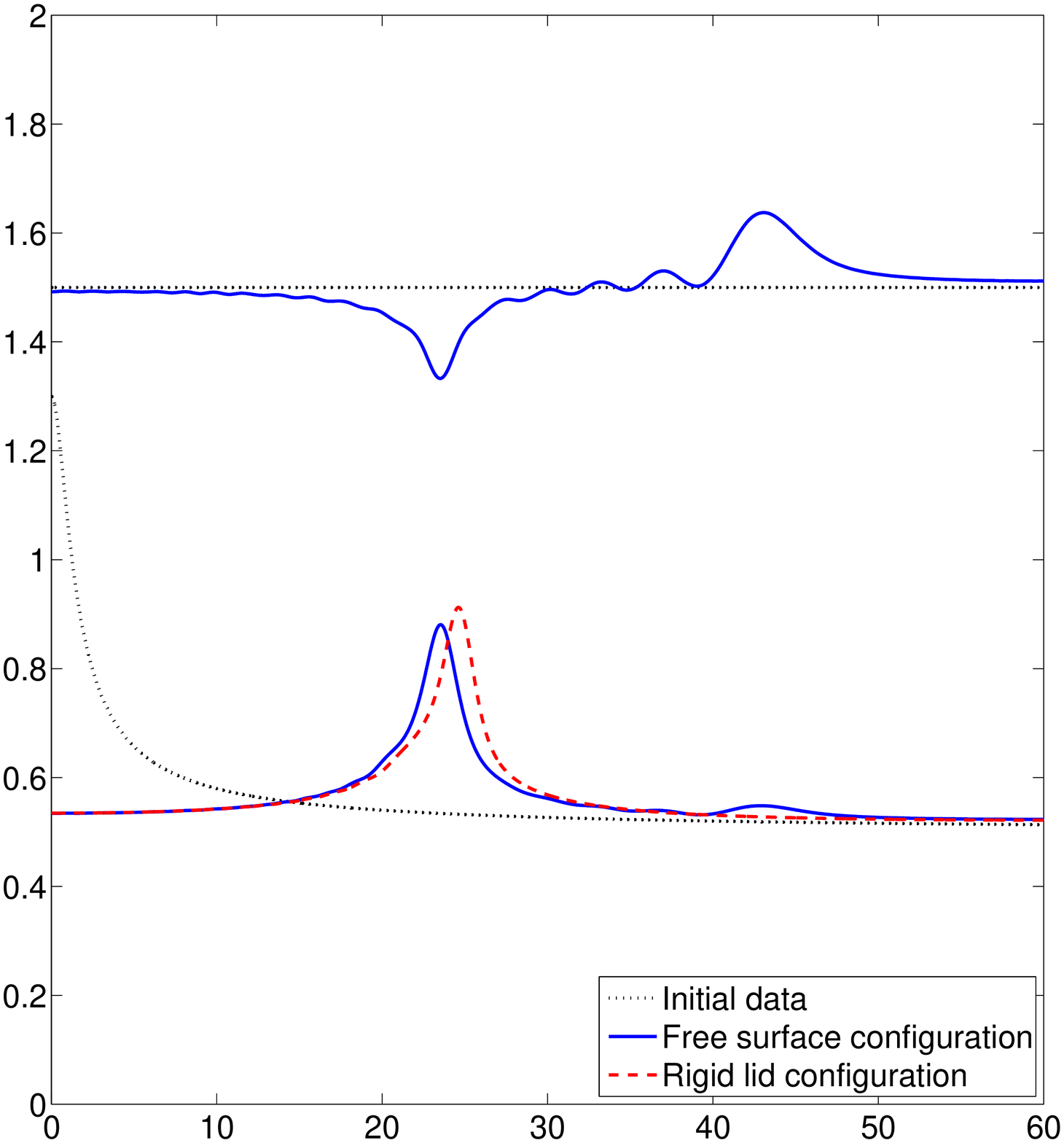}} 
 \caption{Solution of the KdV approximations, in both rigid lid and free surface configurations, with the upper fluid of greater depth, at times (a) t=20, (b) t=40.}
 \label{fig:rigidlidB}
\end{figure}
 \begin{figure}[hp]
 \centering
\subfigure[$t=20$. $\epsilon=0.1,\delta=1/4,\gamma=1/4$.]{\includegraphics [width=0.45\textwidth,keepaspectratio=true]{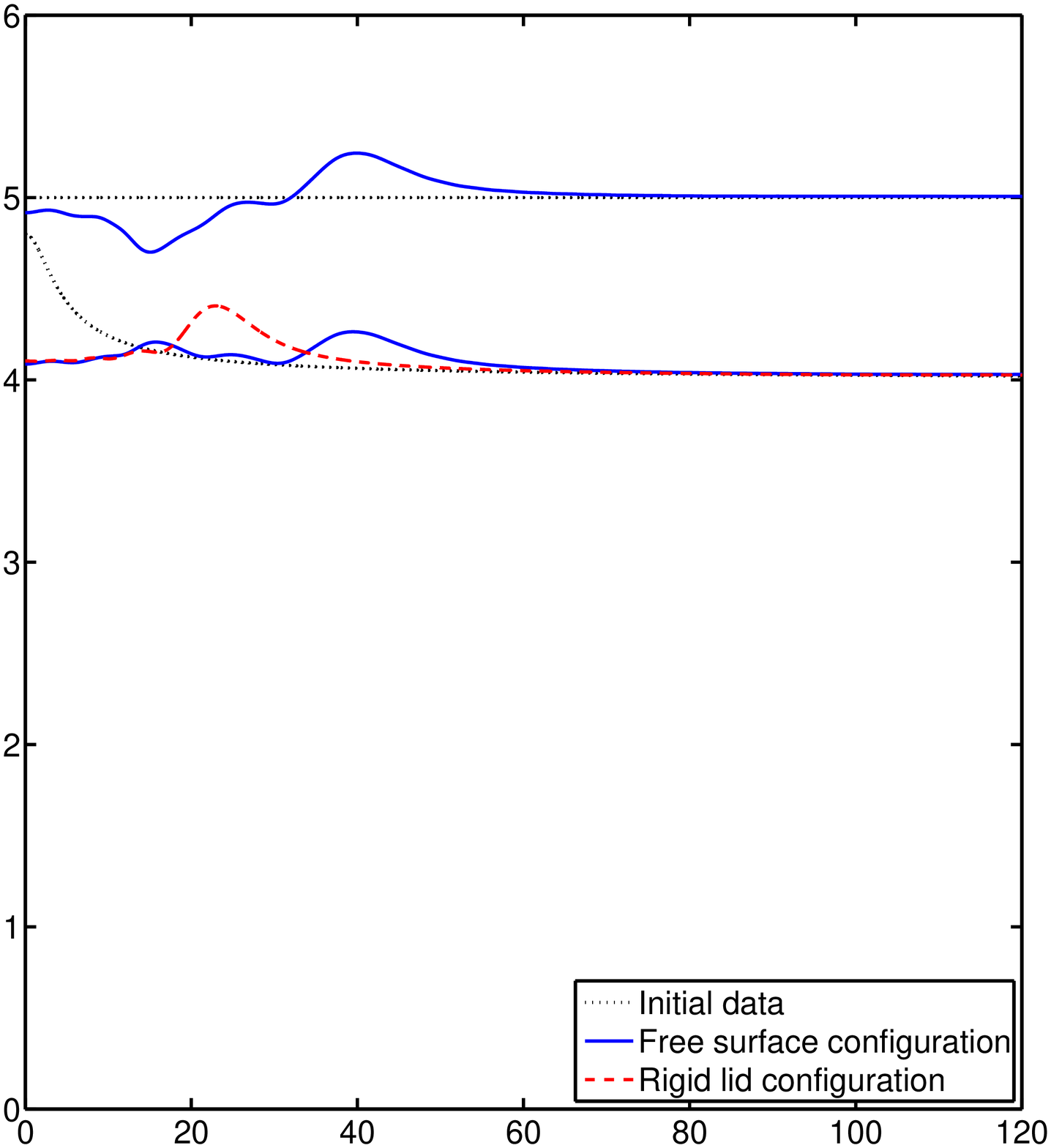}} 
\subfigure[$t=40$. $\epsilon=0.1,\delta=1/4,\gamma=1/4$.]{\includegraphics [width=0.45\textwidth,keepaspectratio=true]{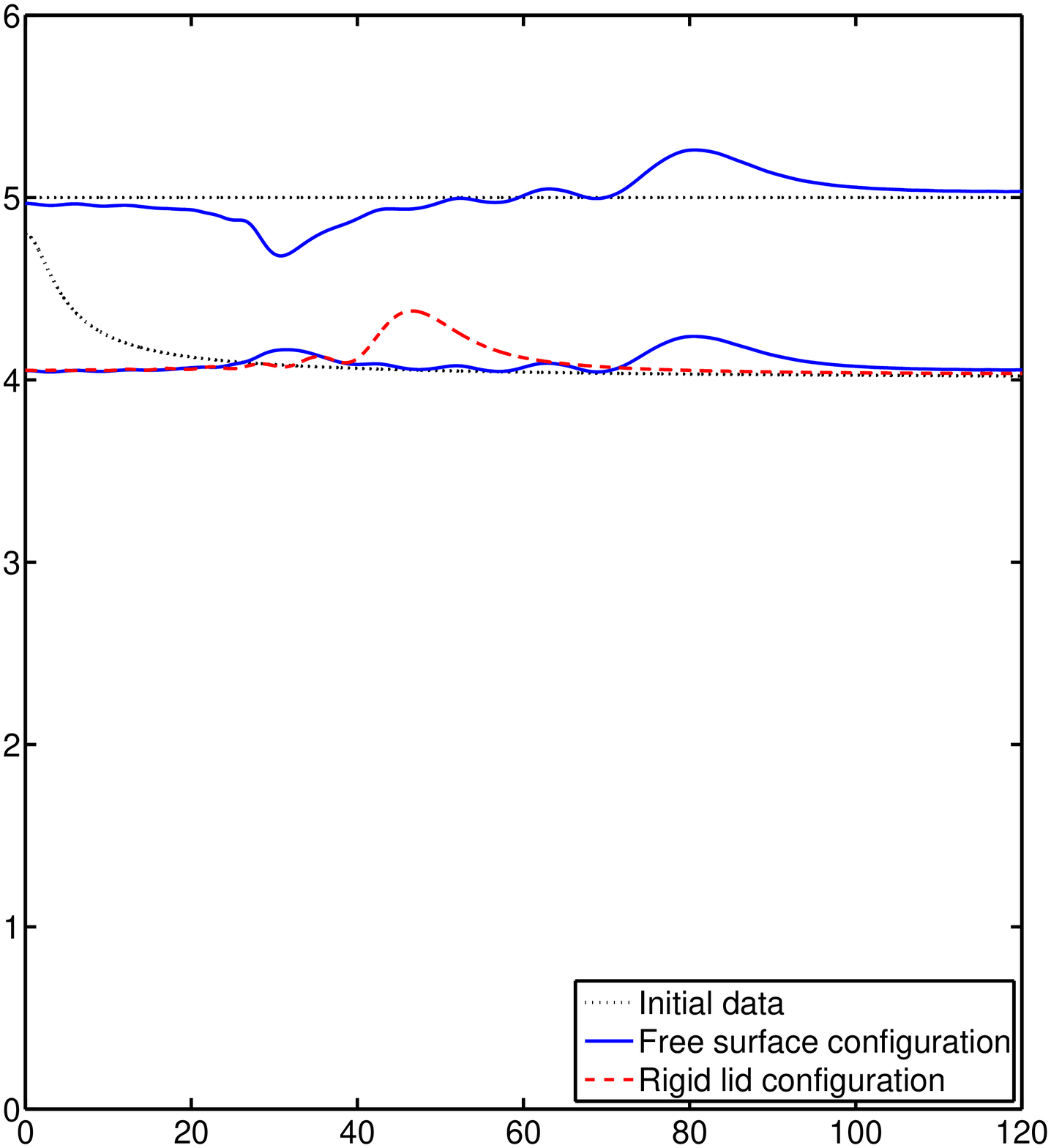}} 
 \caption{Solution of the KdV approximations, in both rigid lid and free surface configurations, with the lower fluid of greater depth, at times (a) t=20, (b) t=40.}
 \label{fig:rigidlidD}
\end{figure}
The Figures~\ref{fig:solitonA} and~\ref{fig:clocheA} provide the snapshots of the previous simulations for $\epsilon=0.1$, and at times $t=20$ and $t=40$. The results of the two schemes are plotted in the same figure, as well as the initial data and the (ten times) emphasized error at the interface, for readability. We only show the right side of the simulations, as the left side is obtained symmetrically.

\medskip

We see that in the case of a solitary wave, the Boussinesq/Boussinesq system produces an almost perfect soliton, that is predicted by the KdV approximation (Figure~\ref{fig:solitonA}). The difference between the two solutions is located in space. However, when the initial data is not sufficiently decreasing at infinity, there is a big qualitative difference between the solution of the Boussinesq model, and the KdV approximation (Figure~\ref{fig:clocheA}), since their difference is not exclusively located in the area of the solitary waves. Therefore, apart from a good initial localization in space, there is a significant interaction between the traveling waves of different wave modes, that is captured by the Boussinesq/Boussinesq system, and not by the KdV approximation.

The same simulations are produced for different parameters, {\it i.e.} $\delta=2$ and $\gamma=1/4$, and are displayed in Figures~\ref{fig:solitonB} and~\ref{fig:clocheB}. We see that the exact same phenomenon appears.

\medskip

Finally, we present from Figures~\ref{fig:rigidlidA} to~\ref{fig:rigidlidD} numerical simulations of the KdV approximation, in both the free surface and rigid lid configurations. The simulation of the Boussinesq/Boussinesq models leads to almost identical results, so that we do not plot them for readability. Again, we set $\epsilon=0.1$, and the different times of the snapshots are $t=20$ and $t=40$. In Figures~\ref{fig:rigidlidA} and~\ref{fig:rigidlidC}, we choose the parameters $\delta=1$, and respectively $\gamma=1/4$ and $\gamma=9/10$. In Figures~\ref{fig:rigidlidB} and~\ref{fig:rigidlidD}, we choose the parameters $\gamma=1/4$, and respectively $\delta=2$ and $\delta=1/4$. Each time, the initial data consists in a bell curve in the interface, with no velocities and with a flat top. One clearly sees that, as discussed in Section~\ref{sSec:AnalysisCoefs}, the rigid lid assumption is satisfactory as an approximation of the free surface, only in the case of small difference of densities ($\gamma \sim 1$), or when the depth ratio becomes large. 

\appendix
\section{Proof of Proposition~\ref{Prop:WPSBOUSS}}
\label{Sec:ProofPropWPSBOUSS}
The proof is made of three steps. First, we introduce an energy of the system, and obtain an a priori estimate on this energy. Then, using this estimate and regularization operators, we prove the existence of a solution of our problem. Finally, using the energy estimate on the difference of two solutions, we get the uniqueness of the solution.

\subsection{Energy estimate}
\label{ProofEstimate}
In the following, we denote by $U$ a solution of~\eqref{SBOUSS} on $[0,T]$, with ${U\in L^\infty([0,T];H^{s+1})^4}$. When we multiply~\eqref{SBOUSS} by $\Lambda^{s}$, the system becomes
\begin{equation} \label{SsBOUSS}\begin{array}{r}\Big(S_0+\epsilon S_1(U)-\epsilon S_2\partial_x^2\Big)\partial_t \Lambda^s U+\epsilon[\Lambda^s, S_1(U)]\partial_t U+\Big(\Sigma_0+\epsilon \Sigma_1(U)-\epsilon \Sigma_2 \partial_x^2\Big)\partial_x \Lambda^s U\\
+\epsilon [\Lambda^s, \Sigma_1(U)] \partial_x U=0.                                 
                                \end{array}\end{equation}
We then introduce the energy associated to the system
\begin{dfntn}
\label{DefEnergy}
Let $U\in H^{s+1}(\RR^4)$. We define the energy of the function $U$ associated to the system~\eqref{SBOUSS} as
 \[E_s(U)\equiv 1/2(S_0 \Lambda^s U,\Lambda^s U)+\epsilon/2(S_1(U) \Lambda^s U,\Lambda^s U)+\epsilon/2(S_2 \Lambda^s \partial_x U,\Lambda^s \partial_x U).\]
 If there is no risk of confusion, we simply write $E_s$.
\end{dfntn}

From Assumption~\ref{Hyp}, $S_0$ and $S_2$ are definite positive (with eigenvalues depending on $\gamma$ and $\delta$). Moreover, since $s>1/2$, one has by Sobolev embeddings ${\big|S_1(U)\big|_{L^\infty}\leq C_0\big\vert U\big\vert_{H^s}}$. Hence, there exists $\alpha=C_0(\frac1{\gamma(1-\gamma)},\delta+\frac1\delta)$ such that if $\epsilon \big\vert U \big\vert_{H^s} < \frac1{C_0(\alpha)}$, then
\begin{equation}\label{EsHseps}
 \frac{1}{\alpha}\big\vert U\big\vert_{H^{s+1}_\epsilon}^2\leq  E_s(U) \leq \alpha \big\vert U\big\vert_{H^{s+1}_\epsilon}^2.
\end{equation}

\medskip

Let us multiply~\eqref{SsBOUSS} on the right by $\Lambda^s U$ and integrate. One obtains
\begin{equation}\label{dEs/dt}\begin{array}{r}
 \displaystyle\frac{d}{dt}E_s= \epsilon/2 (S_1(\partial_t U)\Lambda^s U,\Lambda^s U)-\epsilon ([\Lambda^s,S_1(U)]\partial_t U,\Lambda^s U)+\epsilon/2 ((\Sigma_1(\partial_x U) \Lambda^s U),\Lambda^s U)\\-\epsilon ([\Lambda^s,\Sigma_1(U)]\partial_x U,\Lambda^s U).
\end{array}\end{equation}
Now, one has thanks to Cauchy-Schwarz inequality and Sobolev embeddings
\begin{equation}\label{esti1}|(\Sigma_1(\partial_x U) \Lambda^s U,\Lambda^s U)|\ \leq\ \big| \Sigma_1(\partial_x U)\big|_{L^\infty}\big|\Lambda^s U\big|_{L^2}^2\ \leq\ C_0  \big|U\big|_{H^s}^3. \end{equation}
We then use classical Kato-Ponce commutator estimate~\cite{KatoPonce88}: 
\begin{lmm} \emph{(Kato-Ponce)}\label{Kato-Ponce}
 For $s\geq 0$, if $f\in H^s$ and $g\in H^{s-1}$, then one has the estimate
 \[\big|[\Lambda^s,f]g\big|_{L^2} \leq C_0  \big| f\big|_{H^s}\big| g\big|_{L^\infty}+C_0 \big| \partial_x f\big|_{L^\infty}\big| g\big|_{H^{s-1}}.\]
\end{lmm}
In our case, it leads to the commutator estimate
\begin{equation}\label{esti2}\big|([\Lambda^s,\Sigma_1(U)]\partial_x U,\Lambda^s U)\big|\leq \big|([\Lambda^s,\Sigma_1(U)]\partial_x U\big|_{2}\big|\Lambda^s U\big|_{2}\leq C_0  \big|U\big|_{H^s}^3.\end{equation}

\medskip

In order to deal with the nonlinear terms with a time-derivative $\partial_t$, we will need the following Lemma, using the elliptic form of the operator $P_\epsilon$ defined in~\eqref{OpBOUSS}:
\begin{lmm} \label{LemmaOp}
Let $P_\epsilon(U,\partial)$ and $Q_\epsilon(U,\partial)$ the differential operators defined in~\eqref{OpBOUSS} and $s>1/2$. Then there exists $C_0=C_0(\frac1{\gamma(1-\gamma)},\delta+\frac1\delta)$ such that if $\epsilon \big\vert U \big\vert_{H^{s}} <1/C_0$, one has the following properties:
\begin{enumerate}
\item $P_\epsilon(U,\partial):H^1\to H^{-1}$ is one-to-one, and onto.
\item For $k\in[0,s]$, if one has $P_\epsilon(U,\partial)V\in H^{k-1}$, then $V\in H^{k+1}$. 
\item For $1/2<k\leq s$ and $V\in H^{k+1}$, one has \[1/C_0\big\vert V\big\vert_{H^{k+1}_\epsilon}\leq ( P_\epsilon(U,\partial)V,V)\leq C_0 \big\vert V\big\vert_{H^{k+1}_\epsilon}.\]
\item For $1/2<k\leq s$, the operator $P_\epsilon(U,\partial_x)^{-1}Q_\epsilon(U,\partial_x)$ is bounded $H^{k+1}_\epsilon \rightarrow H^{k+1}_\epsilon$, uniformly with respect to $\epsilon$.
\end{enumerate}
\end{lmm}
\begin{proof}
\begin{enumerate}
\item Using the fact that $S_0$ and $S_2$ are symmetric definite positive (with eigenvalues depending on $\gamma$ and $\delta$) in the formulation of $P_\epsilon$, it is obvious to see that one can choose $C_0=C_0(\frac1{\gamma(1-\gamma)},\delta+\frac1\delta)>0$ such that for $\epsilon\big|U\big|_{L^\infty}<1/C_0$ and any $V\in H^1$,
\[(P_\epsilon(U,\partial)V,V)=(S_0 V,V)+\epsilon(S_1(U)V,V)+\epsilon (S_2\partial_x V,\partial_x V) \geq \frac{\epsilon}{C_0} \big\vert V\big\vert_{H^1}^2.\]
In the same way, the bilinear form $a(V,W)=(P_\epsilon(U,\partial)V,W)$ is clearly continuous on $H^1\times H^1$:
\[(P_\epsilon(U,\partial)V,V)=(S_0 V,V)+\epsilon(S_1(U)V,V)+\epsilon (S_2\partial_x V,\partial_x V) \leq C_0 \big\vert V\big\vert_{H^1}^2.\]
 Using Lax-Milgram lemma, we obtain that for all $F\in H^{-1}$, there exists a unique $V\in H^1$ such that for all $W\in H^1$, $a(V,W)=(F,W)$, and hence there exists a unique variational solution of 
\[P_\epsilon(U,\partial)V=F.\] 

\item We will prove the second point by induction. The result is known for $k=0$ from the previous point. 
Then, we remark that if $P_\epsilon(U,\partial)V=W$ with $W\in L^2$ (and $V\in H^1$), then
\[\epsilon S_2\partial_x^2 V=S_0 V+\epsilon S_1(U)V-W.\]
Therefore, since $S_2$ is invertible, we get $\partial_x^2 V \in L^2$, and thus $V\in H^2$. The second point is therefore proved for $k=1$, and the intermediate values $k\in[0,1]$ follow by interpolation.

Let now assume that the result is known for a $k\in\mathbb{N}^*$ and let $P_\epsilon(U,\partial)V \in H^k\subset H^{k-1}$. From the induction hypothesis, one has $V\in H^{k+1}$, and moreover
\[\epsilon S_2\partial_x^2 \Lambda^k V=S_0 \Lambda^k V+\epsilon \Lambda^k (S_1(U) V)-\Lambda^k P_\epsilon(U,\partial)V.\]
We know that $\Lambda^k V\in L^2$ and $\Lambda^k P_\epsilon(U,\partial)V \in L^2$. Moreover, since $U\in H^{s}$ and $s\geq k>1/2$, we know from the Sobolev embedding that  $\Lambda^k (S_1(U)V) \in L^2$ and $\big\vert \Lambda^k (S_1(U)V) \big\vert_2 \leq C_0 \big\vert U \big\vert_{H^k}\big\vert V \big\vert_{H^{k}}$.
Therefore, $\partial_x^2 \Lambda^k V \in L^2$, and as we know from the induction hypothesis that $V\in H^{k+1}$, one has indeed $V\in H^{k+2}$. The second point is thus proved for all $k\in \mathbb{N}\cap[0,s]$, and we obtain the result for $0\leq k\leq s$ in the same way, starting the induction from $k-\lfloor k \rfloor \in [0,1)$.\par~\par

\item Since $S_0$ and $S_2$ are symmetric definite positive, one has for ${\epsilon\big|U\big|_{H^s}<1/C_0}$ with $C_0=C_0(\frac1{\gamma(1-\gamma)},\delta+\frac1\delta)$, and for any $V\in H^{k+1}$,
\begin{align*}(\Lambda^k (P_\epsilon(U,\partial)V),\Lambda^kV)\ =\ &(S_0 \Lambda^kV,\Lambda^kV)+\epsilon(\Lambda^k (S_1(U)V),\Lambda^kV) +\epsilon (S_2\partial_x \Lambda^kV,\partial_x \Lambda^kV) \\
 \ \geq \ & \frac1{C_0}\big\vert V\big\vert_{H^k}^2+\frac1{C_0}\epsilon\big\vert V\big\vert_{H^{k+1}}^2.
\end{align*}
The second inequality is straightforward.\par~\par
 
\item Let $W\in H^{k+1}$. It is obvious that, since $s\geq k>1/2$, one has \[Q_\epsilon(U,\partial_x )W=\big(\Sigma_0+\epsilon \Sigma_1(U)\big) W-\epsilon \Sigma_2 \partial_x^2 W\in H^{k-1}.\] Hence, we know from the previous points that $V\equiv P_\epsilon(U,\partial_x )^{-1}Q_\epsilon(U,\partial_x )W \in H^{k+1}$, and  
\[ \frac1{C_0}\big\vert V\big\vert_{H^{k+1}_\epsilon}^2\leq(\Lambda^k (P_\epsilon(U,\partial)V),\Lambda^kV)\leq C_0\big\vert V\big\vert_{H^{k+1}_\epsilon}^2.\]
Therefore, one has:
\begin{align*}
 \big\vert V\big\vert_{H^{k+1}_\epsilon}^2&\leq C_0(\Lambda^k (Q_\epsilon(U,\partial)W),\Lambda^k V)\\
 &\leq C_0(\Lambda^k (\Sigma_0+\epsilon \Sigma_1(U))W),\Lambda^kV)+C_0(\partial_x\Lambda^k \Sigma_2 W),\partial_x\Lambda^k V)\\
&\leq C_0(\gamma,\frac1\delta,\big\vert U\big\vert_{H^{k}})\big(\big\vert W\big\vert_{H^{k}}\big\vert V\big\vert_{H^{k}} +\epsilon \big\vert W\big\vert_{H^{k+1}}\big\vert V\big\vert_{H^{k+1}}\big).
\end{align*}
Finally, one has
\begin{align*}\big\vert V\big\vert_{H^{k+1}_\epsilon}^2 \leq & \frac12\big(C_0\big\vert W\big\vert_{H^{k}}-\big\vert V\big\vert_{H^{k}}\big)^2 +\frac\epsilon2 \big(C_0\big\vert W\big\vert_{H^{k+1}}-\big\vert V\big\vert_{H^{k+1}}\big)^2 +\frac12\big\vert V\big\vert_{H^{k+1}_\epsilon}^2+ \frac12C_0^2 \big\vert W\big\vert_{H^{k+1}_\epsilon}^2\\
\leq & \frac12\big\vert V\big\vert_{H^{k+1}_\epsilon}^2+ \frac12C_0^2 \big\vert W\big\vert_{H^{k+1}_\epsilon}^2.\end{align*}
The operator $P_\epsilon(U,\partial_x)^{-1}Q_\epsilon(U,\partial_x)$ is thus bounded $H^{k+1}_\epsilon \rightarrow H^{k+1}_\epsilon$ by $C_0^2$, which ends the proof of the Lemma.
\end{enumerate}\end{proof}

Using this Lemma, we immediately obtain that for $U\in H^{s+1}$ satisfying~\eqref{OpBOUSS} and the hypothesis of the Lemma, one has $\partial_t U \in H^{s}$, and 
\begin{equation}\label{dtvsdx}
 \big|\partial_t U\big|_{H^{s}_\epsilon}=\big|P_\epsilon(U,\partial_x)^{-1}Q_\epsilon(U,\partial_x) \partial_x U\big|_{H^{s}_\epsilon}\leq C_0  \big|\partial_x U\big|_{H^{s}_\epsilon}\leq C_0 \big| U\big|_{H^{s+1}_\epsilon}.
\end{equation}

Therefore, we can use the same calculations as in~\eqref{esti1} and~\eqref{esti2}, and obtain
\begin{equation}\label{esti3}|(S_1(\partial_t U) \Lambda^s U,\Lambda^s U)|\leq \big| S_1(\partial_t U)\big|_{L^\infty}\big|\Lambda^s U\big|_{L^2}^2\leq C_0 \big|U\big|_{H^s}^2\big|U\big|_{H^{s+1}_\epsilon} \end{equation}
(since $s-1>1/2$ such that $\big|\partial_t U\big|_{L^\infty} \leq C_0 \big|\partial_t U\big|_{H^{s-1}}\leq C_0 \big| U\big|_{H^{s+1}_\epsilon}$), and with Kato-Ponce theorem,
\begin{align}\label{esti4}|([\Lambda^s,S_1(U)]\partial_t U,\Lambda^s U)|&\leq \big|([\Lambda^s,S_1(U)]\partial_t U\big|_{L^2}\big|\Lambda^s U\big|_{L^2}\nonumber \\
 &\leq C_0  \big|U\big|_{H^s}^2 \big|\partial_t U\big|_{H^{s-1}} \leq C_0  \big|U\big|_{H^s}^2\big|U\big|_{H^{s+1}_\epsilon}.\end{align}

Finally, one deduces from~\eqref{EsHseps},~\eqref{dEs/dt} and the estimates~\eqref{esti1}-\eqref{esti4}:
\[\frac{d}{dt}E_s \leq \epsilon C_0  \ \big|U\big|_{H^s}^2\big|U\big|_{H^{s+1}_\epsilon}\leq \epsilon C_0 E_s^{3/2},\] 
with $C_0=C_0(\frac1{\gamma(1-\gamma)},\delta+\frac1\delta)$, and providing the fact that $\epsilon\big|U\big|_{L^\infty H^s}<1/C_0$.

From Gronwall-Bihari's inequality and~\eqref{EsHseps}, it follows 
\begin{equation}
\big| U\big|_{H^{s+1}_\epsilon}\leq C_0 E_s^{1/2} \leq C_0\frac{{E_s^{1/2}}\id{t=0}}{1-C_0 \epsilon t {E_s^{1/2}}\id{t=0} } \leq C_0 \frac{\big| U^0\big|_{H^{s+1}_\epsilon}}{1-C_0\epsilon \big| U^0\big|_{H^{s+1}_\epsilon} t}.
\end{equation}
Thus, one sees that there exists $C_0$ such that if $\epsilon\big|U_0\big|_{H^{s+1}_\epsilon}<1/C_0$, then one can choose $T=T(C_0)$ such that the smallness assumption $\epsilon\big|U\big|_{H^s}<1/C_0$ remains valid for any $t\in[0,T/\epsilon]$. Hence, in the frame of the proposition, one has the following estimate:
\begin{equation}\label{est}
\big| U\big|_{L^\infty([0,T/\epsilon];H^{s+1}_\epsilon)} \leq C_0  \frac{\big| U^0\big|_{H^{s+1}_\epsilon}}{1-C_0 \epsilon \big| U^0\big|_{H^{s+1}_\epsilon} t},
\end{equation}
with $C_0=C_0(\frac1{\gamma(1-\gamma)},\delta+\frac1\delta)$ and $T>0$, independent of $\epsilon$.

Let us note that by Lemma~\ref{LemmaOp}, we know that we can also control the time derivative $\big| \partial_t U\big|_{L^\infty([0,T/\epsilon] ; H^{s}_\epsilon)}$ with the same bound. 

\subsection{Existence of a solution}
\label{ProofExistence}
 We can deduce from the energy estimate~\eqref{est} the existence of a maximal solution $U\in L^{1,\infty}_t([0,T/\epsilon] ; H^{s+1})$ of~\eqref{SBOUSS} for any initial data $U^0\in H^{s+1}$. 
 We follow the classical Friedrichs proof, using the regularization operators defined thanks to the Fourier transform as below:
 \begin{equation}\label{DefJnu}
  \forall v\in L^2, \quad \forall \xi\in\RR, \quad \widehat{J_\nu v(\xi)}\equiv\varphi(\nu\xi) \widehat v(\xi),
 \end{equation}
with $\varphi$ a smooth numeric function with compact support, such that $\varphi(0)=1$. These operators have the following classical properties
\begin{lmm} \label{Lem:Jnu}
 \begin{enumerate}
 \item $J_\nu$ is bounded $H^s\to H^s$: there exists $C_0(s,\nu)$ such that for any $v\in H^s$,
\begin{equation}
  \big\vert J_\nu v \big\vert_{H^s}\leq C_0(s,\nu) \big\vert  v \big\vert_{H^s}
 \end{equation}
 \item $J_\nu$ commutes with $\Lambda^s$, and is a self-adjoint operator
 \item There exists $C_0$ independent of $\nu$ such that 
\begin{equation}
  \big\vert J_\nu v \big\vert_{2}\leq C_0 \big\vert  v \big\vert_{2}
 \end{equation}
 \end{enumerate}
\end{lmm}
We obtain then a solution of~\eqref{OpBOUSS} as the limit of $U_\nu$ the solutions of
\begin{equation}
 \label{OpBOUSSnu}
 \partial_t U_\nu +J_\nu P_\epsilon(J_\nu U_\nu,\partial_x)^{-1} Q_\epsilon(J_\nu U_\nu,\partial_x) J_\nu \partial_x U_\nu =0,
\end{equation}
with ${U_\nu}\id{t=0}=U^0$.

\bigskip

From Lemma~\ref{Lem:Jnu}, one has that~\eqref{OpBOUSSnu} is an ordinary differential equation on the Banach space $H^s$. Thus, thanks to Cauchy-Lipschitz Theorem, we know that there exists a unique maximal solution $U_\nu \in C([0,T_\nu),H^s)$.

Using the calculations of Section~\ref{ProofEstimate} and the properties of $J_\nu$, one obtains the energy estimate
\[\big| U_\nu\big|_{H^{s+1}_\epsilon}\leq C_0  \frac{\big| U^0\big|_{H^{s+1}_\epsilon}}{1-C_0 \epsilon \big| U^0\big|_{H^{s+1}_\epsilon} t},\]
with $C_0$ independent of $\nu$. Moreover, $U_\nu$ satisfies~\eqref{dtvsdx}:
\begin{equation}
 \big|\partial_t U_\nu\big|_{H^{s-1}}\leq C_0  \big|\partial_x U\big|_{H^{s-1}} \leq C_0 \big| U\big|_{H^{s}}.
\end{equation}
Thus one can find $T>0$ independent of $\nu$ such that the solution $U_\nu$ does not blow up for $t\in [0,T/\epsilon]$. In particular, one has $T_\nu>T/\epsilon>0$.

\medskip

Now, since $(U_\nu)_{\nu>0}$ is uniformly bounded in $L^{1,\infty}_t([0,T/\epsilon],H^{s+1}_\epsilon)$, one can extract a subsequence that weakly converges towards a function $U\in L^{1,\infty}_t([0,T/\epsilon],H^{s+1}_\epsilon)$. We want to use Ascoli theorem, but we need the injection $H^s(\RR)\subset H^{s-1}(\RR)$ to be compact, which is not true since $\RR$ is unbounded. However, one can easily circumvent this problem, using weighted Sobolev spaces for example. Finally, one obtains that there exists $U\in C^0([0,T/\epsilon),H^{s-1})$ such that $U_\nu$ converges strongly towards $U$ as the subsequence $\nu$ tends to 0. Also, by interpolation inequalities, $U_\nu$ converges strongly towards $U$ in $C^0([0,T/\epsilon),H^{s-\alpha})$, for $0<\alpha\leq 1$. Then, since one can find $\alpha$ such that $H^{s-\alpha}$ injects continuously in $C^1(\RR)$, one proves that $Q_\epsilon(J_\nu U_\nu,\partial_x) J_\nu \partial_x U_\nu$ converges to  $Q_\epsilon( U,\partial_x)\partial_x U$ and $J_\nu P_\epsilon(J_\nu U_\nu,\partial_x)\partial_t U$ converges to $P_\epsilon(U_\nu,\partial_x)\partial_t U$ as $\nu\to 0$. Hence, $U$ is indeed a solution of~\eqref{OpBOUSS}. 

\medskip

From Section~\ref{ProofEstimate}, we know that $U\in L^{1,\infty}_t([0,T/\epsilon] ; H^{s+1}_\epsilon)$, and one can prove (see~\cite{Taylor97} XVI.1.4 for example) that $U\in  C^0([0,T/\epsilon) ; H^{s+1}) \cap C^1([0,T/\epsilon) ; H^{s})$.

\subsection{Uniqueness of the solution}
\label{ProofUnicity}
Let $U_1, U_2 \in  C^0([0,T/\epsilon) ; H^{s+1}) \cap C^1([0,T/\epsilon) ; H^{s})$ be two solutions of the Cauchy problem~\eqref{SBOUSS} with initial data ${{U_1}\id{t=0}={U_2}\id{t=0}=U^0}$. One can immediately check that ${R\equiv U_1-U_2}$ satisfies
\begin{align} \label{RRBOUSS}\Big(S_0+\epsilon S_1(U_1)-\epsilon S_2\partial_x^2\Big)\partial_t \Lambda^s R+\Big(\Sigma_0+\epsilon \Sigma_1(U_1)-\epsilon \Sigma_2 \partial_x^2\Big)\partial_x \Lambda^s R \nonumber \\
+\epsilon[\Lambda^s, S_1(U_1)]\partial_t R + \epsilon [\Lambda^s, \Sigma_1(U_1)] \partial_x R=\epsilon F,\end{align}
with $F=-\Lambda^s \Big(S_1(R)\partial_t U_2+\Sigma_1(R)\partial_x U_2\Big)$. Then, we can carry out the same calculations as in Section~\ref{ProofEstimate} on $R$, and obtain the equivalent energy estimate
\[\frac{d}{dt}E_s(R) \leq \epsilon  C_0 ( \big|U_1 \big|_{H^{s}}+ \big|U_2 \big|_{H^{s}})E_s,\] 
with $C_0=C_0(\frac{1}{\gamma(1-\gamma)},\delta+\frac1\delta,\big| U^0\big|_{H^{s+1}_\epsilon})$.

\medskip

From Gronwall-Bihari's inequality and the estimate~\eqref{est} on $U_1$ and $U_2$, and since $E_s(R)\id{t=0}=0$, one has immediately $E_s(R)=0$ on $[0,T/\epsilon]$, and finally $U_1=U_2$.

\begin{acknowledgement}
\noindent{\bf Acknowledgements.}
This work received the support of the Agence Nationale de la Recherche (project ANR-08-BLAN-0301-01). The author would also like to thank David Lannes                          for very helpful and stimulating discussions, and Florent Chazel for his help on numerical simulations.
\end{acknowledgement}
\bibliographystyle{abbrv}
\bibliography{Biboceano}

\end{document}